\documentstyle[bezier]{article}
\begin{document}

\input amssym.def \input amssym
\font\tencmmib=cmmib10 \font\sevencmmib=cmmib7 \font\fivecmmib=cmmib5
\newfam\cmmibfam \textfont\cmmibfam=\tencmmib
\scriptfont\cmmibfam=\sevencmmib \scriptscriptfont\cmmibfam=\fivecmmib

\font\mff=msbm10 at 15pt \font\mathfontzwoelf=msbm10 at 12pt
\font\mathfontzehn=msbm10 at 10pt \font\mathfontsieben=msbm7
\font\mathfontfuenf=msbm7 at 7pt \font\gothic=eufm10 at 12 pt
\font\script=eusm10 at 12 pt \font\romcur=eurm10 at 12 pt
\font\helv=cmssbx10 at 10 pt
\def\nz{\hfill\break}
\def\isom{\ifmmode\ \cong\ \else$\isom$\fi}

\def\aa{\hbox{\gothic a}} \def\bb{\hbox{\gothic b}} \def\pp{\hbox{\gothic
    p}} \def\dd{\hbox{\gothic d}} \def\Gg{\hbox{\gothic g}}
\def\ff{\hbox{\gothic f}} \def\hh{\hbox{\gothic h}} \def\mm{\hbox{\gothic
    m}} \def\nn{\hbox{\gothic n}} \def\rr{\hbox{\gothic r}}
\def\tt{\hbox{\gothic t}} \def\kk{\hbox{\gothic k}} \def\ll{\hbox{\gothic
    l}} \def\Ss{\hbox{\gothic s}} \def\cc{\hbox{\gothic c}}
\def\uu{\hbox{\gothic u}} \def\vv{\hbox{\gothic v}} \def\ee{\hbox{\gothic
    e}} \def\zz{\hbox{\gothic z}} \def\qq{\hbox{\gothic q}}
\def\oo{\hbox{\gothic o}}

\def\AA{\hbox{\gothic A}} \def\BB{\hbox{\gothic B}} \def\HH{\hbox{\gothic
    H}} \def\DD{\hbox{\gothic D}} \def\MM{\hbox{\gothic M}}
\def\LL{\hbox{\gothic L}} \def\JJ{\hbox{\gothic J}} \def\WW{\hbox{\gothic
    W}} \def\UU{\hbox{\gothic U}} \def\OO{\hbox{\gothic O}}
\def\PP{\hbox{\gothic P}}

\def\VV{\hbox{\gothic V}} \def\YY{\hbox{\gothic Y}} \def\ZZ{\hbox{\gothic
    Z}} \def\CC{\hbox{\gothic C}}

\def\oM{\hbox{\romcur o}_M} \def\oG{\hbox{\romcur o}_G}

\def\fO{{\Bbb O}} \def\fA{{\Bbb A}}

\def\nat{\ifmmode%
         \mathchoice{\hbox{\mathfontzwoelf N}}%
                     {\hbox{\mathfontzehn N}}%
                     {\hbox{\mathfontsieben N}}%
                     {\hbox{\mathfontfuenf N}}%
          \else$\nat$%
          \fi} \def\fN{\nat}
\def\integer{\ifmmode%
          \mathchoice{\hbox{\mathfontzwoelf Z}}%
                     {\hbox{\mathfontzehn Z}}%
                     {\hbox{\mathfontsieben Z}}%
                     {\hbox{\mathfontfuenf Z}}%
          \else$\integer$%
          \fi} \def\fZ{\integer}

\def\fatf{\ifmmode%
          \mathchoice{\hbox{\mathfontzwoelf F}}%
                     {\hbox{\mathfontzehn F}}%
                     {\hbox{\mathfontsieben F}}%
                     {\hbox{\mathfontfuenf F}}%
          \else$\fatf$%
          \fi} \def\fF{\fatf}

\def\komp{\ifmmode%
          \mathchoice{\hbox{\mathfontzwoelf C}}%
                     {\hbox{\mathfontzehn C}}%
                     {\hbox{\mathfontsieben C}}%
                     {\hbox{\mathfontfuenf C}}%
          \else$\komp$%
          \fi} \def\fC{\komp}

\def\fatd{\ifmmode%
     \mathchoice{\hbox{\mathfontzwoelf D}}%
      {\hbox{\mathfontzehn D}}%
       {\hbox{\mathfontsieben D}}%
         {\hbox{\mathfontfuenf D}}%
         \else$\fatd$\/\fi} \def\fD{\fatd}

\def\proj{\ifmmode%
          \mathchoice{\hbox{\mathfontzwoelf P}}%
                     {\hbox{\mathfontzehn P}}%
                     {\hbox{\mathfontsieben P}}%
                     {\hbox{\mathfontfuenf P}}%
                     \else$\proj$\/\fi} \def\fP{{\Bbb P}}
\def\ball{\ifmmode%
          \mathchoice{\hbox{\mathfontzwoelf B}}%
                     {\hbox{\mathfontzehn B}}%
                     {\hbox{\mathfontsieben B}}%
                     {\hbox{\mathfontfuenf B}}%
                     \else$\ball$\/\fi} \def\fB{\ball}
\def\sieg{\ifmmode%
          \mathchoice{\hbox{\mathfontzwoelf S}}%
                     {\hbox{\mathfontzehn S}}%
                     {\hbox{\mathfontsieben S}}%
                     {\hbox{\mathfontfuenf S}}%
                     \else$\sieg$\/\fi} \def\fS{\sieg}
\def\uki{\ifmmode%
          \mathchoice{\hbox{\mathfontzwoelf U}}%
                     {\hbox{\mathfontzehn U}}%
                     {\hbox{\mathfontsieben U}}%
                     {\hbox{\mathfontfuenf U}}%
                     \else$\uki$\/\fi} \def\fU{\uki}
\def\rat{\ifmmode%
          \mathchoice{\hbox{\mathfontzwoelf Q}}%
                     {\hbox{\mathfontzehn Q}}%
                     {\hbox{\mathfontsieben Q}}%
                     {\hbox{\mathfontfuenf Q}}%
                     \else$\rat$\/\fi}

                   \def\fQ{\rat}
\def\real{\ifmmode%
          \mathchoice{\hbox{\mathfontzwoelf R}}%
                     {\hbox{\mathfontzehn R}}%
                     {\hbox{\mathfontsieben R}}%
                     {\hbox{\mathfontfuenf R}}%
                     \else$\real$\/\fi}
\def\fR{\real}%
\def\fatH{\ifmmode%
          \mathchoice{\hbox{\mathfontzwoelf H}}%
                     {\hbox{\mathfontzehn H}}%
                     {\hbox{\mathfontsieben H}}%
                     {\hbox{\mathfontfuenf H}}%
                     \else$\fatH$\/\fi}
\def\fH{\fatH}%

\def\fatG{\ifmmode%
          \mathchoice{\hbox{\mathfontzwoelf G}}%
               {\hbox{\mathfontzehn G}}%
                    {\hbox{\mathfontsieben G}}%
                          {\hbox{\mathfontfuenf G}}%
                          \else$\fatG$\/\fi}
\def\fG{\fatG}%
\def\fX{\hbox{\mathfontzwoelf X}} \def\fY{\hbox{\mathfontzwoelf Y}}

\def\itemm{\itemitem}
\def\itemmm{\par\indent\indent\hangindent3\parindent\textindent}
\def\bs{\backslash}

\def\gD{\Delta} \def\gd{\delta} \def\gG{\ifmmode {\Gamma} \else$\gG$\fi}
\def\gg{\gamma} \def\gs{\sigma} \def\gS{\Sigma} \def\gz{\zeta}
\def\gff{\varphi} \def\gt{\theta} \def\gT{\Theta} \def\ga{\alpha}
\def\gb{\beta} \def\gk{\kappa} \def\gl{\lambda} \def\gL{\Lambda}
\def\gO{\Omega} \def\go{\omega} \def\gm{\mu} \def\gn{\nu} \def\gr{\rho}
\def\grr{\varrho} \def\ge{\varepsilon} \def\gx{\xi} \def\gy{\eta}

\def\bfo{\go} \def\del{\partial} \def\bfa{\ga}

\def\delbar{\overline{\partial}} \def\-#1{\overline{#1}}
\def\~#1{\tilde{#1}} \def\o{\over} \def\ds{\displaystyle} \def\cQ{$\cal Q \
  $} \def\cO{\ifmmode {\cal O} \else$\cO$\fi}

\def\f-#1{\mathchoice{\hbox{\mathfontzwoelf #1}}%
                     {\hbox{\mathfontzehn #1}}%
                     {\hbox{\mathfontsieben #1}}%
                     {\hbox{\mathfontfuenf #1}}}%


                   \def\cA{\ifmmode {\cal A} \else$\cA$\fi}

                   \def\cB{\ifmmode {\cal B} \else$\cB$\fi}

                   \def\cOK{\ifmmode {\cal O}_K \else$\cOK$\fi}

                   \def\cC{\ifmmode {\cal C} \else$\cC$\fi}

                   \def\cD{\ifmmode {\cal D} \else$\cD$\fi}

                   \def\cDh{\ifmmode {\check{\cal D}} \else$\cDh$\fi}

                   \def\cDs{\ifmmode {{\cal D}^*} \else$\cDs$\fi}

                   \def\cE{\ifmmode {\cal E} \else$\cE$\fi}

                   \def\cF{\ifmmode {\cal F} \else$\cF$\fi}

                   \def\cG{\ifmmode {\cal G} \else$\cG$\fi}

                   \def\cH{\ifmmode {\cal H} \else$\cH$\fi}

                   \def\cI{\ifmmode {\cal I} \else$\cI$\fi}

                   \def\cJ{\ifmmode {\cal J} \else$\cJ$\fi}

                   \def\cK{\ifmmode {\cal K} \else$\cK$\fi}

                   \def\cL{\ifmmode {\cal L} \else$\cL$\fi}

                   \def\cM{\ifmmode {\cal M} \else$\cM$\fi}

                   \def\cN{\ifmmode {\cal N} \else$\cN$\fi}

                   \def\cO{\ifmmode {\cal O} \else$\cO$\fi}

                   \def\cP{\ifmmode {\cal P} \else$\cP$\fi}

                   \def\cR{\ifmmode {\cal R} \else$\cR$\fi}

                   \def\cQ{\ifmmode {\cal Q} \else$\cQ$\fi}

                   \def\cS{\ifmmode {\cal S} \else$\cS$\fi}

                   \def\cT{\ifmmode {\cal T} \else$\cT$\fi}

                   \def\cU{\ifmmode {\cal U} \else$\cU$\fi}
                   \def\cW{\ifmmode {\cal W} \else$\cW$\fi}
                   \def\cV{\ifmmode {\cal V} \else$\cV$\fi}

                   \def\cY{\ifmmode {\cal Y} \else$\cY$\fi}

                   \def\cX{\ifmmode {\cal X} \else$\cX$\fi}

                   \def\cZ{\ifmmode {\cal Z} \else$\cZ$\fi}

                   \def\cSC{\ifmmode {\cS\cC(\gG)} \else$\cSC$\fi}


                   \def\scA{\ifmmode \hbox{{\script A}} \else$\scA$\fi}
                   \def\scB{\ifmmode \hbox{{\script B}} \else$\scB$\fi}

                   \def\scC{\ifmmode \hbox{{\script C}} \else$\scC$\fi}

                   \def\scD{\ifmmode \hbox{{\script D}} \else$\scD$\fi}

                   \def\scDh{\ifmmode {scheck{\script D}} \else$\scDh$\fi}

                   \def\scDs{\ifmmode {{\script D}^*} \else$\scDs$\fi}

                   \def\scE{\ifmmode \hbox{{\script E}} \else$\scE$\fi}

                   \def\scF{\ifmmode \hbox{{\script F}} \else$\scF$\fi}

                   \def\scG{\ifmmode \hbox{{\script G}} \else$\scG$\fi}

                   \def\scH{\ifmmode \hbox{{\script H}} \else$\scH$\fi}

                   \def\scI{\ifmmode \hbox{{\script I}} \else$\scI$\fi}

                   \def\scJ{\ifmmode \hbox{{\script J}} \else$\scJ$\fi}

                   \def\scK{\ifmmode \hbox{{\script K}} \else$\scK$\fi}

                   \def\scL{\ifmmode \hbox{{\script L}} \else$\scL$\fi}

                   \def\scM{\ifmmode \hbox{{\script M}} \else$\scM$\fi}

                   \def\scN{\ifmmode \hbox{{\script N}} \else$\scN$\fi}

                   \def\scO{\ifmmode \hbox{{\script O}} \else$\scO$\fi}

                   \def\scP{\ifmmode \hbox{{\script P}} \else$\scP$\fi}

                   \def\scQ{\ifmmode \hbox{{\script Q}} \else$\scQ$\fi}

                   \def\scS{\ifmmode \hbox{{\script S}} \else$\scS$\fi}

                   \def\scT{\ifmmode \hbox{{\script T}} \else$\scT$\fi}

                   \def\scU{\ifmmode \hbox{{\script U}} \else$\scU$\fi}

                   \def\scV{\ifmmode \hbox{{\script V}} \else$\scV$\fi}

                   \def\scY{\ifmmode \hbox{{\script Y}} \else$\scY$\fi}

                   \def\scZ{\ifmmode \hbox{{\script Z}} \else$\scZ$\fi}

                   \def\rB{\ifmmode {\cyr B} \else$\rB$\fi}

                   \def\rC{\ifmmode {\cyr C} \else$\rC$\fi}

                   \def\rD{\ifmmode {\cyr D} \else$\rD$\fi}

                   \def\rI{\ifmmode {\cyr I} \else$\rI$\fi}

                   \def\rL{\cyr L}

                   \def\rP{\ifmmode {\cyr P} \else$\rP$\fi}

                   \def\rQ{\ifmmode {\cyr Q} \else$\rQ$\fi}

                   \def\rU{\ifmmode {\cyr U} \else$\rU$\fi}

                   \def\rV{\ifmmode {\cyr V} \else$\rV$\fi}

                   \def\rZ{\ifmmode {\cyr Z} \else$\rZ$\fi}
                   \def\Ad{\hbox{Ad}} \def\ad{\hbox{ad}}
                   \def\dim{\hbox{dim}} \def\deg{\hbox{deg}}
                   \def\Pic{\hbox{Pic}} \def\Jac{\hbox{Jac}}
                   \def\Aut{\hbox{Aut}} \def\Im{\hbox{Im}}
                   \def\mod{\hbox{mod}} \def\diag{\hbox{diag}}
                   \def\det{\hbox{det}} \def\Ker{\hbox{Ker}}
                   \def\Hess{\hbox{Hess}} \def\rank{\hbox{rank}}

                   \def\Sym{\hbox{Sym}} \def\im{\hbox{Im}}

                   \def\isom{\ifmmode\ \cong\ \else$\isom$\fi}
                   \def\storth{\underline{\perp}} \def\TT{ {^{^(} T^{^)} }}

                   \def\STB{$\cal ST(B)$} \def\elin{$\ge-$line}

                   \def\nni{\supset} \def\vvs{${\gG \bs \cal D}^*$}
                   \def\ovv{\ifmmode \overline{\gG \bs \cal D}\
                     \else$\ovv$\fi} \def\ovvone{$\overline{\gG_1 \bs {\cal
                         D}_1}$} \def\ovvtwo{$\overline{\gG_2 \bs {\cal
                         D}_2}$}

                   \def\cDone{${\cal D}_1$} \def\cDtwo{${\cal D}_2$}

                   \def\xg{\ifmmode {X_{\gG}} \else$\xg$\fi}
                   \def\xgeq{\ifmmode {\xg = \gG\bs\cD} \else$\xgeq$\fi}

                   \def\xgs{\ifmmode {X_{\gG}^*} \else$\xgs$\fi}

                   \def\xgc{\ifmmode {\overline{X}_{\gG}} \else$\xgc$\fi}

                   \def\dcom{\ifmmode {\overline{D}} \else$\dcom$\fi}

                   \def\xs{\ifmmode {X^*} \else$\xs$\fi}

                   \def\Xbar{\ifmmode {\overline{X}} \else$\Xbar$\fi}

                   \def\Ybar{\ifmmode {\overline{Y}} \else$\Ybar$\fi}

                   \def\Dbar{\ifmmode {\overline{D}} \else$\Dbar$\fi}

                   \def\Feo{F_{\gl}^{\ge,0}} \def\Veo{V_{\gl}^{\ge,0}}

                   \def\piD{\pi_1(D_i)}

                   \def\pigD{\pi_1((\gD_{F_{\gl}})_{\ge}\bs\gD_{F_{\gl}})}

                   \def\igD{(\gD_{F_{\gl}})_{\ge}\bs\gD_{F_{\gl}}}

                   \def\dedelt{D\cup (\gD_{\ge}\bs\gD)}

                   \def\Delep{(\gD_{F_{\gl}})_{\ge}\bs\gD_{F_{\gl}}}

                   \def\Dibar{\ifmmode {\overline{D_i}} \else$\Dibar$\fi}

                   \def\Diemd{\gD_{i,\ge}\bs \gD_i}
                   \def\Djemd{\gD_{j,\ge}\bs \gD_j}
                   \def\Piemp{\Pi_{i,\ge}\bs \Pi_i}
                   \def\Pjemp{\Pi_{j,\ge}\bs \Pi_j} \def\Demd{\gD_{\ge}\bs
                     \gD} \def\DDem{D\cup \Demd} \def\Pemp{\Pi_{\ge}\bs
                     \Pi}

                   \def\Pieo{\Pi_{\gl}^{\ge,0}}
                   \def\Pieoeq{\Pi_{\gl}^{\ge,0}=(\Pi_{\gl})_{\ge}\bs\Pi_{\gl}
                     }

                   \def\Dieo{\gD_{\gl}^{\ge,0}}
                   \def\Dieoeq{\gD_{\gl}^{\ge,0}=(\gD_{\gl})_{\ge}\bs\gD_{\gl}
                     }

                   \def\key{\pi_1(\xg)=\pi_1((\gD_{\gG})_{\ge}\bs\gD_{\gG}\cup
                     D)} \def\inneq{\varsubsetneqq} \def\GQ{\ifmmode
                     {G_{\fQ}} \else $\GQ$ \fi} \def\MQ{\ifmmode {M_{\fQ}}
                     \else $\MQ$ \fi} \def\cTG{\ifmmode {{\cal T}(G_{\fQ})}
                     \else$\cTG$\fi}

                   \def\cTg{\ifmmode {{\cal T}(\gG)} \else$\cTg$\fi}

                   \def\cSG{\ifmmode {{\cal S}(G_{\fQ})} \else$\cSG$\fi}

                   \def\cSg{\ifmmode {{\cal S}(\gG)} \else$\cSg$\fi}

                   \def\cTSG{\ifmmode {{\cal TS}(G_{\fQ})} \else$\cTSG$\fi}

                   \def\cTSg{\ifmmode {{\cal TS}(\gG)} \else$\cTSg$\fi}

                   \def\prodbuild{$\cT(G_{\fQ,1}) \otimes \cT(G_{\fQ,2})$}
                   \def \PB{\ifmmode {\prodbuild} \else $\prodbuild$ \fi}

                   \def\balln{\ifmmode {\gG_K(n) \bs \fB_3}
                     \else$\balln$\fi}
\def\cIf{\ifmmode {\cal I}_5%
  \else$\cIf$\/\fi} \def\cIt{\ifmmode {\cal I}_{10} \else$\cIt$\fi}
\def\JG{\ifmmode {G_{25,920}} \else$\JG$\fi}

\def\HG{G_{648}}

\def\pts{\ifmmode \hbox{\helv P}^2_6 \else $\pts$\fi} \def\pthrees{\ifmmode
  \hbox{\helv P}^3_6 \else $\pthrees$\fi}

\def\pos{\ifmmode \hbox{\helv P}^1_6 \else $\pos$\fi}

\def\pkn{{\helv P}$^k_n$\ } \def\helva{\hbox{\helv A}}

\def\pid{\pi_1(\gD_{\ge}\bs \gD \cup D)} \def\pip{\pi_1(\Pi_{\ge}\bs \Pi
  \cup E)}
\def\GKk{G_{K|k}}

\def\hra{\hookrightarrow} \def\hla{\hookleftarrow}
\def\lra{\longrightarrow} \def\sura{\twoheadrightarrow}
\def\lla{\longleftarrow} \def\rar{\rightarrow}
\def\llra{\longleftrightarrow} \def\ra{\rightarrow} \def\Ra{\Rightarrow}
\def\La{\Leftarrow} \def\bs{\ifmmode {\setminus} \else$\bs$\fi}
\def\hat{\widehat} \def\tilde{\widetilde} \def\nin{\not\in}
\def\inn{\subset} \def\nni{\supset} \def\und{\underline }
\def\ove{\overline}

\def\ende{\hfill $\Box$ \vskip0.25cm }

\def\p{\prime} \def\sdprod{\rtimes}



\newtheorem{theorem}{Theorem}[section] \newtheorem{claim}{Claim}[theorem]
\newtheorem{Lemma}{Lemma}[theorem]
\newtheorem{Corollary}{Corollary}[theorem]
\newtheorem{lemma}[theorem]{Lemma}
\newtheorem{proposition}[theorem]{Proposition}
\newtheorem{corollary}[theorem]{Corollary}
\newtheorem{observation}[theorem]{Observation}
\newtheorem{question}[theorem]{Question} \newcounter{def}
\newenvironment{definition}{\refstepcounter{theorem}
  \vskip.2cm\noindent{\bf Definition \thesection.\arabic{theorem}
    }}{\vskip.2cm\noindent} \newenvironment{fact}{\refstepcounter{theorem}
  \vskip.2cm\noindent{\bf Fact \thesection.\arabic{theorem}
    }}{\vskip.2cm\noindent }
\newenvironment{problem}{\refstepcounter{theorem} \vskip.2cm\noindent{\bf
    Problem \thesection.\arabic{theorem} }}{\vskip.2cm\noindent}
\newenvironment{notations}{\refstepcounter{theorem} \vskip.2cm\noindent{\bf
    Notations \thesection.\arabic{theorem} }}{\vskip.2cm\noindent}
\newenvironment{facts}{\refstepcounter{theorem} \vskip.2cm\noindent{\bf
    Facts \arabic{section}.\arabic{theorem} }}{\vskip.2cm\noindent}
\newenvironment{example}{\refstepcounter{theorem} \vskip.2cm\noindent{\bf
    Example \thesection.\arabic{theorem} }}{\vskip.2cm\noindent}
\newenvironment{condition}[1]{ \vskip.2cm\noindent{\bf Condition #1:
    }}{\vskip.2cm} \newenvironment{examples}{\refstepcounter{theorem}
  \vskip.2cm\noindent{\bf Examples \thesection.\arabic{theorem}:
    }}{\vskip.2cm\noindent}
\newenvironment{remark}{\refstepcounter{theorem} \vskip.2cm\noindent{\bf
    Remark \thesection.\arabic{theorem} }}{\vskip.2cm\noindent}
\newenvironment{remarks}{\refstepcounter{theorem} \vskip.2cm\noindent{\bf
    Remarks \thesection.\arabic{theorem} }}{\vskip.2cm\noindent}

\pagestyle{headings} \setcounter{secnumdepth}{4} \setcounter{tocdepth}{1}

\title{A Gem of the modular universe} \author{Bruce Hunt \\ FB Mathematik,
  Universit\"at \\ Postfach 3049 \\ 67653 Kaiserslautern \\ {\small e-mail
    hunt@mathematik.uni-kl.de}} \maketitle

\begin{center}
  {\bf Introduction}
\end{center}

There is a very beautiful set of modular varieties, studied in \cite{J},
which have not only nice interpretations as moduli spaces, but also have
projective embeddings as very special hypersurfaces in $\fP^4$. The first
is the Segre cubic $\cS_3$, the unique cubic hypersurface in $\fP^4$
(threefold) with ten ordinary double points. This is the moduli space
$\pos$ of six marked points in $\fP^1$. The dual variety, the Igusa quartic
$\cI_4$, is the moduli space of six points on a conic in $\fP^2$,
birational to $\cS_3$, but with singular locus consisting of 15 singular
lines which meet in 15 singular points. These two varieties also turn out
to be moduli spaces of Jacobians of certain curves:
\begin{itemize}\item $\cS_3$ is the Satake compactification of the
  arithmetic quotient of the complex three-ball $\fB_3$ by the principal
  congruence subgroup $\gG_{\sqrt{-3}}(\sqrt{-3})$ of $U(3,1;\cO_K)$, where
  $K=\fQ(\sqrt{-3})$. This is the moduli space of Jacobians of trigonal
  (also called Picard) curves of genus 4, with a $\sqrt{-3}$-level
  structure.
\item $\cI_4$ is the Satake compactification of the arithmetic quotient of
  the Siegel space of degree 2, $\fS_2$, by the principal congruence
  subgroup of level 2, $\gG(2)\inn Sp(4,\fZ)$. This is the moduli space of
  Jacobians of hyperelliptic curves of genus two with a marking of the six
  branch points.
\end{itemize}
In \cite{J} the moduli questions were treated in detail, but the projective
embeddings were only mentioned. In this paper we will provide proofs (some
of which are known but difficult to find in the liturature) on the
projective realisations $\cS_3$ and $\cI_4$ of these moduli spaces.  These
give rise to other interesting hypersurfaces: the Hessian varieties of
$\cS_3$ and $\cI_4$, hypersurfaces $\cN_5$ and $\Hess(\cI_4)$ of degrees 5
and 10, respectively, and another degree 10 variety $\cW_{10}$, closely
related to $\Hess(\cI_4)$.  $\cW_{10}$ has the same degree, symmetry group
and singular locus of $\Hess(\cI_4)$, and may actually coincide, a point we
were unable to clarify. The variety $\cW_{10}$ is the image of $\cN_5$
under the birational map of $\fP^4$ into itself given by the Jacobian ideal
of $\cS_3$. $\cN_5$ has been studied recently by Barth and Nieto \cite{BN}.

Another beautiful modular variety is the Burkhardt quartic $\cB_4$, which
can be shown to be the unique quartic hypersurface threefold with 45
ordinary double points \cite{JVS}. The variety $\cB_4$ again has a moduli
interpretation:
\begin{itemize}\item $\cB_4$ is the Satake compactification of the
  arithmetic quotient of the three-ball $\fB_3$ by the principal congruence
  subgroup of level 2, $\gG_{\sqrt{-3}}(2)$ in $U(3,1;\cO_K)$, $K$ as
  above.  This is the moduli space of Jacobians of trigonal (also called
  Picard) curves of genus 4, with a level 2 structure.
\end{itemize}
In this paper we introduce an algebraic fourfold, a quintic hypersurface
$\cI_5$ in $\fP^5$, which is related to all the above mentioned varieties;
this is the gem of the title. But it turns out to be also related to the
moduli space $\pts$ of six marked points in $\fP^2$, hence also to the
moduli space of marked cubic surfaces in $\fP^3$. The moduli space $\pts$
has a projective realisation, which we denote by $\cY$; it is a double
cover of $\fP^4$ branched along the Igusa quartic $\cI_4$, $\cY\lra \fP^4$.
We may also consider the double cover of $\fP^4$ branched along the variety
$\cW_{10}$, which we denote by $\cW\lra \fP^4$. The relations to
the above mentioned varieties (with the exception of $\cB_4$, whose
relation to $\cI_5$ is much more subtle and complex, making it impossible
for us to discuss it here) are as follows.

The Segre cubic is the resolving divisor of 36 triple points of $\cI_5$,
while there are 36 hyperplane sections (dual to the 36 triple points)
isomorphic to $\cN_5$. Similarly, $\cI_4$ and $\cW_{10}$ are related to the
dual variety. But by far the most intriguing is the relation to the Coble
variety $\cY$, which is the main result of this paper. Let $\cY\lra \fP^4$
and $\cW\lra \fP^4$ be the double covers branched over $\cI_4$ and
$\cW_{10}$, respectively; let $\cZ=\cW\times_{\fP^4}\cY$ be the fibre
product. Considering the projection from a triple point $p\in \cI_5$
displays $\cI_5$ blown up in $p$ as a double cover of $\fP^4$, branched
along $\cS_3\cup \cN_5$. Dually this is the double cover branched over
$\cI_4\cup \cW_{10}$. We prove

\begin{theorem}[Corollary \ref{c157.1}] $\cI_5$ sits $\gS_6$-equivariantly
  birationally in the center of the diagram \unitlength1cm
  $$\begin{picture}(2,2)\put(0,0){${\cal W}$}
    \put(.5,0.15){\vector(1,0){1.3}} \put(1.9,0){${\Bbb P}^4$.}
    \put(0,1.6){${\cal Z}$} \put(.5,1.75){\vector(1,0){1.3}}
    \put(1.9,1.6){${\cal Y}$}
    \put(.25,1.5){\vector(0,-1){1}}\put(1.95,1.5){\vector(0,-1){1}}
    \put(.4,1.6){\vector(1,-1){.5}} \put(.9,.9){${\cal Z}'$}
    \put(1.2,.8){\vector(1,-1){.5}}
\end{picture}$$
This shows the relation between the quintic $\cI_5$ and the Coble variety
$\cY$. \end{theorem} So we see that $\cI_5$ is {\em very near} to being the
moduli space of marked cubic surfaces, but not quite. Let us remark here
that one suspects this diagram to have the following property: all
varieties $\cZ, \cW, \cZ', \cY$ are quotients of the symmetric domain of
dimension four and type $\bf IV$ (or equivalently, type $\bf I_{2,2}$).
This is true for $\cY$ (see Theorem \ref{t133.1}). In \cite{mat} a diagram
of subgroups commensurable with the group giving $\cY$ as quotient is given
on p.~389. These groups are in the group $U(2,2)$, which gives the domain
of type $\bf I_{2,2}$ which by one of the exceptional isomorphisms is the
same domain as type $\bf IV_4$.  The following diagram of groups:
$$\begin{array}{rcl} S\gG(i+1) & \hra & S\gG_T(i+1) \\ \downarrow & &
  \downarrow \\ \gG_M(i+1) & \hra & \gG_T(i+1) \end{array}$$ may be
conjectured to correspond to the diagram of varieties above; all arrows are
injective inclusions of subgroups of index two. For the variety $\cY$ this
is proved in \cite{MSY}. All these groups are arithmetic subgroups $\gG\inn
G_{\fR}$, such that the arithmetic quotient $\gG\bs \cD$ is a moduli space
of certain K3 surfaces with a marking on $H^2$ of some special type.

As a by-product we get the following result concerning moduli spaces of
Calabi-Yau threefolds:
\begin{theorem}[Corollary \ref{c158.1}]
  The family of hyperplane sections of $\cI_5$ passing through one of the
  36 triple points $p$ is, via projection, a family of Calabi-Yau
  threefolds which are degenerations of double octics, branched along the
  union $V_3\cup V_5$ of a cubic and a quintic surface which meet in 15
  lines.
\end{theorem}
In other words, the family of quintic hypersurfaces is connected with the
moduli space of double octics, along a four-dimensional sublocus.

Let us now make a few remarks about the entire set of examples. First of
all, they are (almost) all hypersurfaces. For a Baily-Borel embedding to be
a hypersurface, the singularities on $\xgs$ must be hypersurface
singularities --- which they usually are not.  Note that all of the
examples where $\xgs$ is a normal hypersurface (i.e., excluding perhaps the
dual of $\cB_4$), the arithmetic groups contain torsion. Roughly speaking,
this is experimental evidence of a statement like: $\gG$ has torsion
$\llra$ the singularities of $\xgs$ are very mild (hypersurface, complete
intersection); $\gG$ torsion free $\llra$ the singularities are not so
mild. Whether this has some general validity would seem to be quite a
difficult problem.

Secondly, all examples have interesting, decomposing hyperplane sections,
and the components of these reducible sections are modular subvarieties.
Let $\cX$ be one of our varieties, and let:
\begin{eqnarray*}\mu & = & \hbox{\# of hyperplanes $H$, such that $H\cap \cX$
    decomposes} \\ \nu & = & \hbox{\# of linear subspaces on $\cX$ cut out
    by them} \\ \tau & = & \hbox{\# of linear subspaces in each $H\cap
    \cX$} \end{eqnarray*} Then the results are summed up in Table
    \ref{table0}.
\begin{table}\begin{center}
    \caption{\label{table0} Special hyperplane sections in the examples}

    $$\begin{array}{|c|c|c|c|c|}\hline \cX & \mu & \nu & \tau & \hbox{rest
        of $H\cap \cX$} \\ \hline\hline \cS_3 & 15 & 15 & 3 & - \\ \hline
      \cI_4 & 10 & 10\ \fP^1\times \fP^1{'s} & 1 & - \\ \hline \cN_5 & 15 &
      15 & 3 & \hbox{quadric surface}\\ & 15 & 15 & 5 & - \\ \hline
      \cW_{10} & 10 & 10\ \fP^1\times \fP^1{'s} & 1 & - \\ & ? & ? & ? & ?
      \\ \hline \cB_4 & 40 & 40 & 4 & - \\ \hline \cB_4^{\vee} & 45 & 45\
      \fP^1\times \fP^1{'s} & 1 & ? \\ \hline \cI_5 & 27 & 45 & 5 & - \\
      \hline \cI_5^{\vee} & 36 & 120 & 10 & \cI_4 \\ \hline
\end{array}$$
{\small In this table we have included in addition to linear subspaces,
  also quadric surfaces, which are products of linear spaces, in the column
  titled ``$\nu$''.}
\end{center}
\end{table}

So this behavior looks like it can be expected in general, at
least for those $\gG$ which have torsion.  Some other curiosities of the
examples are the following.  The {\it Hessian variety} $\Hess(V)$ of a
hypersurface $V$ meets $V$ along the parabolic divisor, which is the locus
of points where the tangent hyperplane sections has worse than an ordinary
double point at the point of tangency; it gets blown down to a singular
locus of the dual variety under the duality map. The following strange
behavior occurs in our examples:
\begin{itemize}\item[i)] $\cS_3\cap \Hess(\cS_3)$ = 15 planes;
\item[ii)] $\cI_4\cap \Hess(\cI_4)$ = ten quadric surfaces;
\item[iii)] $\cB_4\cap \Hess(\cB_4)$ = 40 planes;
\item[iv)] $\cI_5\cap \Hess(\cI_5)$ = 45 $\fP^3$'s.
\end{itemize}
Furthermore, each of these intersections consists of modular subvarieties,
viewing the varieties as arithmetic quotients\footnote{it has not yet been
  rigorously proved that $\cI_5$ is an arithmetic quotient, as mentioned
  above}.  Another classical notion is that of {\it the Steinerian} of a
hypersurface: it is the locus of singular points of the quadric polars with
respect to $V$, which are singular. (The Hessian is the locus of points for
which the quadric polar is singular, and the vertices of these cones cut
out the Steinerian.) The following is even stranger:
\begin{itemize}\item[i)] $\cB_4$ is self-Steinerian;
\item[ii)] $\cI_5$ is self-Steinerian.
\end{itemize}
It is really not understood what self-Steinerian means geometrically, but
has something to do with the parabolic divisor. For i) we have a conceptual
proof, due to Coble; in the latter case this has just been checked
computationally.

As the expressed purpose of this paper is to point out the relation between
a set of (established) modular varieties and the quintic $\cI_5$, we spend
some time first recalling these varieties; their equations, moduli
interpratations and birational models. Although much of this material is
known to experts, we recall these for the convenience of the reader;
moreover there are several facts for which we prefer to give our own
proofs, making the paper in fact quite self-contained. This prepatory
material is collected in Part I, while the results on $\cI_5$ are given in
Part II.

We end this introduction by recalling a few facts about cubic surfaces and
the 27 lines on them, which is fundamental to all that follows.

\begin{center}
  {\bf Cubic surfaces}
\end{center}

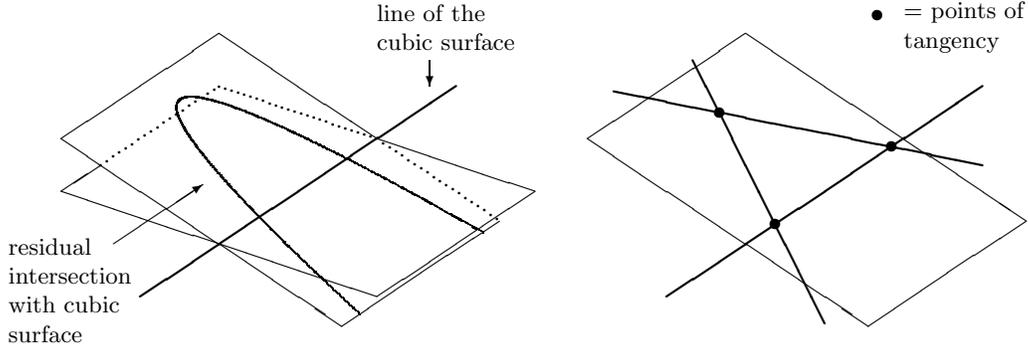
\begin{figure}[htb]

  $$ \unitlength0.7cm
\begin{picture}(10,8)

  \thicklines
\put(2.5,1.5){\line(3,2){6}}                               

\thinlines
\put(1.0,3.5){\line(3,-1){6}}                              
\put(1.0,4.5){\line(3,-2){5.33}} \put(4.0,6.5){\line(3,-2){3.0}}
\multiput(7.0,4.5)(0.12,-0.08){20}{\circle*{0.0001}}
\put(7.0,4.5){\line(3,-1){3}}
\multiput(4.0,5.5)(0.12,-0.04){25}{\circle*{0.0001}}

\put(1.0,3.5){\line(3,2){0.75}}                            
\multiput(1.75,4.0)(0.12,0.08){19}{\circle*{0.0001}}
\put(1.0,4.5){\line(3,2){3}} \put(7.0,1.5){\line(3,2){3}}
\put(6.33,0.944){\line(3,2){3}}

\bezier{600}(6.67,1.167)(-1.33,8.56)(9,2.722)             

\put(8.0,6.0){\vector(0,-1){0.5}} \put(7.0,6.5){\parbox{2cm}{\small line of
    the\\ cubic surface}} \put(2.0,2.5){\vector(3,2){1.7}}
\put(0.0,1.5){\parbox{2cm}{\small residual\\ intersection\\ with cubic\\
    surface}}

\end{picture}
%
%
\unitlength0.7cm
\begin{picture}(9.5,8)

  \thicklines
\put(2.5,1.5){\line(3,2){6}}                               
\put(3.5,5){\line(1,-2){2}} \put(3.5,5){\line(-1,2){0.5}}
\put(3.5,5){\line(5,-1){5}} \put(3.5,5){\line(-5,1){2}} \thinlines
\put(1.0,4.5){\line(3,-2){5.33}} \put(4.0,6.5){\line(3,-2){5.33}}

\put(1.0,4.5){\line(3,2){3}} \put(6.33,0.944){\line(3,2){3}}
\put(3.5,5){\circle*{0.2}} \put(4.56,2.88){\circle*{0.2}}
\put(6.77,4.35){\circle*{0.2}}

\put(6.5,6.85){\circle*{0.2}} \put(7.0,6.5){\parbox{2cm}{\small = points
    of\\ tangency}}


\end{picture}
$$
\caption[Sections of a cubic surface]{\label{f1}
  (a) generic section through a line\hspace{1.5cm} (b) tritangent section}
\end{figure}

Consider a smooth cubic surface $S$ in $\fP^3(\fC)$, given by the vanishing
of a cubic form $f$. There are exactly 27 lines on this surface, a fact
first pointed out by George Salmon in 1849, in response to a letter from
Arthur Cayley.  The argument he used is reproduced for example in
\cite{SF}, and is purely algebraic.  For the convienience of the reader we
sketch this briefly.  Consider hyperplane sections $H_t$ passing through
one of the lines.  Since $H_t\cap S=\{$3rd degree plane curve$\}$, which
already contains a line, it consists generically of a quadric and a line.
Since this intersection contains two double points, it follows that {\em
  every} such plane is a bitangent, i.e., tangent to $S$ in two points. For
a finite number of hyperplanes $H_t$ this intersection degenerates into
three lines, as illustrated in Figure \ref{f1}.

Such planes $H_t$ which are tangent to $S$ at three points (and contain
three lines) are called accordingly {\em tritangent planes}.  Now fix a
line $L$ on the cubic surface given in local coordinates by $x_3=x_4=0$,
say.  Then the equation $f=0$ can be written $x_3U+x_4V=0$ with $U,\ V$
quadratic. To find the tritangent planes one puts $x_3=\gm x_4$ into this
relation, divides by $x_4$ and forms the discriminant of the ensuing
equation. Viewing this as an equation in $\gm$, one sees easily that it has
degree 5, i.e., there are five values of $\gm$ for which the corresponding
plane is a tritangent, or in other words the given line lies in five
tritangent planes.  Now count: starting with a given tritangent plane, it
meets 4$\cdot$3=12 other tritangents, in each of which there are two other
lines, which gives 12$\cdot$2+3=27 lines on a smooth cubic surface.  Since
each line is contained in five tritangents, this leads to
27$\cdot$5$\div$3=45 tritangent planes to the cubic surface $S$.  In the
book \cite{SF} several other proofs of the magic number 27 are given.

We consider now the group of permutations, Aut(\cL), of the 27 lines (or of
the 45 planes), by which we mean the permutations of the lines preserving
the intersection behavior of the lines.  For this it is useful to consider
the famous double sixes and the notation for the 27 lines introduced by
Schl\" afli.  A {\em double six} is an array
$$N= \left[ \matrix{ a_1 & a_2 & a_3 & a_4 & a_5 & a_6 \cr b_1 & b_2 & b_3
    & b_4 & b_5 & b_6 \cr } \right] $$ of 12 of the 27 lines with the
property that two of these 12 meet if and only if they are in different
rows and columns.  (This notation distinguishes this particular set of 12,
although any such double six is equivalent to it under Aut(\cL)).  The
other lines are given by the $6 \choose 2$=15 $c_{ij}=a_ib_j\cap a_jb_i$,
where $a_ib_j$ denotes the tritangent spanned by those two lines.  There
are 36 double sixes, namely the $N$ above, 15 $N_{ij}$ and 20 $N_{ijk}$:
\begin{equation}\label{eB2.2}
\begin{minipage}{14cm}\begin{center}
    $N_{ij}=\left[ \matrix{ a_i & b_i & c_{jk} & c_{jl} & c_{jm} & c_{jn}
        \cr a_j & b_j & c_{ik} & c_{il} & c_{im} & c_{in} \cr } \right], $

    $N_{ijk}=\left[ \matrix{ a_i & a_j & a_k & c_{mn} & c_{ln} & c_{lm} \cr
        c_{jk} & c_{ik} & c_{ij} & b_l & b_m & b_n \cr} \right]. $
\end{center}
\end{minipage}
\end{equation}
Since a double six describes, by definition, the intersection behavior of
the lines, we see immediately that $\gS_6$ (the symmetric group on six
letters) acts by permutations on a double six and a $\fZ_2$ acts by
exchanging rows.  Since there are 36 double sixes, we see $|Aut({\cal
  L})|=|\gS_6|\cdot2\cdot$36=51,840. A natural question arising here is:
how many lines do two double sixes have in common?  The answer is twofold:
\begin{itemize}
\item either: four (like $a_1,\ a_2,\ b_1,\ b_2$, which $N$ and $N_{12}$
  have in common) which have the property of lying in pairs in planes, the
  pairs being however mutually disjoint;
\item or: six (like $a_1,\ a_2,\ a_3,\ b_1,\ b_2,\ b_3$, which $N$ and
  $N_{123}$ have in common) which form two triples.
\end{itemize}

Following Sylvester one speaks accordingly of {\em syzygetic} and {\em
  azygetic} pairs of double sixes.  A given double six is syzygetic to 15
others and azygetic to 20 others.  A pair of azygetic double sixes form
through the 12 lines they do {\em not} have in common a third double six,
which is azygetic to both.  There are 120 such triads of azygetic
double sixes, and 36$\cdot$15$\div$2=270 pairs of syzygetic double sixes.

We mention here a further geometric curiosity of the 45 planes.  Take two
of the tritangents which do not meet in a line on the cubic surface, say
$\ga_1,\ \ga_2$ (for example $a_1b_2c_{12}$ and $a_3b_4c_{34}$). These two
planes determine three others, denoted $\gb_1,\ \gb_2,\ \gb_3$, by the
property that each $\gb_i$ is determined as the tritangent containing each
a line of $\ga_1$ and of $\ga_2$ (in the example above, $a_1b_4,\ b_2a_3,\
c_{12}c_{34}$). The third line in each of $\gb_i$ which is not one of
$\ga_i$ all lie in a common tritangent (for example here
$c_{14}c_{23}c_{56}$).  Then this is a unique third tritangent denoted
$\ga_3$. The set of six tritangents $(\ga_1,\ga_2,\ga_3),\
(\gb_1,\gb_2,\gb_3)$ have the special property that the nine lines
$\ga_{\gm}\cap\gb_{\gn}$ together with a point (19 conditions) determine
$S$.  Such a set of six tritangents is called a {\em trihedral pair}, and
there are 120 such; this implies that the equation of the cubic surface $S$
can be written in 120 ways as
$$y_1y_2y_3+z_1z_2z_3=0.$$
  One finds the following types:

\begin{equation}\label{eB2.1}
\begin{array}{|ccc|ccc|ccc|} \hline
  & (20) & & & (10) & & & (90) & \\ \hline a_i & b_j & c_{ij} & c_{il} &
  c_{jm} & c_{kn} & a_i & b_j & c_{ij} \\ b_k & c_{jk} & a_j & c_{mn} &
  c_{ik} & c_{lj} & b_l & a_k & c_{kl} \\ c_{ik} & a_k & b_i & c_{jk} &
  c_{ln} & c_{im} & c_{il} & c_{jk} & c_{mn}\\ \hline \end{array}\quad,
\end{equation}
the rows (of each box) giving $\ga_1,\ \ga_2,\ \ga_3 $ and the columns
giving $\gb_1,\ \gb_2,\ \gb_3$.  This configuration is in some sense
complementary to the double sixes: {\em Starting with nine lines lying in
  such a trihedral pair, the remaining 18 lines form a unique azygetic
  triple of double sixes and conversely, the nine lines not contained in a
  given azygetic triple of double sixes always lie in a trihedral pair.}
Furthermore there are 40 {\em triads (triples)
  of trihedral pairs}, such that each such triad contains all 27 lines.
 One is lead to study {\em
  enneahedra}, i.e., 9-gons, by which one means sets of 9 tritangents which
contain all 27 lines.  It turns out that there are two types of such
enneahedra: 40 of the first type which have the property that the nine
tritangents form {\em four different} triads of trihedral pairs, and 160 of
the second type which belong to a unique triad. We sum up the configuration
in Table \ref{table20}.

\begin{table}
\caption{\label{table20}
  Loci assoicated with the 27 lines on a smooth cubic surface}
\begin{center}
\begin{tabular}{|l||l|} \hline
  \# objects & description of the objects \\ \hline \hline 27 & lines on a
  cubic surface \\ \hline 135 & intersection points of two of the lines \\
  \hline 216 & pairs of skew lines \\ \hline 36 & double sixes \\ \hline 45
  & tritangents \\ \hline 120 & trihedral pairs, set of six tritangents
  containing nine lines \\ \hline 40 & triples of trihedral pairs, set of
  18 tritangents containing all 27 lines \\ \hline 120 & triads of azygetic
  double sixes \\ \hline 270 & pairs of syzygetic double sixes \\ \hline
\end{tabular}
\end{center}
\end{table}

The modern point of view is to consider cubic surfaces as del Pezzo
surfaces of degree 6. The combinatorics of the 27 lines, as listed in Table
\ref{table20}, are then encoded in the Picard group of the del Pezzo
surface. In fact, the complement in $\Pic(S)$ to the hyperplane section,
call it $\Pic^0(S)$, is isomorphic to the root lattice of $E_6$. The
equation of the surface is given by the embedding of $S$ by means of the
linear system of elliptic curves through the six given points.

Fix six points in $\fP^2$, say $x=(p_1,\ldots,p_6)$, such that the $p_i$
are in general position, i.e., no three lie on a line, and not all six lie
on a conic.  Let $\hat{\fP}^2_x$ denote the blow up of $\fP^2$ at all six
points, $\grr_x:\hat{\fP}^2_x\lra \fP^2$. Consider the following curves as
classes in $\Pic(\hat{\fP}^2_x)$:
\begin{itemize}\item[i)] $a_1,\ldots,a_6$, the exceptional divisors over
  $p_1,\ldots,p_6$;
\item[ii)] $b_1,\ldots,b_6$, $b_i$ the proper transform of the conic $q_i$
  passing through all points $p_j,\ j\neq i$;
\item[iii)] $c_{ik}$, the proper transform of the line $\-{p_ip_k}$.
\end{itemize}
If we consider the surface $\hat{\fP}^2_x$, we have
$H^2(\hat{\fP}^2_x,\fZ)= [l]\fZ\oplus_i\fZ a_i$. Let $Q$ be the
intersection form on $H^2(\hat{\fP}^2_x,\fZ)$; then the classes $a_i,\
b_i,\ c_{ij}$ fulfill $Q(a_i,a_i) = Q(b_i,b_i) = Q(c_{ij},c_{ij})=-1$.  In
a well-known manner one takes a rank six subset, which is isomorphic to the
root lattice of type ${\bf E_6}$. For details see \cite{DO} and references
therein. Consider the orthocomplement of the anti-canonical class on
$\hat{\fP}^2_x$, and denote this by $\Pic^0(\hat{\fP}^2_x)$. Recall that
the anti-canonical class is $3l+\sum_{i=1}^6 a_i$, and that the
anti-canonical embedding of $\hat{\fP}^2_x$ is as a cubic surface.
Consequently we may view $\Pic^0(\hat{\fP}^2_x)$ as the orthocomplement of
the hyperplane section class of $\Pic(S_x)$, where $S_x$ is the cubic
surface which is the anti-canonical embedding. The following elements
$\gl$ with $Q(\gl,\gl)=-2$ form a basis of $\Pic^0(S_x)$:
\begin{equation}\label{eB3.1} \begin{array}{rcl} \ga_0 & = & l-a_1-a_2-a_3 \\
    \ga_1 & = & a_1-a_2 \\ \ga_2 & = & a_2-a_3 \\ \ga_3 & = & a_3-a_4 \\
    \ga_4 & = & a_4-a_5 \\ \ga_5 & = & a_5-a_6 \end{array}
\end{equation}
These also form a base of a root system of type ${\bf E_6}$, by taking
$\ga_1,\ldots,\ga_5$ as the sub root system of type $\bf A_5$. Since the
classes $a_i,\ b_i,\ c_{ij}$ are exceptional, they all represent elements
of $\Pic^0(S_x)$, and the 45 tritangents give 45 relations like
$a_i+b_j+c_{ij}=0$. This leads to the following exact sequence of
$\fZ$-modules
\begin{equation}\label{eB3.2} 0 \lra \fZ^{24} \lra \fZ^{45} \lra \fZ^{27}
  \lra \Pic^0(S_x) \lra 0
\end{equation}
which we will meet again later.

\noindent{\bf Aknowledgements:} Many of the results in Part II on the
geometry of $W(E_6)$ in $\fP^5$, for example the sequence (\ref{eB3.2}),
were attained (several years ago) in collaboration with Duco v.~Straten, who
unfortunately had to devote himself to more urgent matters. I have also had
discussions with W. Barth, Bert v.~Geemen, Igor Dolgachev, E. Looijenga, K.
Matsumoto, I. Naruki and Masaaki Yoshida on many of the matters presented
here, and it is my pleasure to thank them heartily.

\tableofcontents

\part{Projective embeddings of modular varieties}\label{chapter10}

\section{The tetrahedron in ${\Bbb P}^3$}
\subsection{Arrangements defined by Weyl groups}
Let $\Phi(G,T)\inn \tt^*$ be a root system of a simple group $G$ (over
$\fC$). Using notations as in Bourbaki we have the roots (for those systems
which will be of interest to us in the sequel)
\begin{equation}\label{e108.1}\begin{minipage}{14cm}\begin{itemize}
    \item[$\bf A_n$] $\pm (\ge_i-\ge_j),\ 1\leq i< j \leq n+1$;
    \item[$\bf B_n$] $\pm (\ge_i\pm \ge_j),\ \pm\ge_i, 1\leq i < j \leq n$;
    \item[$\bf C_n$] $\pm (\ge_i\pm \ge_j),\ \pm2\ge_i, 1\leq i < j \leq
      n$;
    \item[$\bf D_n$] $\pm (\ge_i\pm\ge_j),\ 1\leq i < j \leq n$;
    \item[$\bf F_4$] $\pm (\ge_i\pm \ge_j),\ \pm\ge_k,\ \pm{1\over
        2}(\ge_1\pm\ge_2\pm\ge_3\pm\ge_4), 1\leq i < j \leq 4,\
      k=1,\ldots,4$;
    \item[$\bf E_6$] $\pm (\ge_i \pm \ge_j), 1\leq i<j \leq 5,\ \pm{1\over
        2}(\ge_1\pm\ge_2\pm\ge_3\pm\ge_4\pm\ge_5-\ge_6-\ge_7+\ge_8)$, with
      an even number of ``$-$'' signs in the parenthesis;
\end{itemize}
\end{minipage}
\end{equation}
Each root $\ga$ determines an orthogonal plane $\ga^{\perp}$, and for any
arrangement $\bf X_n$,
\begin{equation}\label{e108.10} \cA({\bf X_n}):=\{\ga^{\perp} \Big| \ga
  \hbox{ a root}\}
\end{equation}
is a central arrangement in $\fC^n$, i.e., each of the planes passes
through the origin. This induces a projective arrangement in $\fP^{n-1}$,
as follows. Blow up the origin of $\fC^n$; the exceptional divisor is a
$\fP^{n-1}$. The {\em projective arrangement} is the union of the
intersections $[H]\cap \fP^{n-1}$ in the exceptional divisor, where $[H]$
is the proper transform of the hyperplane $H=\ga^{\perp}$ under the blow
up. The projective arrangements for $\bf B_n$ and $\bf C_n$ coincide, and
these arrangements are given in $\fP^{n-1}$ with projective coordinates
$(x_1:\ldots :x_n)$ as follows:
\begin{equation}\label{e108.2} \begin{minipage}{14cm}\begin{tabbing}
      $\cA({\bf A_n})$: \quad \= $\{ x_i=0,\ i=1,\ldots, n;\ x_i=x_j,\
      1\leq i<j\leq n\}$;\\ $\cA({\bf B_n})$: \> $\{ x_i=0,\ i=1,\ldots,
      n;\ x_i=\pm x_j,\ 1\leq i<j\leq n\}$;\\ $\cA({\bf D_n})$: \> $\{
      x_i=\pm x_j,\ 1\leq i<j\leq n\}$;\\ $\cA({\bf F_4})$: \> $\{ x_i=0,\
      i=1,\ldots, n;\ x_i=\pm x_j,\ 1\leq i<j\leq 4,\ {1\over 2}(x_1\pm
      x_2\pm x_3\pm x_4)\}$;\\ $\cA({\bf E_6})$: \> $\{ x_i=\pm x_j,\ 1\leq
      i<j \leq 5,\ {1\over 2}(x_1\pm x_2 \pm x_3\pm x_4 \pm x_5 + x_6)\}$.
\end{tabbing}
\end{minipage}
\end{equation}
For the arrangement of type $\bf A_n$ we have made the coordinate
transformation $x_1=\ge_1-\ge_{n+1}, \ldots, x_n=\ge_n-\ge_{n+1}$, so
$x_i-x_j=\ge_i-\ge_j$ for $1\leq i< j\leq n$, and for $\bf E_6$ we have
taken $x_6$ to replace $x_8-x_7-x_6$.

The arrangements above are the arrangements defined by the projective
reflection groups $PW({\bf X_n})$. Each hyperplane is the reflection plane
for the reflection on the corresponding root. From this point of view these
arrangements are studied in \cite{OS2}.

\subsection{Rank 4 arrangements}
As described above, the groups $W({\bf A_4})$, $W({\bf B_4})$, $W({\bf
  D_4})$ and $W({\bf F_4})$ give rise to projective arrangements in
$\fP^3$. They consists of ten, 16, 12 and 24 planes, respectively. They
may also be described as follows (see \cite{GS})):

\begin{tabbing}
  $\cA({\bf A_4})$: \quad \= four faces of a tetrahedron plus the six
  symmetry planes; \\ $\cA({\bf B_4})$: \> \parbox{12cm}{six faces of a cube
    plus the nine symmetry planes plus the plane at infinity;}\\ $\cA({\bf
    D_4})$: \> \parbox[t]{12cm}{six faces of a cube plus the six symmetry
    planes through two edges each, OR: eight faces of an octahedron plus
    three symmetry planes containing four vertices each plus the plane at
    infinity;}\\ $\cA({\bf F_4})$: \> \parbox[t]{12cm}{ ``desmic figure'':
    six faces of the cube, eight faces of an inscribed octahedron, nine
    symmetry planes and the plane at infinity; this is also determined by
    the regular 24-cell;}
\end{tabbing}
The combinatorial description of these arrangements can be encoded in
numbers:
\begin{equation}\label{e109.1} t_q(j):=\#\{\fP^j \hbox{'s of the arrangement
    through which $q$ of the reflection planes pass}\}.
\end{equation}
In the case of the above arrangements we have the following data
($t_q:=t_q(0)$, the number of points):
\begin{equation}\label{e109.2}\begin{minipage}{14cm}
\begin{tabbing}
  $\cA({\bf A_4})$: \quad \= $t_6=5,\ t_4=10;\ t_3(1)=10,\ t_2(1)=15.$\\
  $\cA({\bf B_4})$: \> $t_9=4, t_6=8,\ t_5=12,\ t_4=16;\ t_4(1)=6,\
  t_3(1)=16,\ t_2(1)=36$. \\ $\cA({\bf D_4})$: \> $t_6=12,\ t_3=12;\
  t_3(1)=16,\ t_2(1)=18.$ \\ $\cA({\bf F_4})$: \> $t_9=24,\ t_4=96;\
  t_4(1)=18,\ t_3(1)=32,\ t_2(1)=72.$
\end{tabbing}
\end{minipage}
\end{equation}
\begin{definition}\label{d109.1} An arrangement $\cA\inn \fP^n$ is said to be
  in (combinatorial) {\em general position}, if $t_q(j)=0$ for all $q>n-j$.
  All $\fP^j$'s $\inn \cA$ for which $j>n-q$ holds are the {\em
    singularities} of the arrangement. The singularities are {\it genuine}
  if they are not the intersection of higher-dimensional singular loci with
  one of the planes of the arrangement. The union of all genuine
  singularities is the {\em singular locus}.
\end{definition}
In the above arrangements we have the following singular loci:
\begin{equation}\label{e109.3}\begin{minipage}{14cm}
\begin{tabbing}
  $\cA({\bf A_4})$: \quad \= five singular points, ten singular lines;\\
  $\cA({\bf B_4})$: \> 12=4+8 (genuine) singular points, 22=6+16 singular
  lines; \\ $\cA({\bf D_4})$: \> 12 singular points, 16 singular lines; \\
  $\cA({\bf F_4})$: \> 24 (genuine) singular points, 50=18+32 singular
  lines.
\end{tabbing}
\end{minipage}
\end{equation}

\subsection{The tetrahedron}
Consider now the arrangement $\cA({\bf A_4})$ in (\ref{e109.2}). By
(\ref{e109.3}) the singular locus consists of five points and ten lines. We
introduce the following notation: $P_1=(1,0,0,0),\ P_2=(0,1,0,0),\
P_3=(0,0,1,0),\ P_4=(0,0,0,1),\ P_5=(1,1,1,1)$, and $l_{ij}$ will denote
the line joining $P_i$ and $P_j$. Each line contains two of the five
points, and at each of the points four of the ten lines meet. The
arrangement is {\em resolved} by performing the following birational
modification of $\fP^3$:
\begin{equation}\label{e109.4} \begin{minipage}{14cm}\begin{itemize}\item[a)]
      Blow up the five points, $\grr_1:\hat{\fP}^3\lra \fP^3$;
    \item[b)] Blow up the proper transforms of the ten lines,
      $\grr_2:\tilde{\fP}^3\lra \hat{\fP}^3,\ \grr:\tilde{\fP}^3\lra
      \fP^3$.
\end{itemize}
\end{minipage}
\end{equation}
In the resolution 15 exceptional divisors $E_1,\ldots,E_5$ and
$L_{12},\ldots, L_{45}$ are introduced. The $E_i$ are the proper transforms
of the exceptional divisors introduced under $\grr_1$, and are isomorphic
to $\fP^2$ blown up in the four points $(1:0:0),\ (0:1:0),\ (0:0:1),\
(1:1:1)$, as are the proper transforms $H_i$ of the ten planes of the
arrangement. The ten exceptional divisors $L_{ij}$ are each isomorphic to
$\fP^1\times \fP^1$. The symmetry group of $\tilde{\fP}^3$ consists of
projective linear transformations of $\fP^3$ which preserve the arrangement
$\cA({\bf A_4})$, together with certain {\em birational} transformations of
$\fP^3$ which are {\em regular} on $\tilde{\fP}^3$, i.e., which contain the
singular locus (\ref{e109.3}) with simple multiplicity in their
ramification locus. Hence the Weyl group itself, $W({\bf A_4})=\gS_5$
(symmetric group on five letters) is contained in the symmetry group. But
in fact, $\gS_6$ is the symmetry group, and the extra generator is a
permutation of one of the $E_i$ and $H_j$, which clearly can be done {\em
  on } $\tilde{\fP}^3$.

\subsection{A birational transformation}
Note that since each of the ten lines in (\ref{e109.3}) contains two of the
five points which are blown up under $\grr_1$, the normal bundle of the
proper transform of each line on $\hat{\fP}^3$ is $\cO(1-2)\oplus\cO(1-2)=
\cO(-1)\oplus \cO(-1)$. By general results of threefold birational
geometry, it follows that
\begin{equation}\label{e110.1} \begin{minipage}{14cm}\begin{itemize}\item[a)]
      The divisors $L_{ij}$ on $\tilde{\fP}^3$ may be blown down to an
      ordinary threefold rational point (node), i.e., a singularity given
      by the equation $x^2+y^2+z^2+t^2=0$, OR:
    \item[b)] The ten lines on $\hat{\fP}^3$ may be blown down to the nodes
      mentioned in a).
\end{itemize}\end{minipage}\end{equation}
In other words, there is a threefold which we denote by $T$, which contains
ten threefold nodes, with a birational triangle:
\begin{equation}\label{e110.2}
  \unitlength1cm
\begin{picture}(3,2)(0,-.5)
  \put(-.5,1){$\tilde{\fP}^3$}\put(0,1.1){\vector(1,0){1.5}} \put(.7,1.3){$
    \grr_2$} \put(1.8,1){$\hat{\fP}^3 \stackrel{\grr_1}{\lra} \fP^3$}
  \put(-.4,0){$\Pi_2$} \put(-.3,.9){\vector(1,-1){1}} \put(1.8,0){$\Pi_1$}
  \put(1.9,.9){\vector(-1,-1){1}} \put(.7,-.5){$T$}
\end{picture}
\end{equation}
The map $\Pi_2$ blows down the union of ten disjoint ``quadric surfaces''
(i.e., divisors isomorphic to $\fP^1\times \fP^1$) to ordinary nodes, while
$\grr_2$ blows these quadric surfaces down to ten disjoint lines, which
$\Pi_1$ then blows down to the same ten isolated nodes. The 5+10 divisors
$E_i$ and $H_j$ on $\tilde{\fP}^3$ have the following properties:
\begin{equation}\label{e111.1}\begin{minipage}{14cm} \begin{itemize}\item[a)]
      Each is isomorphic to $\fP^2$ blown up in four points;
    \item[b)] Each contains ten lines of intersection with the other 15,
      forming an arrangement in the blown up $\fP^2$ of ten lines meeting
      in 15 points.
    \item[c)] Under the birational map $\Pi_2$ each of the divisors $E_i$
      and $H_j$ are blown down to a $\fP^2$; the image of the ten lines of
      b) lie four at a time in each of these $\fP^2$'s, as the four
      $t_3$-points of the following arrangement, which is the union of the
      intersections of the given $\fP^2$ with the others:
      $$\unitlength1cm
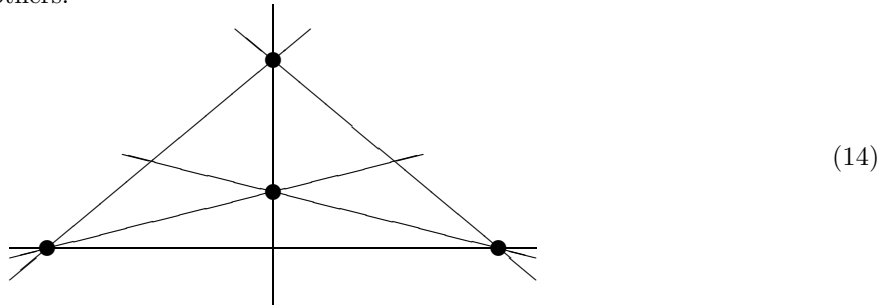
\begin{picture}(8,3.5)(1,0.8)

  \thinlines \put(2,1.5){\line(1,0){6.5}} \put(2,1.5){\line(4,1){5}}
  \put(2,1.5){\line(6,5){3.5}}

  \put(2,1.5){\line(-1,0){0.5}} \put(2,1.5){\line(-4,-1){0.5}}
  \put(2,1.5){\line(-6,-5){0.5}}

  \put(5,.75){\line(0,1){4}} \put(5,4){\line(6,-5){3.5}}
  \put(5,4){\line(-6,5){0.5}} \put(8,1.5){\line(-4,1){5}}
  \put(8,1.5){\line(4,-1){0.5}}

  \put(2,1.5){\circle*{.2}} \put(5,2.25){\circle*{.2}}
  \put(5,4){\circle*{.2}} \put(8,1.5){\circle*{.2}}
\end{picture}$$
\item[d)] The composition $\Pi_2\circ \grr^{-1}$ restricted to each of the
  planes $H_j$ is a usual Cremona transformation, blowing up three non-colinear
  points and blowing down the three joining lines. {\it Proof:} Take a face
  $H_j$ of the tetrahedron; $\grr_1$ blows up the three vertices it
  contains, so $\grr_1^{-1}(H_j)$ (the proper transform of $H_j$) is
  $\fP^2$ blown up in three points. Under $\grr_2$, a fourth point of
  $\grr_1^{-1}(H_j)$ is blown up, but it is blown down again under $\Pi_2$,
  as are the proper transforms of the three lines (in the plane $H_j$)
  joining the three vertices. By symmetry the same holds for all the $H_j$.
\end{itemize}
\end{minipage}
\end{equation}

It follows that on $T$, the images $\tilde{H}_j=\Pi_2(H_j)$ and
$\tilde{E}_i=\Pi_2(E_i)$ are copies of $\fP^2$, each containing four of the
ten nodes of $T$. Furthermore, in each of $\tilde{H}_j$ and $\tilde{E}_j$
we have the four $t_3$-points of the arrangement (\ref{e111.1}), which are
these four nodes of $T$. Finally, since there are 15 $\fP^2$'s, ten nodes
and four of them in each of the 15 $\fP^2$'s, there are five of these
divisors passing through a given node. Explicitly, take the node $n_{ij}$
corresponding to the line $l_{ij}$ in (\ref{e109.3}). Then it meets the
exceptional divisors $\~E_i,\ \~E_j$, as well as the three of the
$\~H_{\nu}$ for which $H_{\nu}$ contains the line $l_{ij}$.

\subsection{Fermat covers associated with arrangements}
Let $\cA\inn \fP^n$ be an arrangement of hyperplanes, i.e., a union
$\cA=\cup _{i=1}^k H_i$ of $k$ hyperplanes, and let $d\geq2 $ be an
integer. To the pair $(\cA,d)$ there is an associated function field
$\cL(\cA,d)$, an algebraic extension of the rational function field
$\cM(\fP^n)$. It defines, in a unique way, a branched cover $Y(\cA,d)\lra
\fP^n$, and a unique desingularisation $\~Y(\cA,d)$. The function field is
defined by:
\begin{equation}\label{e111a.1} \cL(\cA,d)=\fC\left({x_1 \over x_0},\ldots
    ,{x_n\over x_0}\right)\left[
    \sqrt{H_2/H_1}\hspace{-1.5cm}\raisebox{.2cm}{$\scriptstyle d$}
    \hspace{1.5cm},\ldots,
    \sqrt{H_k/H_1}\hspace{-1.5cm}\raisebox{.2cm}{$\scriptstyle d$}
    \hspace{1.5cm}\right],
\end{equation}
and the cover $Y(\cA,d)$ is the so-called Fox closure of the \'etale cover
over $\fP^n-\cA$ which is defined by (\ref{e111a.1}). $Y(\cA,d)$ is smooth
outside of the {\em singular locus} of $\cA$ (Definition \ref{d109.1}), and
the singularities of $Y(\cA,d)$ are resolved by resolving the singularities
of $\cA$. This is done by first blowing up all (genuine, i.e., not near
pencil) singular points, then all singular lines, and so forth. The
resolution (\ref{e109.4}) is a typical example. This is described in more
detail in the author's thesis; the desingularisation $\~Y(\cA,d)$ is the
fibre product in the following diagram:
\begin{equation}\label{e111a.2}
\begin{array}{ccc}\~Y(\cA,d) & \lra & Y(\cA,d) \\ \~{\pi}\downarrow &
  & \downarrow\pi \\ \~{\fP}^n & \stackrel{\grr}{\lra} & \fP^n
\end{array}
\end{equation}
where the horizontal arrows are modifications and the vertical arrows are
Galois covers with Galois group $(\fZ/d\fZ)^{k-1}$, which is the Galois
group of the field extension of (\ref{e111a.1}). $\grr$ is the modification
of $\fP^n$ which resolves the singularities of $\cA$. For example, each
singular point $P$ on $Y(\cA,d)$ is resolved by an {\em irreducible}
divisor $D_P$, which itself is a Fermat cover $Y(\cA',d)$, where $\cA'\inn
\fP^{n-1}$ is the arrangement induced in the exceptional $\fP^{n-1}$ which
resolves the point $P'=\pi(P)$. It consists of $k'$ planes, where $k'$ is
the number of the $k$ hyperplanes which meet at the point $P'$, and in the
process of resolving $Y(\cA,d)$, the cover $Y(\cA',d)$ is resolved also.
Hence on $\~Y(\cA,d)$ there is a {\em smooth} divisor $D_P$ which resolves
the singular point $P$ of $Y(\cA,d)$.

The singular covers $Y(\cA,d)$ can also be realised as complete
intersections, namely the intersections of $N=k-n$ Fermat hypersurfaces in
$\fP^{k-1}$:
\begin{eqnarray}\label{e111a.3}
  F_1 & = & a_{11}x_1^d+\cdots + a_{1k}x_k^d \nonumber\\ \vdots & & \vdots
  \\ F_N & = & a_{N1}x_1^d+\cdots +a_{Nk}x_k^d \nonumber
\end{eqnarray}
where $a_{11}H_1+\cdots +a_{1k}H_k, \ldots , a_{N1}H_1+\cdots + a_{Nk}H_k$
are the $(k-n)$ linear relations among the $k$ hyperplanes $H_i$. The map
$Y(\cA,d)\lra \fP^n$ is realised explicitly by the map
$(x_1,\ldots,x_k)\mapsto (x_1^d,\ldots,x_k^d)$.

\subsection{The hypergeometric differential equation}
The Fermat covers for the arrangement $\cA({\bf A_4})$ are closely related
to solutions of the hypergeometric differential equation on $\fP^3$, which
is an algebraic differential equation with regular singular points, whose
singular locus {\em coincides} with the arrangement $\cA({\bf A_4})$,
meaning that solutions are locally branched along the planes of the
arrangement.

First we introduce a new notation for the 15 surfaces $E_i,\ H_j$. These
can be numbered by pairs $(i,j),\ i<j\in \{0,\ldots,5\}$, with $E_i=H_{0i}$
and
\begin{equation}\label{e111b.1} H_{ij}\cap H_{kl}\neq \emptyset \iff
  i\neq j\neq k\neq l.
\end{equation}
We denote by $0i$ the point $P_i$ in $\fP^3$, and by $0ij$ the singular
line joining $0i$ and $0j$ in $\fP^3$. We then let $L_{0ij}$ denote the
exceptional divisor on $\~{\fP}^3$. We have (in $\fP^3$)
\begin{equation}\label{e111b.2} H_{ij}\cap H_{kl}=0mn \iff \{i,j,k,l\}\cap
  \{0,m,n\}=\emptyset.
\end{equation}
We want to consider branched covers $Y\lra \~{\fP}^3$ (with $\~{\fP}^3$ as
in (\ref{e109.4})), which are branched along the $H_{ij}$ and the
$L_{0ij}$.  Hence we let
\begin{equation}\label{e111b.3} n_{ij}:=\hbox{ branching degree along
    $H_{ij}$};\quad n_{0ij}:=\hbox{ branching degree along $L_{0ij}$},
\end{equation}
and of course $n_{ij},\ n_{0ij}\in \fZ\cup \infty$. (It makes sense to
allow negative branching degrees, as we will see below.)

To define the hypergeometric differential equation we may just as well work
on $\fP^n$ with homogenous coordinates $(x_0:\ldots:x_n)$, and consider the
arrangement $\cA({\bf A_n})$ of (\ref{e108.2}). Let $\gl_i\in \fQ,\ \
i=0,\ldots,n+1, \infty,\ \sum_i\gl_i=n+1$.  The hypergeometric differential
equation is:
\begin{equation}\label{e111b.4} \left\{
\begin{minipage}{14cm}$\ds(x_i-x_j)\del_i\del_jF + (\gl_i-1)(\del_iF-\del_jF) =
  0,\ 1\leq i<j \leq n$ \\ $x_i(x_i-1)\del_i^2F+P_i(x,\gl)\del_iF
  +(\gl_i-1)\sum{x_{\ga}(x_{\ga}-1) \over (x_i-x_{\ga})}\del_{\ga}F +
  \gl_{\infty}(1-\gl_i)F=0,\ 1\leq i \leq n$.
\end{minipage} \right.
\end{equation}
where
$$P_i(x,\gl)=x_i(x_i-1)\sum{1-\gl_{\ga} \over x_i-x_{\ga}} + \gl_0+\gl_i
-3-(2\gl_i+\gl_0+\gl_{n+1})x_i.$$ A solution of (\ref{e111b.4}) turns out
to be a period of an algebraic curve (the periods are many valued, as are
the solutions of (\ref{e111b.4})). The curve is
\begin{equation}\label{e111b.5}
  y^{\nu}=x^{\mu_0}(x-1)^{\mu_{n+1}}(x-t_1)^{\mu_1}\cdots (x-t_n)^{\mu_n},
\end{equation}
where the $\mu_i,\ \nu$ are related to the $\gl_i$ by the relation
\begin{equation}\label{e111b.6} {\mu_i \over \nu} = 1-\gl_i.
\end{equation}

The equation (\ref{e111b.4}) has an $(n+1)$-dimensional solution space,
spanned by $(n+1)$ periods of differentials of the curve (\ref{e111b.5}):
\begin{equation}\label{e111c.1} \go_i=\int_{\gg_i}{dx \over y},\quad
  <\gg_0,\ldots,\gg_n>=H^1(C,\fZ).\end{equation} Taking these gives a
homogenous many valued map
\begin{equation}\label{e111c.2} (\go_0,\ldots, \go_n):D\inn \~{\fP}^n
  \stackrel{\phi}{\lra} \fP^n,
\end{equation}
where $D$ is some Zariski open set (see (\ref{e111d.2}) below).  The map is
well-defined, since not all $\go_i$ vanish simultaneously. For very special
values of the parameters $\gl_i$, the image of $\phi$ is the complex ball
$\fB_n\inn \fP^n$ (this is just the Borel embedding of $\fB_n$ in its
compact dual).  In fact, one has the following theorem:
\begin{theorem}[\cite{DM}, \cite{T}]\label{t111c.1} If the following
  conditions are satisfied, then $\phi(D)=\fB_n$:
  $$\sum\mu_i=2,\quad \forall_{i,j}:\ (1-\mu_i-\mu_j)^{-1}\in \fZ\cup
  \infty.$$ In this case there exists a finite cover
  $$Y\lra D$$ branched along the total transform of $\cA({\bf A_n})$, which
  is a quotient $\gG\bs \fB_n$ with $\gG$ torsion free.
\end{theorem}
The integers $n_{ij}:=(1-\mu_i-\mu_j)^{-1}$ are then just the branching
degrees of $Y\lra D$ along the divisor $H_{ij}$. In fact the numbering
introduced in (\ref{e111b.1}) can be done analogously for any $n$.

In the special case of $\cA({\bf A_4})$ on $\fP^3$, the integers $n_{0ij}$
of (\ref{e111b.3}) are determined by the relation
$$n_{0ij}=2\left({1\over n_{kl}}+{1\over n_{lm}}+{1\over
    n_{km}}\right)^{-1},
$$ where the line $0ij$ is the intersection of $H_{kl},\ H_{lm},\ H_{km}$,
and these together with the $n_{ij}$ describe the branching degrees along
the entire branch locus. The solutions of the equations in Theorem
\ref{t111c.1} are as follows: \renewcommand{\arraystretch}{1.5}
\begin{equation}\label{e111d.1}\begin{array}{cl} 1) &
    {1\over 3}, {1\over 3}, {1\over 3}, {1\over 3}, {1\over 3}, {1\over 3},
    \\ 2) & {1 \over 2}, {1\over 2}, {1\over 4}, {1\over 4}, {1\over 4},
    {1\over 4} \\ 3) & {3 \over 4}, {1\over 4}, {1\over 4}, {1\over 4},
    {1\over 4}, {1\over 4} \\ 4) & {1\over 2}, {1\over 3}, {1\over 3},
    {1\over 3}, {1\over 3}, {1\over 6}\\ 5) & {3\over 8}, {3\over 8},
    {3\over 8}, {3\over 8}, {3\over 8}, {1\over 8} \\ 6) & {5\over 12},
    {5\over 12}, {5\over 12}, {1\over 4}, {1\over 4}, {1\over 4} \\ 7) & {7
      \over 12}, {5 \over 12}, {1\over 4}, {1\over 4}, {1\over 4}, {1\over
      4}
\end{array}
\end{equation}
\renewcommand{\arraystretch}{1.2} The set $D$ of (\ref{e111c.2}) is
determined as the complement of
\begin{equation}\label{e111d.2} H_{\infty}=\{H_{ij}\Big| n_{ij}=\infty; \
  L_{0ij}\Big| n_{0ij}=\infty\}\inn \~{\fP}^3.
\end{equation}
This is the locus which the uniformizing map (\ref{e111c.2}) maps onto the
{\em boundary} of $\fB_3\inn \fP^3$, i.e., $\phi(D)=\fB_n,\
\phi(H_{\infty})\inn \del\fB_n$. This requires of course that the
corresponding covers of the divisors on the cover $\~Y$ be abelian
varieties (as these are the compactification divisors on ball quotients).
This can happen as follows
\begin{equation}\label{e111d.3}\begin{minipage}{12cm}
\begin{itemize}\item[(i)] On one of the
  $H_{ij}$, this can occur if the branching degrees are: 2 for the six
  lines of (\ref{e111.1}), and $-4$ for the four exceptional curves.
\item[(ii)] On $L_{0ij}$, this can happen if $\mu_k+\mu_l+\mu_m=1,\
  \mu_0+\mu_i+\mu_j=1$.
\end{itemize}
\end{minipage}
\end{equation}
In the second case, the surface $S_{0ij}$ covering $L_{0ij}$ is of the form
$C_1\times C_2$, where $C_1\lra \fP^1$ (respectively $C_2\lra \fP^1$) is a
cover, with branching determined by $(\mu_k,\mu_l,\mu_m)$ (respectively
determined by $(\mu_0,\mu_i,\mu_j)$). It is an abelian variety $\iff$ both
curves $C_i$ are elliptic. Note that $Y\lra \~{\fP}^3$ will be a Fermat
cover $\iff$ all $n_{ij}$ coincide $\iff$ all $\mu_i$ conicide. In
particular,
\begin{proposition}\label{p111d.1} The only ball quotient in the list
  (\ref{e111d.1}) which is a Fermat cover which is a ball quotient is the
  solution 1), namely $Y(\cA({\bf A_4}),3)$ is a ball quotient.
\end{proposition}
\begin{remark} We will see later (see I3 following Lemma \ref{lq4.1} below)
  that the solution 4) gives rise also to a Fermat cover which is a ball
  quotient, namely $Y(\cA({\bf D_4}),3)$.
\end{remark}

\section{The Segre cubic ${\cal S}_3$}
In this section we will show that the variety $T$ of (\ref{e110.2}) has a
projective embedding as a cubic hypersurface known as the Segre cubic,
which we denote by $\cS_3$.

\subsection{Segre's cubic primal}
In $\fP^5$ with homogenous coordinates $(x_0:\ldots:x_5)$ consider the
locus
\begin{equation}\label{e111.3} \cS_3:=\{\sum_{i=0}^5x_i=0;\quad
  \sum_{i=0}^5x_i^3=0\}.
\end{equation}
As the first equation is linear, this shows that $\cS_3$ is a hypersurface,
i.e., $\cS_3\inn \fP^4=\{x\in \fP^5\big|\sum x_i=0\}$. Using
$(x_0:\ldots:x_5)$ as projective coordinates, the relation $x_5=-x_0-\cdots
-x_4$ gives the equation of $\cS_3$ as a hypersurface; however, the
equation in $\fP^5$ shows that $\cS_3$ is invariant under the symmetry
group $\gS_6$, acting on $\fP^5$ by permuting coordinates, which is not so
immediate from the hypersurface equation.

It is known that for any degree $d$ there is an upper bound on the number
of ordinary double points which a hypersurface of degree $d$ can have, the
so-called Varchenko bound. For cubic threefolds this number is ten, and it
has been known since the last century that $\cS_3$ is the {\em unique} (up
to isomorphism) cubic with ten nodes. The nodes on $\cS_3$ are given by the
points of $\fP^5$ for which three of the coordinates are 1 and the other
three are $-1$.  This is just the $\gS_6$-orbit of
\begin{equation}\label{e112.0} (1,1,1,-1,-1,-1).
\end{equation}
There is another interesting locus on $\cS_3$. Consider, in $\fP^5$, the
planes $P_{\gs}$ given by
\begin{equation}\label{e112.1} P_{\gs}=\{x_{\gs(0)}+x_{\gs(3)}=x_{\gs(1)} +
  x_{\gs(4)}=x_{\gs(2)}+x_{\gs(5)}=0\},
\end{equation}
where $\gs\in \gS_6$. There are 15 such $P_{\gs}$'s, the $\gS_6$-orbit of
\begin{equation} P_{id}=\{x_0+x_3=x_1+x_4=x_3+x_5=0\}.
\end{equation}
One checks easily that each $P_{\gs}$ contains four of the double points;
for example $P_{id}$ contains the following:
   $$(1,1,-1,1,-1,-1),\ (1,-1,1,1,-1,-1),\ (1,-1,-1,1,1,-1),\
    (1,-1,-1,1,-1,1).$$
Furthermore, the intersection of $P_{id}$ with the other $P_{\gs}$ is the
line arrangement (\ref{e111.1}). It is easily checked that each $P_{\gs}$
is contained entirely in $\cS_3$. One can also argue as follows. Any line
in $\fP^5$ which contains two of the nodes of $\cS_3$ meets $\cS_3$ with
multiplicity 4, hence is contained in $\cS_3$.  Similarly, each $P_{\gs}$
meets $\cS_3$ in the six lines of the arrangement (\ref{e111.1}), hence is
contained in $\cS_3$.

We just remark that the hyperplane sections $\{x_i=0\}$ of $\cS_3$ are
cubic surfaces with equation
\begin{equation}\label{e112.3}
  S_3 = \{ \sum_{i=0}^4x_i=\sum_{i=0}^4x_i^3=0\}.
\end{equation}
This cubic surface is known as the Clebsch diagonal surface and is a
remarkably beautiful object. It is the unique cubic surface having $\gS_5$
as symmetry group. The relation between $S_3$ and the icosahedral group was
studied by Hirzebruch. It turns out that $S_3$ is $A_5$-equivariantly
birational to the Hilbert modular surface for $\cO_{\fQ(\sqrt{5})}$, of
level $\sqrt{5}$.

Other interesting hyperplane sections are given by the hyperplanes
$\cT_{ij}=\{x_i-x_j=0\}$; indeed, $\cT_{ij}$ also contains four of the ten
nodes, hence $\cT_{ij}\cap \cS_3$ is a four-nodal cubic surface. This
four-nodal cubic surface is projectively unique, and is called the Cayley
cubic.

\subsection{A birational transformation}
\begin{theorem}\label{t113.1} The variety $T$ of equation (\ref{e110.2}) is
  biregular to $\cS_3$; the isomorphism $\psi:T\lra \cS_3$ defined below is
  $\gS_6$-equivariant.
\end{theorem}
{\bf Proof:} Following Baker \cite{Baker}, IV, p.~152, we define a
birational map
$$\gb:\fP^3- - \ra \cS_3.$$ Consider all quadric surfaces in $\fP^3$
passing through the points $P_i$ of (\ref{e109.3}). A base of this linear
system is given by the following degenerate quadrics. Let
$(z_0:\ldots:z_3)$ be homogenous coordinates on $\fP^3$, and set
\begin{equation}\label{e113.1}
\begin{array}{lll} \xi=z_0(z_3-z_1), & \eta=z_1(z_3-z_2), &
  \gz=z_2(z_3-z_0); \\ \xi'=z_1(z_3-z_0), & \eta'=z_2(z_3-z_1), &
  \gz'=z_0(z_3-z_2).
\end{array}
\end{equation}
These quadrics satisfy the relations $\xi+\eta+\gz=\xi'+\eta'+\gz'$ and
$\xi\eta\gz=\xi'\eta'\gz'$. Now change coordinates by setting
\begin{equation}\label{e113.a}
\begin{array}{lll} \xi=X+Y, & \eta=Y+Z, & \gz=X+Z; \\
  \xi'=-(X'+Y'), & \eta'=-(Y'+Z'), & \gz'=-(X'+Z').
\end{array}\end{equation}
Then the relations $\xi+\eta+\gz=\xi'+\eta'+\gz'$ and
$\xi\eta\gz=\xi'\eta'\gz'$ become
\begin{eqnarray}\label{e113.2} X+Y+Z+X'+Y'+Z' & = & 0 \\
  X^3+Y^3+Z^3+(X')^3+(Y')^3+(Z')^3 & = & 0. \nonumber
\end{eqnarray}
One sees this is just equation (\ref{e111.3}) of the Segre cubic. This
yields a rational map $\gb:\fP^3- - \ra \cS_3$, $\gb(z_0:z_1:z_2:z_3)=
(X,Y,Z,X',Y',Z')$. The base locus of the linear system of quadrics defining
$\gb$ (\ref{e113.1}) is the five points of (\ref{e109.3}), as the quadrics
all contain these points. It follows that $\gb$ blows up all five points,
the exceptional divisors $E_1,\ldots ,E_5$ being projective planes. Now
consider one of the ten lines of (\ref{e109.3}); for example, the one given
by $z_2=z_3=0$. Then $\eta=\gz=\eta'=\gz'=0$ and $\xi=\xi'=-z_0z_1$. In
other words, $\gb$ maps that line to the point $(1,0,0,1,0,0)$ in the
$(\xi,\eta,\gz,\xi',\eta',\gz')$ space, which is the point
$(1,1,-1,-1,-1,1)$ in the $(X,Y,Z,X',Y',Z')$ space. But that is just one of
the ten nodes of $\cS_3$. From $\gS_6$-symmetry we conclude that
$\gb:\fP^3- - \ra \cS_3$ coincides with the map $\Pi=\Pi_1\circ
\grr_1^{-1}$, with $\Pi_1$ as in (\ref{e110.2}) and $\grr_1$ as in
(\ref{e109.4}). In other words, $\gb=\Pi$ is the composition of morphisms
\begin{equation}\label{e113.3} \unitlength1cm \begin{picture}(3,2)
    \put(1.5,1.66){$\hat{\fP}^3$} \put(1.5,1.5){\vector(-1,-1){.9}}
    \put(1.7,1.5){\vector(1,-1){.9}} \put(.2,.33){$\fP^3$}
    \put(2.76,.33){$T,$} \put(.76,.3){$- - - - \ra$} \put(1.5,.66){$\Pi$}
    \put(.66,1.33){$\grr_1$} \put(2.33,1.33){$\Pi_2$}
\end{picture}
\end{equation}
and since (\ref{e113.2}) states that $\gb=\Pi$ maps onto $\cS_3$, this
gives an isomorphism $T\isom \cS_3$. Explicitly, $t\in T,\ t\mapsto
(\grr_1\circ\Pi^{-1}_2)(t)\mapsto
\gb((\grr_1\circ\Pi^{-1}_2)(t))=\psi(t)\in \cS_3$ is the desired map. The
$\gS_6$-equivariance follows from the fact that the whole diagram
(\ref{e113.3}) is $\gS_6$-equivariant.  \ende

Now consider the Picard group of $\cS_3$. From the explicit form of
birational map as given by Theorem \ref{t113.1} and (\ref{e110.2}), we see
that $\Pic(\cS_3)$ is generated by the image of the hyperplane class, call
it $H$, and the five exceptional classes $E_i$. It follows that
$\Pic(\cS_3)$ has rank 6, and the primitive part $\Pic^0(\cS_3)$, i.e., the
complement of the hyperplane class, has rank 5. The 15 classes $H_{ij}$
introduced in (\ref{e111b.1}) (these are the 15 linear $\fP^2$'s on $\cS_3$
noted in (\ref{e112.1})) give classes in $\Pic(\cS_3)$ and in
$\Pic^0(\cS_3)$. The 15 hyperplanes \begin{equation}
\label{e112b.3} \cH_{ij}=\{x_i+x_j=0\},
\end{equation}
each of which meets $\cS_3$ in three of the 15 $\fP^2$'s, give 15 {\em
  relations} in $\Pic^0(\cS_3)$: since $\cH_{ij}\cap \cS_3$ is a hyperplane
section, the sum of the three $\fP^2$'s cut out by $\cH_{ij}$, i.e.,
$H_{i_1,j_1}+H_{i_2,j_2}+H_{i_3.j_3}=\cH_{ij}\cap \cS_3$, is linearly
equivalent to the hyperplane class. This yields the following exact
sequence of $\fZ$-modules:
\begin{equation}\label{e112b.1} \begin{array}{ccccccccccc} 1 & \lra & K &
    \lra & \fZ\{\cH_{ij}\} & \lra & \fZ\{H_{ij}\} & \lra & \Pic^0(\cS_3) &
    \lra & 1 \\ & & \Big\|\wr & & \Big\|\wr & & \Big\|\wr & & \Big\|\wr & &
    \\ 1 & \lra & \fZ^5 & \lra & \fZ^{15} & \lra & \fZ^{15} & \lra & \fZ^5
    & \lra & 1. \end{array}\end{equation}
\begin{lemma}\label{l112b.1} In the sequence (\ref{e112b.1}), all
  $\fZ$-modules are $\gS_6$-modules, i.e., the exact sequence is one of
  $\gS_6$-modules.
\end{lemma}
{\bf Proof:} This is visible for the right three entries of the first
sequence in (\ref{e112b.1}), and it then follows for $K$. \ende Now
consider a generic hyperplane section of $\cS_3$; this is a smooth cubic
surface. Let $\nu:S=\cS_3\cap H\hra \cS_3$ denote the inclusion of the
section, and let $\nu*:H^2(\cS_3,\fZ)\lra H^2(S,\fZ)$ be the induced map on
cohomology. Then by the Lefschetz hyperplane theorem, this map is {\em
  injective}, and since both $S$ and $\cS_3$ are regular (i.e., not
irregular, that is, have no holomorphic one forms), we may view this as an
injective map of the Picard groups: $\Pic(\cS_3)\hra \Pic(S)$, and a
corresponding inclusion on the primitive part. Recall also that we have on
the cubic surface 27 generators (the 27 lines), 45 relations among these
(the 45 tritangents), and an exact sequence on $\Pic^0(S)$ as in
(\ref{eB3.2}).  All in all we get the following map of sequences as in
(\ref{e112b.1}):
\begin{equation}\label{e112b.2}\begin{array}{ccccccccccc}
    1 & \lra & \fZ^5 & \lra & \fZ^{15} & \lra & \fZ^{15} & \lra & \fZ^5 &
    \lra & 1 \\ & & \downarrow & & \downarrow & & \downarrow & & \downarrow
    & & \\ 1 & \lra & \fZ^{24} & \lra & \fZ^{45} & \lra & \fZ^{27} & \lra &
    \fZ^6 & \lra & 1. \end{array}\end{equation} where the right hand groups
are $\Pic^0(\cS_3)$ and $\Pic^0(S)$, respectively, and the down arrows are
inclusions (by Lefschetz). Note that this corresponds to a symmetry
breaking. Indeed, on the first sequence there is a symmetry group $\gS_6$
acting, as already noted, while on the group $\Pic^0(S)$, in fact on the
whole second sequence, the group $W(E_6)$ acts naturally, as is well-known.

\begin{proposition}\label{p112b.1} The ideal $\scI(10)$ of the ten nodes is
  generated by the five quadrics $\cR_{\gl}$ of the Jacobian ideal of
  $\cS_3$.
\end{proposition}
{\bf Proof:} The inclusion $Jac(\cS_3)\inn \scI(10)$ is obvious, and the
five elements of $Jac(\cS_3)$ are clearly independent. The fact that
$\scI(10)$ has rank 5 has been verified by standard basis computations
(with the algebra program Macaulay). \ende

\begin{corollary}\label{c112b.1} The ideal of the ten nodes of $\cS_3$,
  $\scI(10)$, is the Jacobian ideal of $\cS_3$. \ende
\end{corollary}

\subsection{Uniformisation}
In this section we will show that the Segre cubic $\cS_3$ is actually the
Satake compactification of a Picard modular variety. Let $K=\fQ(\sqrt{-3})$
be the field of Eisenstein numbers, and consider the $\fQ$-group
$G=U(3,1;K)$, the unitary group of a hermitian form on a four-dimensional
$K$-vector space with signature (3,1). Consider the arithmetic group
$\gG:=G_{\fZ}=U(3,1;\cO_K)\inn G(K)$, where $\cO_K$ denotes the ring of
integers in $K$. It acts on the three-ball with non-compact quotient $\xg$.
Consider the principal congruence subgroups $\gG(\sqrt{-3})$ and $\gG(3)$,
as defined in \cite{J}. These determine a corresponding level structure in
the sense of Definition 2.5 of \cite{J}.  Now note the following well-known
isomorphisms:
\begin{equation}\label{e114.1} \gG/\gG(\sqrt{-3})=\gS_6, \quad
  \gG(3)/\gG(\sqrt{-3})=(\fZ/3\fZ)^9.
\end{equation}
It follows from this that the corresponding quotients $X(a):=X_{\gG(a)}$,
$a=1,\sqrt{-3},3$, yield Galois covers
\begin{equation}\label{e114.2} X(3)\stackrel{(\fZ/3\fZ)^9}{\lra} X(\sqrt{-3})
  \stackrel{\gS_6}{\lra} X(1),
\end{equation}
which explicitly describe the level structures involved. As usual let
$X(a)^*$ denote the Satake compactification.
\begin{theorem}\label{t114.1} There is a commutative diagram
\begin{equation}\label{e115.0}
  \unitlength.4cm
\begin{picture}(18,18)
  \put(8,1){$\cS_3$} \put(3,5.7){$\tilde{\fP}^3$}
  \put(14,2.5){$X(\sqrt{-3})^*$} \put(8.8,7.7){$\-X(\sqrt{-3})$}
  \put(7.6,12){$Y^{\wedge}$} \put(3,17){$\tilde{Y}$}
  \put(14.5,14){$X(3)^*$} \put(9,19){$\-X(3)$}
  \put(3.5,16.5){\vector(1,-1){4}} \put(8.5,18.5){\vector(-3,-1){4.3}}
  \put(10,18.5){\vector(1,-1){4}} \put(14,14){\vector(-3,-1){4.3}}
  \put(9,18.5){\vector(0,-1){9.5}} \put(3,16.5){\vector(0,-1){9.5}}
  \put(8,11.5){\vector(0,-1){9.5}} \put(15,13.5){\vector(0,-1){9.5}}
  \put(8.5,7.5){\vector(-3,-1){4}} \put(9.5,7.5){\vector(1,-1){4}}
  \put(13.5, 2.5){\vector(-3,-1){4}} \put(3.5,5.5){\vector(1,-1){4}}
\end{picture}
\end{equation}
where the horizontal maps from right to left are isomorphisms, those from
left to right are birational, and the vertical maps are $(\fZ/3\fZ)^9$
covers.
\end{theorem}
{\bf Proof:} First we have, over $\fP^3$, a singular cover $T_{DM}$ defined
by the solution 1) of (\ref{e111d.1}). This is desingularised by blowing up
the $\fP^3$ along the singular locus of the arrangement, $\tilde{\fP}^3
\lla \tilde{T}_{DM}$. From the fact that all $\mu_i=1/3$, we see that all
$n_{ij}$ and $n_{0ij}$ are equal to three, that is $T_{DM}$ is the Fermat
cover $Y(\cA({\bf A_4}),3)$, and $\tilde{T}_{DM}$ is its desingularisation
$\tilde{Y}:=\tilde{Y}(\cA({\bf A_4}),3)$ as in (\ref{e111a.2}); see also
Proposition \ref{p111d.1}. By Theorems \ref{t111c.1} and
 \ref{t113.1}, $\tilde{Y}$ is the
desingularisation of the ball quotient $\gG'\bs \fB_3$, for some torsion
free group $\gG'$. Blowing $\tilde{Y}$ down from $\tilde{\fP}^3$ to $\cS_3$
gives the singular variety $Y^{\wedge}$, which we will see in a minite is
the Satake compactification of the ball quotient.  Hence we only need to
identify the groups and check the compactifications coincide. As to the
first, we start with
\begin{Lemma}\label{l115.1} Let $\GQ$ be an isotropic
  $\fQ$-form of $U(3,1)$, $G_{\fQ}\sim U(3,1;L)$, $L$ imaginary quadratic
  over $\fQ$, and let $\gG\inn G_{\fQ}$ a torsion free arithmetic subgroup
  with arithmetic quotient $\xg$, Baily-Borel compactification $\xgs$ and
  toroidal compactification $\xgc$.  Then the isomorphism class of a single
  compactification divisor determines the field $L$, and hence $G_{\fQ}$ up
  to isogeny.
\end{Lemma}
{\bf Proof:} First note that for $U(3,1)$ the parabolic (there is only one
conjugacy class of parabolics, as the $\fR$-rank is one) takes on the
particularly simple form
\begin{equation}\label{e115.1} \begin{minipage}{14cm} \hspace*{\fill} $
    P \isom (\cR\cK)\sdprod \cZ V, \quad \cR\isom \fR^{\times}, \quad
    \cK=SU(2)\times U(1)$ \hspace*{\fill}

    \hspace*{\fill} $\cZ= \fR, \quad V=\fC^2 $\hspace*{\fill}
\end{minipage}
\end{equation}
For the $\fQ$-form of $P$, it follows that $V_{\fQ}\isom L^2$ for some
imaginary quadratic field $L$, and for the arithmetic parabolic $\gG_P\inn
P,\ \gG_P\cap V_{\fQ}\inn (\cO_L)^2$ is some lattice. Furthermore, the
theory of toroidal embeddings shows that a compactification divisor of
$\xgc$ is of the form $\fC^2/(\gG_P\cap V_{\fQ})$, which has complex
multiplication by $L$, so its isomorphism class determines $L$, which was
to be shown. \ende Now an easy calculation shows what the compactification
divisors on $\tilde{Y}$ are. Namely, these are the irreducible components
of the inverse image in $\tilde{Y}$ of the exceptional divisors
$L_{0ij}\inn \tilde{\fP}^3,\ L_{0ij}\isom \fP^1\times \fP^1$. The local
geometry of the arrangement shows the branch locus in $L_{0ij}$ is of the
form $p_1^*(\cO(3)) \otimes p_2^*(\cO(3))$, i.e., of the form $\{0\}\times
\fP^1,\ \{1\}\times \fP^1,\ \{\infty\}\times \fP^1$ and $\fP^1\times
\{0\},\ \fP^1\times \{1\},\ \fP^1\times \{\infty\}$. It is well-known that
the elliptic curve $E\lra \fP^1$, branched at $(0,1,\infty)$ to degree 3,
with Galois group $\fZ/3\fZ$, is the elliptic curve $E_{\grr}=\fC/\fZ\oplus
\grr\fZ,\ \grr=e^{2\pi i / 3}$.  From this it follows
\begin{Lemma}\label{l115.2} The compactification divisors $\gD_i$ of
  $\tilde{Y}$ are products
  $$ \gD_i\isom E_{\grr}\times E_{\grr},$$ where $E_{\grr}$ is the unique
  elliptic curve with $\fZ/6\fZ$ as automorphism group, i.e.,
  $E_{\grr}=\fC/\fZ\oplus \grr\fZ=\{x^3+y^3+z^3=0\}$ and
  $\Aut(E_{\grr})=<\pm 1,\pm\grr, \pm\grr^2>$.\ende
\end{Lemma}
Note that the morphism $\tilde{Y}\lra Y^{\wedge}$ blows down the $\gD_i$ to
singular points (just as $\tilde{\fP}^3\lra \cS_3$ blows down the $L_{0ij}$
to the nodes of $\cS_3$) which lie over the nodes of $\cS_3$.  From this
and the well-known fact that the Satake compactification of a ball quotient
has only isolated, zero-dimensional singularities, which are resolved in a
torus embedding by means of complex tori or quotients thereof, we get the
following
\begin{Corollary}\label{c115.1} The variety $Y^{\wedge}$ is the Satake
  compactification of the quotient $\gG'\bs \fB_3$, with $\gG'\inn G_{\fQ}$
  and $G_{\fQ}$ isogenous to $U(3,1;K)$, $K$ the field of Eisenstein
  numbers as above. \ende
\end{Corollary}

Now that it is established that $\gG'$ is (isogenous to) an arithmetic
subgroup of $U(3,1;K)$, group-theoretic methods can be applied to determine
the arithmetic subgroup. This is done in detail in \cite{J}, Lemma 2.9 and
Theorem 2.11. The result is: $\gG'=P\gG(3)$, and the group $\gG_{\cS_3}$
giving rise to the Segre cubic is $\gG_{\cS_3}=P\gG(\sqrt{-3})$. This
yields the statements of the theorem on the arithmetic groups. The
compactification divisors of $\tilde{Y}$ coincide by Lemma \ref{l115.2}
with those of $\-X_{\gG'}$, and these are blown down under $\tilde{Y}\lra
Y^{\wedge}$ to the singularities on the Satake compactification which is
$Y^{\wedge}$, $Y^{\wedge}\isom X_{\gG'}^*$. The cover $\tilde{Y}\lra
\tilde{\fP}^3$ (respectively $Y^{\wedge}\lra \cS_3$) is now readily
identified with $\-X(3)\lra \-X(\sqrt{-3})$ (respectively
with $X(3)^*\lra X(\sqrt{-3})^*$) of (\ref{e114.2}), from the fact
that the branching loci,
degrees and group actions coincide. Details can be found in \cite{J}.
\ende

\subsection{Moduli interpretation}
Now applying Shimura's theory we get the following moduli description of
$\cS_3$ (see \cite{J}, \S2 for details).
\begin{theorem}\label{t116.1} Any point $x\in \cS_3-\{\hbox{ten nodes}\}$
  determines a unique isomorphism class of principally polarised abelian
  fourfolds with complex multiplication by $K=\fQ(\sqrt{-3})$ and a level
  $\sqrt{-3}$ structure. The signature of the complex multiplication is
  (3,1). Any point $x\in Y^{\wedge}-\{\hbox{inverse image under $\phi$ of
    (\ref{e115.0}) of the ten nodes}\}$ determines a unique isomorphy class
  of abelian fourfolds as above with a level 3 structure.
\end{theorem}
Moreover, the moduli interpretation of the 15 $\fP^2$'s on $\cS_3$ is given
in \cite{J}.
\begin{proposition}\label{p116.1} The 15 $\fP^2$'s on $\cS_3$ are
  compactifications of two-dimensional ball quotients which are moduli
  spaces of those abelian fourfolds $A_x^4$ as above which split:
  $$A_x^4\isom A_x^3\times E_{\grr}.$$ The intersections of the 15 planes
  determine moduli points of $A_x^4$ which further decompose, i.e., $A_x^3$
  splits.
\end{proposition}

\begin{remark} It is natural to ask whether, given a point $x\in \cS_3$, one
  can give the equations defining the abelian variety $A_x$ occuring in
  Theorem  \ref{t116.1}. In some sense one can.
  First it turns out the $A_x$ is the
  Jacobian of an algebraic curve, as described by the hypergeometric
  equation as in equation (\ref{e111b.5}). Since the parameters are by
  Proposition \ref{p111d.1} the set 1) in (\ref{e111d.1}), these curves
  have the form:
\begin{equation}\label{e116.1} C_{\tau}=\{y^3=\prod_{i=1}^6(x-t_i(\tau))\};
\end{equation}
$C_{\tau}$ obviously has an automorphism of order 3, given by $y\mapsto
\grr y$ with the third root of unity $\grr$. This yields an automorphism of
the Jacobian of $C_{\tau}$. Without much difficulty one finds
\begin{itemize}\item[(i)] $\Jac(C_{\tau})=A_{\tau}$ has complex
  multiplication by $\cO_K$, the signature is (3,1).
\item[(ii)] The automorphism group is $\cO_K^*$, and is given by
  multiplication by $\pm\grr$ in $\cO_K$.
\end{itemize}
The most direct way to see this is to write down the Jacobian of the curve
(\ref{e116.1}) and show that its periods have the complex multiplication. A
basis of the (1,0) differentials on $C_{\tau}$ written in the normal form
\begin{equation}\label{e116a.1} y^3=x(x-1)(x-t_1)(x-t_2)(x-t_3)
\end{equation}
is given by
\begin{equation}\label{e116a.3} \int{dx \over
    \sqrt[3]{x(x-1)(x-t_1)(x-t_2)(x-t_3)}};
\end{equation}
choosing a base of $H_1(C_{\tau},\fZ)$ and taking the integrals over the
elements of that base gives the Jacobian; the multiplication by $\grr$ is
then evident. Hence one may invoke Shimura's theory to conclude:
\begin{lemma}\label{l116a.1} The isomorphism classes of the Jacobians of the
  curves (\ref{e116a.1}) are given as the points of the arithmetic quotient
  $PU(3,1;\cO_K)\bs \fB_3$. Putting a $\sqrt{-3}$ level structure on the
  Jacobians yields the moduli space $\gG(\sqrt{-3})\bs \fB_3$.
\end{lemma}
The latter space has already been identified with the open subset of smooth
points on $\cS_3$.

The precise relation between the moduli {\em point} $\tau\in \fB_3$ and the
{\em values} of the $t_i$ has been derived for surfaces, i.e., for $\tau$
in one of the subballs covering one of the 15 $\fP^2$'s on $\cS_3$, by
Holzapfel.  The result is: there are automorphic forms $G_2,G_3$ and $G_4$
of indicated weights on $\fB_2$ such that
\begin{equation}\label{e116a.2}
  C_{\tau}=\{y^3=x^4-G_2(\tau)x^2-G_3(\tau)x-G_4(\tau)\},
\end{equation}
much akin to the Weierstra\ss\ equation for an elliptic curve. (The
variable $x$ in (\ref{e116a.2}) is of course different than that in
(\ref{e116a.1})).  There is no doubt a similar expression for $\tau\in
\fB_3$.
\end{remark}

\section{The Igusa quartic ${\cal I}_4$}
This variety has been known since the last century, and it is related to
the configuration in $\fP^4$ which is dual to the 15 hyperplanes of
(\ref{e112b.3}) and the 15 planes of (\ref{e112.1}) which they cut out on
$\cS_3$, and in fact $\cI_4$ is just the dual variety of $\cS_3$. It was
also known in the last century that the tangent hyperplane sections of
$\cI_4$ are Kummer surfaces, giving $\cI_4$ a moduli interpretation. Igusa,
in the 1960's, made this rigorous and showed that $\cI_4$ is the Satake
compactification of $\gG(2)\bs \fS_2$, the Siegel modular threefold of
level 2.  We begin by discussing the projective variety, then turn to
Igusa's results.

\subsection{The quartic locus associated to a configuration of 15 lines}
Let $l_{\gs}$ be the line dual in $\fP^4$ to the $\fP^2$ of (\ref{e112.1}),
and let $h_{ij}$ denote the point dual to $\cH_{ij}$ of (\ref{e112b.3}).
Then these 15 lines meet at the 15 points $h_{ij}$, and three of the 15
lines meet at each, corresponding to the three $\fP^2$'s which are
contained in each $\cH_{ij}$.  Furthermore, each of the 15 lines contains
three of the 15 points, as each $\fP^2$ is contained in three of the
$\cH_{ij}$. It is useful to introduce the following notation: each line is
given a notation $(ij)$, and two such lines $(ij),\ (kl)$ meet if and only
if the sets $(ij),\ (kl)$ are disjoint. Hence the 15 points are numbered by
{\em synthemes} $(ij,kl,mn)$ and the three lines meeting each point are the
indicated {\it duads} (pairs) $(ij),\ (kl),\ (mn)$.  Then there are ten
sets such as 23, 31, 12 and 56, 64, 45 with the property that the first and
last three do not meet, but each of the first meets each of the last.
Therefore the six lines are generators of a quadric surface
\begin{equation}\label{e117a.1}
  Q_{ijk}=\parbox{12cm}{quadric with $(ij), (jk), (ik)$ in one ruling and
    $(lm), (mn), (ln)$ in the other ruling}
\end{equation}
Then $Q_{ijk}$ lies in a $\fP^3$, and there are ten such, corresponding to
the ways of dividing the six numbers into two {\it triads} (triples).  Let
us denote the corresponding $\fP^3$ by $K_{ijk}$, so
\begin{equation}\label{e117a.1a} Q_{ijk}\inn K_{ijk}.
\end{equation}
Then each of the 15 lines is contained in four of the $K_{ijk}$, and six of
the $K_{ijk}$ meet at each of the 15 points.

Consider now a set of four mutually skew of the 15 lines, for example 12,
23, 24, 25. Then there will be a two-dimensional space of $\fP^2$'s which
meet all four lines (as we are in $\fP^4$, generically a plane and a line
will not intersect). Of all of these planes, there are exactly two passing
through a given point of space $x\in \fP^5$. The locus we are interested in
is:
\begin{equation}\label{e117a.2} \cQ:=\left\{x\in \fP^5 \left|
      \parbox{6cm}{ the two planes meeting four skew lines of the 15
        $(ij)$ and passing through $x$ {\em coincide}}\right.\right\}.
\end{equation}
If, as in (\ref{e113.1}), we take coordinates
$\xi,\eta,\gz,\xi',\eta',\gz'$ satisfying $\xi+\eta+\gz=\xi'+\eta'+\gz'$ as
coordinates on $\fP^4$, then the condition (\ref{e117a.2}) yields a locus
with equation (\cite{Baker}, p.~125):
\begin{equation}\label{e117a.3} \sqrt{(\eta-\gz')(\eta'-\gz)} +
  \sqrt{(\gz-\xi')(\gz'-\xi)} + \sqrt{(\xi-\eta')(\xi'-\eta)}=0.
\end{equation}

To find the dual variety of the locus $\cQ$, Baker does the following.
Letting $a,b,c$ be variables, $a'=(1-a),\ b'=(1-b),\ c'=(1-c)$, consider
the six points which are the vertices of a coordinate simplex in $\fP^5$, and
call them $A, B, C, A', B', C'$. Then any point of our $\fP^4$ can be
written as $x=A/bc'+B/ca'+ C/ab'+ A'/b'c+ B'/c'a+ C'/a'b$. Calculating the
tangent plane of $\cQ$ at a point $x\in \cQ$ which satisfies
(\ref{e117a.3}), in terms of the coordinates used in (\ref{e117a.3}), one
gets:
\begin{equation}\label{e117a.4}
  bc'\xi+ca'\eta+ab'\gz-b'c\xi'-c'a\eta'-a'b\gz'=0.
\end{equation}
Now putting $u=bc', v=ca', w=ab', u'=-b'c, v'=-c'a, w'=-a'b$, the equation
becomes
\begin{equation}\label{e117a.5} u\xi+v\eta+w\gz+u'\xi'+v'\eta'+w'\gz'=0,
\end{equation}
with the two identities
\begin{equation}\label{e117a.6} u+v+w+u'+v'+w'=0,\quad uvw+u'v'w'=0.
\end{equation}
Since the identities (\ref{e117a.6}) do not depend on the point, it follows
that these equations define the dual variety. Now comparing with
(\ref{e113.1}), we have
\begin{proposition}\label{p117a.1} The dual variety of the quartic locus
  $\cQ$ is the Segre cubic $\cS_3$.
\end{proposition}

It is easy to see that $\cQ$ is singular along the 15 lines.  It was also
noted classically that a tangent hyperplane section of $\cQ$ is a Kummer
quartic surface, with 16 nodes, 15 from the intersections with the 15
singular lines, and one from the point of tangency.

\subsection{Igusa's results}
The relation to the Kummer quartic surfaces is correctly understood by
studying theta constants for the theta functions with 1/2-characteristics.
This was done by Igusa in \cite{igusa}, and we now recall some of his
results.

\subsubsection{Theta functions}
Let $\tau\in \fS_g=\{M\in M_g(\fC)\big| \tau={^t(\tau)},\, \Im(\tau) \hbox{
  positive definite}\}$, $z\in \fC^g$, and $m=(m',m'')\in \fQ^{2g}$. Note
that $\fS_g$ is a hermitian symmetric space of type $\bf III_g$.
\begin{definition}\label{d117.1} The {\em theta function of degree $g$ and
    characteristic $m$} is defined by the power series
  $$\gt_m(\tau,z)=\sum_{n\in \fZ^g}\exp\left({1\over 2}{^t(n+m')}\tau(n+m')
    + {^t(n+m')}(z+m'')\right).$$
\end{definition}
As a function of $\tau$ the series $\gt_m$ converges precisely for $\tau\in
\fS_g$, while as functions of $z$ by fixed $\tau$ these are theta functions
on $A_{\tau}=\fC^g/(\fZ^g+\tau\fZ^g)$. As such the zeros on $A_{\tau}$ are
determined by the characteristic $m$. The corresponding {\em theta
  constant} is \begin{equation}\label{e117.2} \gt_m(\tau):=\gt_m(\tau,0).
\end{equation}
Igusa has studied in \cite{igusa} these theta constants, in particular the
theta functions with characteristics $m\in {1\over 2^n}\fZ$. Some of his
results are the following.
\begin{lemma}\label{l3.1.2} $\gt_m(\tau)\equiv 0 \iff m\hbox{\em mod}(1)$
  satisfies $\exp(4\pi i{(^tm')}m'')=-1.$
\end{lemma}
The Siegel modular group $\gG_g(1)=Sp(2g,\fZ)$ acts on the arguments
$(\tau, z)$ as follows:
\begin{equation}\label{e117.3} M=\left(\begin{array}{cc}A & B \\ C & D
\end{array}\right),\ M(\tau,z)=\left((A\tau+B)(C\tau+D)^{-1}, (C\tau +
D)^{-1}z\right),
\end{equation}
and on the characteristic itself by \begin{equation}\label{e117.4}
  M(m)=\left(\begin{array}{cc}D & -C \\ -B & A \end{array}\right)\cdot m
  +{1 \over 2} \left(\begin{array}{c} \diag(C^tD) \\ \diag(A^tB)\end{array}
  \right).
\end{equation}
The behavior of the theta functions under $M$ is given by
\begin{lemma}\label{l117.1} Let $M\in \gG_g(1)$ act on $(\tau,z)$ as in
  (\ref{e117.3}) and on the characteristic $m$ as in (\ref{e117.4}). Then
  the theta functions transform according to the rule:
\begin{equation} \gt_{M(m)}(M(\tau,z))  =  \gk(M)\exp(2\pi i
  \phi_m(M))\det(C\tau+D)^{1/2} \times \exp(\pi i {^tz}(C\tau
  +D)^{-1}Cz)\gt_m(\tau,z),
\end{equation}
where $\gk(M)$ is some eighth root of unity and $\phi_m(M)$ is defined by
the formula
$$\phi_m(M)=-{1\over2}{^tm'}BDm' + {^tm''}{^tA}Cm'' - 2{^tm'}{^tB}Cm''
-{^t\diag}(A^tB)(Dm'-Cm'').$$
\end{lemma}
In particular for the theta constants the formula becomes
\begin{equation}\label{e117.5} \gt_{M(m)}(M\tau)=\gk(M)\exp(2\pi i
  \phi_m(M))\det(C\tau+D)^{1/2}\gt_m(\tau).
\end{equation}

What the equation (\ref{e117.5}) says for $g=2$ is that up to an eighth
root of unity, non-vanishing theta constants with 1/2-characteristics are
automorphic forms of weight 1/2 for the main congruence subgroup of level 2
in $Sp(4,\fZ)$.  Indeed, for $M\in \gG(2)$, it holds that $e^{2\pi
  i\phi_m(M)}=1$, as Igusa shows. There are 16 characterstics $m$; six are
{\em odd} (i.e., $\gt_m(\tau,z)=-\gt_m(\tau,-z)$) so give rise to vanishing
theta constants, while ten are even. The fourth powers $\gt_m^4$ are
genuine automorphic forms for $\gG(2)$, and determine a morphism
\begin{equation}\label{e118.2} f:\gG(2)\bs \fS_2 \lra \fP^9=(\gt_{m_1}^4:
  \cdots :\gt_{m_{10}}^4),
\end{equation}
where $m_1,\ldots,m_{10}$ are the ten even characteristics.

\subsubsection{The ring of automorphic forms}
Among the ten coordinate theta functions there are five linear relations,
the Riemann relations. This implies that the map $f$ in (\ref{e118.2}) maps
into a $\fP^4$, displaying the quotient $X_{\gG(2)}$ as a hypersurface. In
fact, since this is an embedding by means of automorphic functions whose
closure $X_{\gG(2)}^*\inn \fP^4$ is normal (see below), it follows that $f$
gives a Baily-Borel embedding of the arithmetic quotient. The proof that
$f$ is an embedding given by Igusa is quite deep, involving showing that
the ring of modular forms of $\gG(2)$ is the integral closure of the ring
generated by the said theta functions. More precisely, his result is
\begin{theorem}[\cite{igusa},p.~397]\label{t118.1}
  Take as coordinates in $\fP^4$ the following theta constants:
  $$y_0=\gt^4_{(0110)}(\tau),\ y_1=\gt^4_{(0100)}(\tau),\
  y_2=\gt^4_{(0000)}(\tau),$$
  $$y_3=\gt^4_{(1000)}(\tau)-\gt^4_{(0000)}(\tau),\
  y_4=-\gt^4_{(1100)}(\tau)-\gt^4_{(0000)}(\tau),$$ where we let
  $(ijkl)$ denote the characteristic $({i\over 2}{j\over 2}{k\over
    2}{l\over 2})$.  Set also $$\chi_{10}=\prod_{\hbox{even
      $m$}}\gt_m^2.$$ Then the ring of modular forms of $\gG(2)$ is given
  by:
  $$R(\gG(2))=\fC[y_0,\ldots,y_4,\chi_{10}]/\cE,$$ where $\cE$ is the ideal
  generated by the following two relations:

\begin{minipage}{14cm}$\cE=\left\{\begin{array}{rcl}
      R_1 & = & (y_0y_1+y_0y_2+y_1y_2-y_3y_4)^2-4y_0y_1y_2(\sum y_i) \\ R_2
      & = & \chi_{10}^2-{1\over 4}s(y_0,\ldots,y_4),\ s \hbox{ homogenous
        of degree 5}
\end{array}\right.$
\end{minipage}
\end{theorem}
However, the formula $R_1$ relating the theta functions was known long
before Igusa. Since the five linear relations determining the image $\fP^4$
of $f$ are known, it is sufficient to give a single relation of minimal
degree among the $\gt^4_m$ to determine the image. This relation can be
found as early as in the 1887 paper of Maschke \cite{maschke},
p.~505\footnote{the equation is somewhat hidden: ``...da\ss\ dagegen die
  symmetrische Function vierter Dimension sich bis auf einen Zahlenfactor
  als das Quadrat der zweiten Dimension erweist.''}. In terms of the theta
constants above, this equation is
\begin{equation}\label{e118.1}
  \left(\sum\gt^8_m\right)^2-4\left(\sum \gt_m^{16}\right)=0,
\end{equation}
which, as can be checked, is the same quartic as that given by $R_1$ in
\ref{t118.1}, as well as that given by (\ref{e117a.3}).
\begin{definition}\label{d118.1} The {\em Igusa quartic} $\cI_4$ is the
  quartic threefold defined in $\fP^4$ by the relation $R_1$ of Theorem
  \ref{t118.1} or the equation (\ref{e118.1}).
\end{definition}
As a corollary we have
\begin{corollary}\label{c118.1} The Igusa quartic
  $\cI_4$ and the quartic locus $\cQ$ of (\ref{e117a.3}) coincide, and this
  quartic is the Satake compactification of $X_{\gG(2)}$.
\end{corollary}

Hence we have described $X_{\gG(2)}^*$ as a singular quartic hypersurface
in $\fP^4$. There are the two interesting loci:
\begin{itemize}\item[(i)] the singular locus, which is the boundary of the
  Baily-Borel embedding of $X_{\gG(2)}$;
\item[(ii)] the intersection of $\cI_4$ with the coordinate hyperplanes in
  $\fP^9$, which are the modular subvarieties $\-Y_m(2)$ of \cite{J}, Thm.
  3.19; these are quotients of symmetric subdomains isomorphic to a product
  of discs.
\end{itemize}
As already mentioned, the singular locus of $\cI_4$ consists of 15 lines;
this can be directly calculated from the equation. Alternatively, applying
general formula for the number of cusps (see for example \cite{Yam}) we see
that $X_{\gG(2)}$ has 15 one-dimensional boundary components and 15
zero-dimensional boundary components; by \ref{c118.1} this is then the
singular locus of $\cI_4$. (That these boundary components are rational
curves is obvious ($\gG(2)\bs \fS_1$ is rational); that they are actually
{\em lines} is not so obvious, but an easy calculation). This line of
reasoning also requires the result, also due to Igusa, that, although
$\gG(2)$ is not torsion-free, there are nonetheless no singularities on
$X_{\gG(2)}$.

\subsection{Moduli interpretation}
The embedding (\ref{e118.2}) of $X_{\gG(2)}^*$ as the quartic $\cI_4$ shows
that $\cI_4$ has a moduli interpretation. In fact, $X_{\gG(2)}$ is a rough
moduli space of principally polarised abelian surfaces with a level 2
structure. However, $\gG(2)$ contains torsion, namely the element $-1$, so
$X_{\gG(2)}$ is {\em not} a fine moduli variety. This corresponds to the
fact that the automorphism $z\mapsto -z$ of $A_{\tau}$ {\em preserves} the
level 2 structure, hence the actual {\em object} which is parameterised by
$X_{\gG(2)}$ is the {\em quotient} $A_{\tau}/(z\mapsto -z)$. This is just
the Kummer quartic surface which already occured above. The precise
relation is given by
\begin{theorem}\label{t119.1} For a point $x\in \cI_4-\{\hbox{intersections
    of $\cI_4$ with the ten coordinate planes in (\ref{e118.2})}\}$, the
  corresponding Kummer quartic surface $K_x=A_{\tau}/\{\pm1\}$, where
  $x=p(\tau)$ for the natural projection $p:\fS_2\lra X_{\gG(2)}$, is the
  intersection of $\cI_4$ with the tangent hyperplane at $x$, $T_x\cI_4$:
  $$K_x=\cI_4\cap T_x\cI_4.$$
\end{theorem}
This statement can be found for example in \cite{Baker}. It amounts to the
fact, true in any dimension, that for $n\geq3$ the theta functions with
characteristics $\in \fZ/n\fZ$ on a fixed $A_{\tau}$ give an embedding of
$A_{\tau}$, while for $n=2$ they map onto the Kummer variety.

The reason one must exclude the ten hyperplane sections in Theorem
\ref{t119.1} is the following result.
\begin{proposition}\label{p120.1} The ten hyperplane sections $\{\gt_m^4=0\}
  \cap \cI_4$ are tangent hyperplane sections, i.e., the intersection is of
  degree 2 and multiplicity 2.
\end{proposition}
A proof, based only on the equation of $\cI_4$, can be found in
\cite{Baker}.  To understand the meaning of this, note that a general
hyperplane section meets $\cI_4$ in a quartic surface, while the
intersections here are quadric surfaces, hence to preserve degree must be
counted twice (i.e., multiplicity 2). Consider the symmetric subdomain
$\fS_1\times \fS_1\inn \fS_2$, which in this case is the set of reducible
matrices:
$$\fS_1\times \fS_1 = \left\{\left(\begin{array}{cc} \tau_1 & 0 \\ 0 &
      \tau_2\end{array} \right)\right\}\inn \left\{\left(\begin{array}{cc}
      \tau_1 & \tau_{12} \\ \tau_{12} & \tau_2 \end{array}\right)\right\}
=\fS_2.$$ Then an easy calculation shows that the theta function of
Definition \ref{d117.1} is a {\em product} of two theta functions of a
single variable (i.e., $z\in \fC$).  This is equivalent to the fact that
for reducible $\tau \in \fS_2$, the abelian surface $A_{\tau}$ is a product
of two elliptic curves, $A_{\tau}=E_1\times E_2$. In this case, the map
given onto the ``product Kummer'' variety is a map $s:E_1\times E_2\lra
E_1/\{\pm1\}\times E_2/\{\pm1\}=\fP^1\times \fP^1$, and this $\fP^1\times
\fP^1$ is the quadric surface occuring in \ref{p120.1}. Since $\fP^1\times
\fP^1$ has no moduli, we see that {\em formally} the statement of Theorem
\ref{t119.1} remains true for all $x\in \cI_4-\{\hbox{15 singular
  lines}\}$, if we consider product Kummer varieties instead of the usual
ones, and the hyperplane section is the quadric surface of Proposition
\ref{p120.1}. Note however, that this quadric surface, being a modular
subvariety, can also be described as:
\begin{equation}\label{e120.1} E_1/\{\pm1\}\times E_2/\{\pm1\}\isom
  (\gG_1(2)\bs \fS_1)^*\times (\gG_1(2)\bs \fS_1)^*,
\end{equation}
describing the product Kummer surface of a reducible abelian surface as a
compactification of an arithmetic quotient, that is, as a Janus-like
variety. We then get the following moduli interpretation of the quadric
surfaces.
\begin{proposition}\label{p120.2} The ten quadric surfaces of Proposition
  \ref{p120.1} are modular subvarieties which correspond to abelian
  surfaces which split. More precisely, for any $x$ on one of the quadric
  surfaces, but not on any of the singular lines (there are six such
  singular lines on each quadric surface, see (\ref{e117a.1})), determines
  a smooth abelian surface which splits, with a level 2 structure.
\end{proposition}
Finally we note that this geometry can be described, as discussed already
in \cite{J} and many other places, in terms of the finite geometry of
\begin{equation}\label{e121.1} V=(\fZ/2\fZ)^4.
\end{equation}
Let $<\ ,\ >$ denote the induced symplectic form on $V$; every vector $v\in
V$ is isotropic with respect to $<\ ,\ >$. Since there are 15 non-zero
vectors, there are 15 one-dimensional boundary components. Similarly, there
are 15 isotropic planes in $V$, giving 15 zero-dimensional boundary
components.  The modular subvarieties of Proposition \ref{p120.2}
correspond in this setting to non-singular pairs $\{\gd,\gd^{\perp}\}$,
where $\gd$ is a two-dimensional subspace of $V$ on which $<\ ,\ >$ is non
degenerate, and $\gd^{\perp}$ denotes the orthocomplement with respect to
$<\ ,\ >$. Of these there are exactly ten, as is easily checked. We leave
further details to the reader.

\subsection{Birational transformations}\label{section3.4}
We have seen above in Proposition \ref{p117a.1} that $\cS_3$ and $\cI_4$
are dual varieties. It follows from general theory that they are then in
fact {\em birational}. In this section we describe the ensuing birational
map explicitly. We consider the following modifications of $\fP^4$.
\begin{itemize}\item[a)] Blow up the ten nodes (\ref{e112.0}) of $\cS_3$;
  denote this by $\grr_1:\hat{\fP}^4\lra \fP^4$. There are ten exceptional
  $\fP^3$'s, each with normal bundle $\cO_{\fP^3}(-1)$. Consider one of the
  15 hyperplanes $\cH_{ij}$ of (\ref{e112b.3}). Since each hyperplane
  contains $4+2+1=7$ nodes, its proper transform on $\hat{\fP}^4$ is a
  $\fP^3$ blown up in those seven points; each of the 15 $\fP^2$'s of
  (\ref{e112.1}) lying on $\cS_3$ contains four of the nodes, so their
  proper transforms are copies of $\fP^2$ blown up in four points. Finally,
  let $\hat{\cS}_3\inn \hat{\fP}^4$ denote the proper transform of $\cS_3$
  in $\hat{\fP}^4$; $\hat{\cS}_3$ is smooth, and
  ${\grr_1}_{|\hat{\cS}_3}:\hat{\cS}_3\lra \cS_3$ is a desingularisation of
  $\cS_3$, replacing each node with a quadric surface $\isom \fP^1\times
  \fP^1$.
\item[b)] Let $\scI(15)$ denote the ideal of the 15 singular lines of
  $\cI_4$; blow up $\scI(15)$, and let $\grr_2:\hat{\hat{\fP}}^4\lra \fP^4$
  denote this modification. Under $\grr_2$, each of the lines is replaced
  by a $\fP^2$-bundle over that line, and each point is replaced by a union
  of $\fP^1$'s, one each for each {\em pair} $(l_1,l_2)$ of {\em lines}
  meeting at the point; this mentioned $\fP^1$ is then the intersection of
  the fibre $\fP^2$ of $\grr_2$ at that point with the (two) exceptional
  $\fP^2$-bundles over the lines $l_1$ and $l_2$. Note that the proper
  transforms of the ten quadrics of Proposition \ref{p120.1} on $\cI_4$ are
  still biregular to $\fP^1\times \fP^1$, while the proper transforms of
  each of the lines turns out to be a {\em Kummer modular surface}, that
  is, $\fP^2$ blown up in four points. Let $\hat{\cI}_4$ be the proper
  transform of $\cI_4$ in $\hat{\hat{\fP}}^4$; then
  ${\grr_2}_{|\hat{\cI}_4}:\hat{\cI}_4\lra \cI_4$ is a desingularisation of
  $\cI_4$.
\end{itemize}
\begin{theorem}\label{t123.1} The varieties $\hat{\cS}_3$ and $\hat{\cI}_4$
  are biregular, and the explicit birational map $\gff:\cS_3- - \ra \cI_4$
  is the birational morphism completing the following diagram:
  $$\begin{array}{ccccc} & \hat{\cS}_3 & \stackrel{\tilde{\gff}}{\lra} &
    \hat{\cI}_4 & \\ \grr_1 & \downarrow & & \downarrow & \grr_2 \\ & \cS_3
    & \stackrel{\gff}{- - \ra} & \cI_4. & \end{array}$$ Moreover, $\gff$ is
  $\gS_6$-equivariant.
\end{theorem}
{\bf Proof:} As $\grr_1$ and $\grr_2$ are $\gS_6$-equivariant, the second
statement follows from the first. Let $D\inn \cS_3$ be the open set:
\begin{equation}\label{e123.1} D=\cS_3-\{\hbox{15 hyperplanes $P_{\gs}$ of
    (\ref{e112.1})}\};
\end{equation}
here we may take the regular map of $D$ onto the set of tangent hyperplanes
(now viewing $\cI_4$ as the projective dual of $\cS_3$), and set
\begin{eqnarray}\label{e123.2} \gff_{|D}:D & \lra & D'\inn \cI_4 \\
  x & \mapsto & (\fP^3)_x=\hbox{tangent hyperplane to $\cS_3$ at $x$}
  \nonumber
\end{eqnarray}
\begin{Lemma} The subset $D'\inn \cI_4$ is: $D'=\cI_4-\{\hbox{10 quadric
    surfaces of Proposition \ref{p120.2}}\}$.
\end{Lemma}
{\bf Proof:} Suppose $x\in D$; then $(\fP^3)_x$ meets $D$ in an irreducible
cubic (the union of the $P_{\gs}$ are {\em all} the linear subspaces
contained in $\cS_3$, so outside of this locus $(\fP^3)_x\cap \cS_3$ cannot
have a linear factor, so, being cubic, must be irreducible), while the ten
quadric surfaces are the locus of the tangent hyperplanes meeting $\cS_3$
in one of the nodes, all of which are excluded in $D$. \ende Now we glue
$D$ onto the rest of $\hat{\cS}_3$, and $D'$ onto the rest of
$\hat{\cI}_4$. The locus $\gL_1=\hat{\cS}_3-D$ coincides with
$\gL_2=\hat{\cI}_4-D'$, as follows from the descriptions of the rational
maps $\grr_1$ and $\grr_2$ above. Both the $\gL_i$ consist of ten
$\fP^1\times \fP^1$'s and 15 rational surfaces, each isomorphic to $\fP^2$
blown up in four points. Hence we can complete $\gff_{|D}$ to a biregular
isomorphism $\gff:\hat{\cS}_3\lra \hat{\cI}_4$, by fixing an isomorphism
$\gff_{\gD}:\gL_1\lra \gL_2$, and setting
$$\gff(x)=\left\{\parbox{6cm}{$\gff_{|D}(x),$ if $x\in D$ \\
    $\gff_{\gD}(x)$, if $x\in \gL_1$}\right.$$ completing the proof of
Theorem \ref{t123.1}. \ende The following description is more concrete. If
$x$ is one of the nodes of $\cS_3$, there is a quadric cone of tangent (to
$\cS_3$) hyperplanes at $x$; so closing up $\gff$ maps $x$ to the quadric
surface over which the above is a cone, i.e., $x$ is blown up. If $x$ is
{\em not} a node, then there is a unique tangent hyperplane $T_x\cS_3$,
determining a point of $\cI_4$.  Furthermore, $T_x\cS_3$ and $T_y\cS_3$
{\em coincide} for $x\neq y$, if and only if $x$ and $y$ are contained in a
common Segre plane (\ref{e112.1}), and the line joining $x$ and $y$ in that
Segre plane passes through one of the four nodes, say $N$, in that Segre
plane. This is because $T_x\cS_3\cap \cS_3=\cH\cup Q_x$, where $Q_x$ is a
residual quadric cone, and the quadric cone is the intersection of
$T_x\cS_3$ with the cone $C_N$ which is the tangent cone of the node $N$ in
the Segre plane. So if $x$ and $y$ lie on a line through $N$, $Q_x$ and
$Q_y$ coincide, so $T_x\cS_3$ and $T_y\cS_3$ coincide also.
\begin{theorem}\label{t122a.1} The duality map $d:\cS_3- - \ra \cI_4$ is
  given by the linear system of quadrics \ref{p112b.1}, i.e., by the
  elements of the ideal $\scI(10)$ of the ten nodes: $d=\gff$.
\end{theorem}
{\bf Proof:} It suffices to check that $d$, viewed as a modification of
$\cS_3$, coincides with the birational map $\gff$ of Theorem \ref{t123.1}.
But this is easy. As the base locus is the set of nodes, these are blown
up.  As just explained, $x$ and $y$ in one of the Segre planes map to the
same point on the image line precisely when the line joining them passes
through one of the nodes in the Segre plane. As these lines are precisely
what the map $\gff$ blows down, $d$ certainly coincides with $\gff$. \ende
We also have the following analogue of Corollary \ref{c112b.1}.
\begin{lemma}\label{l122a.1} The ideal $\scI(15)$ of the 15 singular lines of
  the Igusa quartic coincides with the Jacobian ideal of $\cI_4$.
\end{lemma}
{\bf Proof:} Once again the inclusion $\scJ ac(\cI_4)\inn \scI(15)$ is
obvious, and the inverse inclusion can be verified by means of standard
basis computations, namely that $\scI(15)$ is generated by five cubics.
\ende Along the same lines as Theorem \ref{t122a.1} we then get
\begin{theorem} \label{t122b.1} The duality map $d:\cI_4- - \ra\cS_3$ is given
  by the system of cubics containing the 15 lines, i.e., by the Jacobian
  ideal of $\cI_4$: $d=\gff^{-1}$.
\end{theorem}
{\bf Proof:} As above, it suffices to show that $d$, viewed as a
modification of $\cI_4$, coincides with the map $\gff^{-1}$ of Theorem
\ref{t123.1}. This is readily verified, as the base locus, the 15 lines,
are blown up, while the tangent planes for any two points $x$ and $y$ in a
common quadric of $\cI_4$ (of Proposition \ref{p120.1}) coincide, blowing
down the quadric surface to a node. \ende

\subsection{The Siegel modular threefold of level 4}
{}From the general theory of congruence subgroups, $X_{\gG(4)}\lra
X_{\gG(2)}$ is a Galois cover, with Galois group $\gG(2)/\gG(4)\isom
(\fZ/2\fZ)^9$.  Indentifying $X_{\gG(2)}^*$ with $\cI_4$ and identifying
$\cI_4$ birationally with $\cS_3$, we can consider Fermat covers over
$X_{\gG(2)}^*$, i.e., given by a diagram
\begin{equation}\label{e122b.1} \begin{array}{cccccc} Z({\bf A_4},n) &
    \stackrel{\-{\gff}^{-1}}{- - \ra} & Y^{\wedge}({\bf A_4},n) & \lla &
    \tilde{Y}({\bf A_4},n) & \\ \downarrow & & \downarrow & & \downarrow &
    (\fZ/2\fZ)^9 \\ \cI_4 & \stackrel{\gff^{-1}}{- - \ra} & \cS_3 & \lla &
    \tilde{\cS}_3 &
\end{array}
\end{equation}
where $\-{\gff}^{-1}$ is {\em induced} by $\gff^{-1}$, that is,
(\ref{e122b.1}) is a fibre square (cf. (\ref{e115.0}), where
$\tilde{\cS}_3$ is denoted $\tilde{\fP}^3$).
\begin{theorem}\label{t122b.1} The Fermat cover $Z({\bf A_4},2)$ is the
  Satake compactification of the Siegel modular threefold of level 4.
\end{theorem}
{\bf Proof:} It suffices to show that $\~Y({\bf A_4},2)$ is the induced
cover over $\hat{\cI}_4$, where $\grr_2:\hat{\cI}_4\lra \cI_4$ is the
desingularisation of $\cI_4$ of Theorem \ref{t123.1}. Now the
identification can be reduced to identifying what is in the branch locus of
$\~Y({\bf A_4},2)\lra \~{\cS}_3$. There are two kinds of components:
\begin{itemize}\item[a)] covers $\~Y({\bf A_3},2)$ of blown up $\fP^2$'s, the
  $H_{ij}$ of (\ref{e111b.1});
\item[b)] covers $\~Y({\bf A_2},2)\times \~Y({\bf A_2},2)$ of $\fP^1\times
  \fP^1$'s, the $L_{0ij}$.
\end{itemize}
\begin{Lemma}\label{l122b.1} $\~Y({\bf A_3},2)\isom S(4)$, Shioda's elliptic
  modular surface of level 4.
\end{Lemma}
{\bf Proof:} This is well-known. $\~Y({\bf A_3},2)$ is K3 since it is a
Fermat cover branched over six lines. One constructs structures of fibre
space $\~Y({\bf A_3},2)\lra \fP^1$ with elliptic curves as fibres by taking
the cover of the pencil of lines through a node (each such line meets four
of the six lines outside the node, so the cover is branched at four points,
i.e., is elliptic). The six fibres of type $I_4$ are readily identified, as
are the 16 sections. \ende
\begin{Lemma}\label{l122c.1} The cover $\~Y({\bf A_2},2)\lra \fP^1$ coincides
  with the cover $(\gG_1(4)\bs \fS_1)^*\lra \fP^1$, by which we mean the
  Galois actions coincide.
\end{Lemma}
{\bf Proof:} This is even more well-known. \ende The theorem now follows,
provided we accept that $\~Y({\bf A_4},2)$ is a quotient of $\fS_2$ at all,
i.e., that the cover $\fS_2\lra \cI_4^0$ factorises (here
$\cI_4^0=\cI_4-\{\hbox{15 lines}\}$), $\tilde{Y}({\bf A_4},2)^0:= \~Y({\bf
  A_4},2)-q^{-1}(\hbox{15 lines})$:
\begin{equation}\label{e122c.1}\unitlength1.6cm \begin{picture}(1.6,1.4)
    \put(.2,1.33){$\fS_2$} \put(.6,1.4){\vector(1,0){1}}
    \put(1.8,1.33){$\cI_4^0$} \put(.4,1.25){\vector(1,-1){.6}}
    \put(.7,.45){$\~Y({\bf A_4},2)^0.$} \put(1.2,.66){\vector(1,1){.6}}
    \put(1.7,.8){$q$}
\end{picture}
\end{equation}
But there is an easy way to see that this is the case: we can, for any
given $x\in \~Y({\bf A_4},2)-\gD$ and $y=q(x)$, put a level 4 structure on
$A_y$, such that the Galois group just permutes the level 4 over level 2
structures, that is, we make the identification $\gG(2)\bs\gG(4)\isom
(\fZ/2\fZ)^9\isom$ the Galois group of the cover.
So $\~Y({\bf A_4},2)$, being a moduli space as in Shimura's
theory, is a quotient of $\fS_2$. \ende

One could also imagine arguing with uniqueness of Galois covers, since we
know the branch locus, branch degrees and Galois group. However there is in
general no such uniqueness of covers, so we have to be careful. In our
situation, there are two possible approaches to show uniqueness:
\begin{itemize}\item[1)] Since the modular subvarieties determine, on the
  group-theoretic side, generators of the corresponding arithmetic group,
  we could conclude, from the isomorphisms \ref{l122b.1} and \ref{l122c.1},
  the desired result.
\item[2)] Since the branch divisors are totally geodesic with respect to
  the Bergmann metric, on the cover the metric retains its symmetry
  property.
\end{itemize}
Method 1) has been applied in \cite{J}, and 2) can be carried out for ball
quotients.

\section{The Hessian varieties of ${\cal  S}_3$ and ${\cal I}_4$}
\subsection{The Nieto quintic}\label{section4.1}
Let $(x_0:\ldots:x_5)$ be the projective coordinates on $\fP^5$ used to
define $\cS_3$ in (\ref{e111.3}), and let $\gs_i=\gs_i(x_0,\ldots,x_5)$ be
the $i$-th elementary symmetric function $\gs_{\gl}=\sum_{i_1<\ldots <
  i_{\gl}}x_{i_1}\cdots x_{i_{\gl}}$ in $(x_0:\ldots:x_5)$. Define the {\em
  Nieto quintic} $\cN_5$ by the equations\begin{equation}\label{e124.1}
  \cN_5=\left\{\begin{array}{l}\gs_1 = 0 \\ \gs_5=0\end{array} \right. \inn
  \fP^4=\{\gs_1=0\}\inn \fP^5.
\end{equation}
The symmetry of $\cN_5$ under the symmetric group $\gS_6$ is evident from
the equation. This quintic was discovered in the thesis \cite{Nie} and
further studied in \cite{BN}, which will be our general reference for this
section. We just briefly describe the geometry of $\cN_5$ without
discussing details.

The singular locus is relatively easy to determine, just by calculating the
Jacobian of (\ref{e124.1}). The result is
\begin{proposition}[\cite{BN}, 3.1]\label{p124.1} $\cN_5$ has the following
  singular locus:
\begin{itemize}\item[(i)] 20 lines $L_{ijk}=\{x_i=x_j=x_k=0=\sum x_i\}$;
\item[(ii)] ten isolated points, the $\gS_6$-orbit of $(1,1,1,-1,-1,-1)$,
  which are the points $P_{ij}=(1,\pm1,\ldots \pm 1)$, with $+1$ in the
  $i$-th and $j$-th positions.
\end{itemize}
\end{proposition}
We will give a different proof of this below, see the discussion following
Proposition \ref{piq8.1}.  Note that the ten points occuring in (ii) are
just the ten nodes of $\cS_3$ (see (\ref{e112.0}), cf.~also
Remark \ref{r126.2} below).  Furthermore, a local
calculation shows that the singularities of $\cN_5$ along the lines of (i)
are of the type $\{\hbox{disc}\}\times A_1$, and at the points of (ii) are
ordinary double points. Hence the former are resolved by a $\fP^1$-bundle
over the line $L_{ijk}$, while the points are resolved, as with the case of
$\cS_3$, by quadric surfaces. The 20 lines $L_{ijk}$ of \ref{p124.1} meet
at the following 15 points:
\begin{equation}\label{e124.2} Q_{ij}=(0,\ldots,1,\ldots,-1,\ldots)
  =\{\gS_6-\hbox{orbit of } Q_{56}=(0:0:0:0:1:-1)\}.
\end{equation}
\begin{lemma}\label{l124.1}
  The 20 lines $L_{ijk}$ of Proposition \ref{p124.1} meet four at a time at
  the 15 points $Q_{ij}$; each line $L_{ijk}$ contains three of the points,
  namely we have $Q_{ij}\in L_{klm}\iff \{i,j\}\cap \{k,l,m\}=\emptyset$.
\end{lemma}
{\bf Proof:} The line $L_{123}$ contains the three points $Q_{46},\ Q_{45}$
and $Q_{56}$, so by $\gS_6$-invariance each line contains three of the
$Q_{ij}$. The point $Q_{56}$ is contained in the four lines $L_{123},\
L_{124},\ L_{134}$ and $L_{234}$, so by $\gS_6$-invariance, each point is
contained in four lines.  \ende Also, $\cN_5$ contains a finite number of
linear planes.
\begin{lemma}\label{l125.1} $\cN_5$ contains the following 30 $\fP^2$'s:
\begin{itemize}\item[(i)] 15 planes $N_{ijkl}=\{x_i+x_j=x_k+x_l=x_m+x_n=0\}$;
\item[(ii)] 15 planes $N_{ij}=\{x_i=x_j=0=\sum_{k\neq i,j}x_k\}$;
\end{itemize}
\end{lemma}
{\bf Proof:} It is immediately verified that these planes satisfy the
equation (\ref{e124.1}).\ende Presumably these are in fact {\em all} the
linear planes contained in $\cN_5$. Note that the $N_{ijkl}$ are just the
15 planes (\ref{e112.1}) lying on the Segre cubic $\cS_3$.

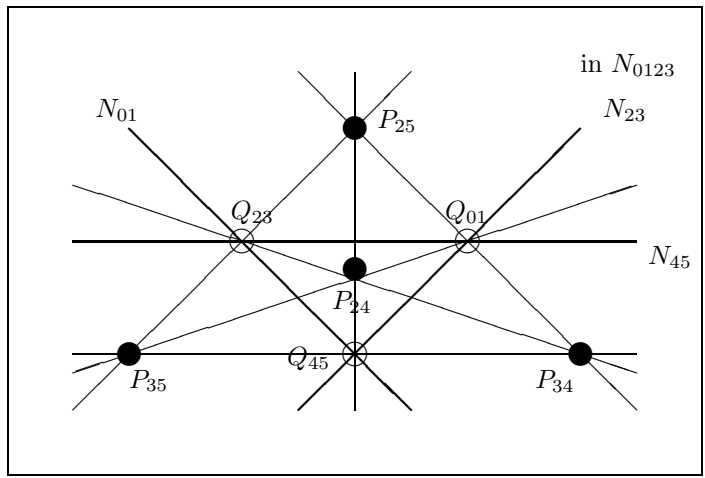
\begin{figure}[thb]
  $$\fbox{\unitlength1.5cm \begin{picture}(6,4) \put(.5,.5){\line(1,1){3}}
      \put(1,1){\line(3,1){4.5}}\put(1,1){\line(-3,-1){.5}}
      \put(1,1){\circle*{.2}}\put(1,.7){$P_{35}$}
      \put(.5,2.5){\line(3,-1){5}} \put(.7,3.1){$N_{01}$}
      \put(2,2){\circle{.2}}\put(1.9,2.2){$Q_{23}$}
      \put(3,1){\circle{.2}}\put(2.4,.9){$Q_{45}$}
      \put(2.5,3.5){\line(1,-1){3}}\put(3,3.5){\line(0,-1){3}}
      \put(3,1.75){\circle*{.2}} \put(2.8,1,4){$P_{24}$}
      \put(3,3){\circle*{.2}}\put(3.2,3){$P_{25}$}
      \put(4,2){\circle{.2}}\put(3.8,2.2){$Q_{01}$}
      \put(5,1){\circle*{.2}}\put(4.6,.7){$P_{34}$}
      \put(5.2,3.1){$N_{23}$}\put(5.6,1.8){$N_{45}$} \put(5,3.5){in
        $N_{0123}$} \put(.5,1){\line(1,0){5}} \thicklines
      \put(1,3){\line(1,-1){2.5}} \put(2.5,.5){\line(1,1){2.5}}
      \put(.5,2){\line(1,0){5}}
\end{picture}}$$
\caption[15 planes of type 1 on ${\cal N}_5$]{\label{Figure1}\small
The plane $N_{0123}$}
\end{figure}

Among the 15 planes $N_{ijkl}$ the common intersections were described in
the discussion of the planes $P_{\gs}$ on the Segre cubic (see
(\ref{e111.1})).
\begin{lemma}\label{l125.2} Each plane $N_{ijkl}$ contains the following four
of
  the ten points of \ref{p124.1}, (ii): $$P_{km},\ P_{kn},\ P_{lm} \hbox{
    and } P_{ln};$$ it also contains the following three of the 15 points
  $Q_{ij}$ of (\ref{e124.2}): $Q_{ij},\ Q_{kl}$ and $Q_{mn}$.
\end{lemma}
{\bf Proof:} Consider $N_{0123}$; it contains the four nodes
$(1:-1:1:-1:1:-1),\ (1:-1:-1:1:1:-1),\ (1:-1:1:-1:-1:1)$ and
$(1:-1:-1:1:-1:1)$ which are the points $P_{24},\ P_{25},\ P_{34}$ and
$P_{35}$, which gives the first statement by $\gS_6$-symmetry (there is an
asymmetry in the notation, since we may take $i<j,k<l$ in the notation for
$N_{ijkl}$, and since the first coordinate of $P_{ij}$ may be assumed to be
$+1$). Similarly, $N_{0123}$ contains the three points $Q_{01},\ Q_{23}$
and $Q_{45}$, giving the second statement by $\gS_6$-symmetry.\ende We now
note that these seven points lie in the plane $N_{0123}$ as in Figure
\ref{Figure1}.
This is in fact easily checked. Note that the lines in $N_{0123}$, i.e.,
the intersections with the other $N_{ijkl}$ are {\em not} the lines of
Proposition \ref{p124.1}; those lines have equations such as
$x_0=x_1+x_2=x_1+x_3=x_4+x_5$. However, in the 15 planes $N_{ij}$ of
\ref{l125.1}, several of the 20 singular lines $L_{ijk}$ {\em do} lie. In
fact, we have
\begin{lemma}\label{l125.3} Each $N_{ij}$ contains the four lines $L_{ijk}.\
  L_{ijl},\ L_{ijm}$ and $L_{ijn}$. There are three planes passing through
  $L_{ijk}$, namely $N_{ij},\ N_{ik}$ and $N_{jk}$. $N_{ij}$ contains none
  of the nodes of \ref{p124.1}, (i), but contains six of the points
  $Q_{ij}$ of (\ref{e124.2}), namely $Q_{kl},\ Q_{km},\ Q_{kn},\ Q_{lm},\
  Q_{ln}$ and $Q_{mn}$. These six points lie three at a time on the
  $L_{ijk}$ and form in each $N_{ij}$ a configuration as shown in Figure
  \ref{Figure2}.
\begin{figure}[tbh]
  $$\fbox{\unitlength1.5cm \begin{picture}(6,4)\put(.5,1){\line(1,0){5}}
      \put(.5,2){\line(1,0){5}} \put(3,3.5){\line(0,-1){3}}
      \put(3,3){\circle{.2}}\put(1,1){\circle{.2}}\put(2,2){\circle{.2}}
      \put(3,1.66){\circle{.2}}\put(4,2){\circle{.2}}\put(5,1){\circle{.2}}
      \thicklines \put(.5,.5){\line(1,1){3}}\put(5.5,2.5){\line(-3,-1){5}}
      \put(.5,2.5){\line(3,-1){5}}\put(2.5,3.5){\line(1,-1){3}}
      \put(.5,2.7){$L_{015}$}\put(1.6,3.6){$L_{012}$}\put(3.6,3.6){$L_{013}$}
      \put(5,2.7){$L_{014}$}\put(5.5,1.2){$N_{0125}$}
      \put(5.5,1.8){$N_{0124}$}\put(2.6,.2){$N_{0123}$}
      \put(1.1,.7){$Q_{25}$}\put(4.6,.7){$Q_{34}$}\put(2.6,1.3){$Q_{23}$}
      \put(1.8,2.2){$Q_{24}$}\put(3.8,2.2){$Q_{35}$}\put(3.2,2.8){$Q_{45}$}
      \put(5,3.5){in $N_{01}$}
\end{picture}}$$
\caption[15 planes of type 2 on ${\cal N}_5$]{\label{Figure2}\small
  The plane $N_{01}$}
\end{figure}
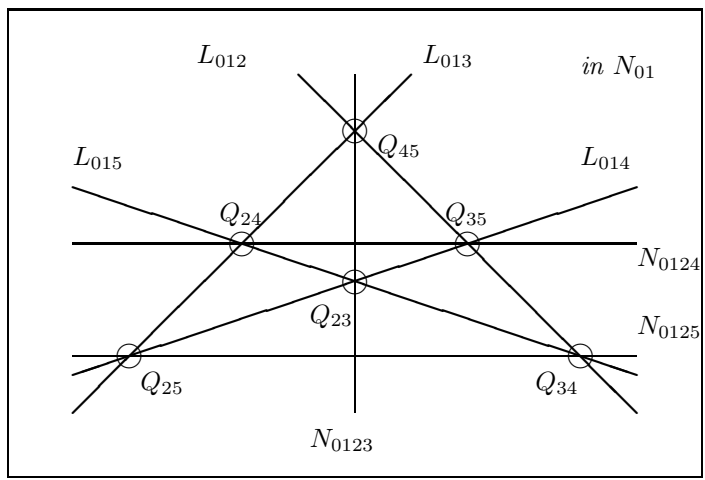
The three light lines are intersection of $N_{01}$ with $N_{ijkl}$ as
indicated.
\end{lemma}
{\bf Proof:} This is once again easily verified. \ende Finally we note that
there are hyperplanes in $\fP^4$ cutting out these $\fP^2$'s on $\cN_5$.
\begin{lemma}\label{l126.1} The six hyperplanes $\cH_{ij}=\{x_i+x_j=0\}$ meet
  $\cN_5$ each in the union of the three planes $N_{ijkl},\ N_{ijkm}$ and
  $N_{ijkn}$, and a residual quadric; the six hyperplanes $x_i=0$ meet
  $\cN_5$ each in the union of five planes $N_{ij},\ N_{ik},\ N_{il},\
  N_{im}$ and $N_{in}$.
\end{lemma}
{\bf Proof:} This is once again just a computation. \ende Now let us
consider the intersection of $\cS_3$ and $\cN_5$. As is obvious from the
above description, they both contain the 15 planes $N_{ijkl}$, and, the
intersection being of degree 15, this is the entire intersection.  From
general arguments on projective varieties, from the fact that the dual of
$\cS_3$, namely the Igusa quartic $\cI_4$, is {\em normal}, it follows that
the parabolic divisor on $\cS_3$, which is the intersection of $\cS_3$ with
the {\em Hessian variety}, must get blown down under the duality map, i.e.,
the intersection $Hess(\cS_3)\cap \cS_3$ consists of the 15 planes on
$\cS_3$! Since the Hessian has degree 5, this is the entire intersection,
and it is natural to ask whether $\cN_5$ and $Hess(\cS_3)$ are related. In
fact, we have
\begin{lemma}\label{l126.2} The Nieto quintic is the Hessian of $\cS_3$,
  i.e., $\cN_5= Hess(\cS_3)$, with equality, not just isomorphism.
\end{lemma}
{\bf Proof:} This is an easy computation (at least for a computer).\ende
\begin{remark}\label{r126.2}
 Since the Hessian variety $\Hess(V)$ aquires nodes where $V$
  has nodes, this ``explains'' the ten isolated singularities on $\cN_5$.
\end{remark}

\subsection{Two birational transformations}
We consider in this section two particularly interesting birational maps
from $\cN_5$.
\subsubsection{ \ }
The first arises through the duality map.  Consider the birational map of
$\fP^4$ given in the following diagram:\begin{equation}\label{e126.1}
\begin{array}{ccccccc} \fP^4 & \stackrel{\grr_1}{\lla}
  & \hat{\fP}^4 & \stackrel{\~{\ga}}{- - - \ra} & \hat{\hat{\fP}}^4 &
  \stackrel{\grr_2}{\lra} & \fP^4 \\ \cup & & \cup & & \cup & & \cup \\
  {\cal S}_3 & \lla & \hat{\cal S}_3 & \stackrel{\sim}{\lra} &
  \hat{\cal I}_4 & \lra & {\cal I}_4
\end{array}\quad;
\end{equation}
$\grr_1$ and $\grr_2$ were described in section \ref{section3.4}, and this
diagram extends the one of Theorem \ref{t123.1} to the ambient rational
fourfolds. It is easily seen that the ensuing rational map of $\fP^4$,
 $\ga:=\grr_2\circ \~{\ga}\circ \grr_1^{-1}:\fP^4\lra \fP^4$, is the
map given by the Jacobian ideal of $\cS_3$, that is, by the
linear system of quadrics on the ten nodes of $\cS_3$ (see Corollary
\ref{c112b.1}). This ``defines'' the map $\tilde{\ga}$; although we could
in principle take any extension of $\tilde{\gff}:\hat{\cS_3}\lra
\hat{\cI_4}$, for our purposes it is convenient to use $\scI(10)$. Then set
\begin{equation}\label{e127.1} \cW_{10}:=\ga(\cN_5)
  \inn \fP^4.
\end{equation}
We now describe $\cW_{10}$ and show it is a hypersurface of degree 10,
explaining the notation. First we have the
\begin{lemma}\label{l25aux} The map $\ga:\fP^4\lra \fP^4$ blows up the ten
  nodes $P_{ij}$ of $\cS_3$, with exceptional divisors $\cE_{ij}$. Let
  $\cC_{ij}$ denote the tangent cone at the point $P_{ij}$, a quadric cone
  fibred in lines passing through $P_{ij}$. Then each line of the cone gets
  blown down to the corresponding point in $\cE_{ij}$.
\end{lemma}
{\bf Proof:} Since all the quadrics of $\scI(10)$
vanish at the $P_{ij}$, these points
are blown up. To say the lines of $\cC_{ij}$ get blown down is to say the
ratios of the quadrics are constant along the line. This follows from the
fact that the quadrics are the partial derivatives of $f$ (the defining
polynomial of $\cS_3$), and the line is tangent to the zero locus of $f$.
\ende From this we get
\begin{lemma} $\ga$, restricted to $\cN_5$,
  is an isomorphism on the complement in $\grr_1^{-1}(\cN_5)$ of the
  intersection of $\cN_5$ with the tangent cones at the ten isolated
  singularities $P_{ij}$.
\end{lemma}
{\bf Proof:} This is clear from construction, taking into account the
following fact, proved in \cite{BN}: the intersection of $\cN_5$ with the
tangent cone of one of the nodes $P_{ij}$ consists of the six Segre planes
through the node, and an irreducible quartic {\em ruled} surface. It is
then clear that these ruled surfaces get blown down, and that outside the
Segre cubic and the ruled quartics, the birational map is a morphism. \ende
{}From this it follows:
\begin{lemma}\label{l127.1} The birational map $\grr_2\circ \~{\ga} \circ
  \grr_1^{-1}:\cN_5- - \ra \cW_{10}$ has image $\cW_{10}$, whose singular
  locus contains the following.
\begin{itemize}\item[(i)] ten singular quadric surfaces (the tangent
  hyperplane intersections of $\cI_4$);
\item[(ii)] 20 singular lines, coming from the singular locus of $\cN_5$.
\end{itemize}
\end{lemma}
{\bf Proof:} This follows from the description above. \ende Now note that
both $Hess(\cI_4)$ and $\cW_{10}$ are symmetric under $\gS_6$, and both are
singular along the ten quadrics of the Igusa $\cI_4$. It may very well be
that the two coincide, but we have not checked this. At any rate,
$\cW_{10}$ meets $\cI_4$ in the union of quadric surfaces, each with
multiplicity 2. Hence
\begin{theorem}\label{t127.1} $\cI_4\cap\cW_{10}$ consists of the ten quadric
  surfaces (\ref{e117a.1}), each with multiplicity 2. Consequently, the
  degree of $\cW_{10}$ is 10, justifying the notation.
\end{theorem}
{\bf Proof:} The intersection has reduced degree 20, and each surface
component is counted twice, hence the degree of the intersection is 40, so
the degree of $\cW_{10}$ is 10. \ende
\begin{problem}\label{p131.1}
  Is $\cW_{10}$ also a compactification of an arithmetic quotient?
\end{problem}

\subsubsection{ \ }
The other birational transformation is the following.
\begin{equation}\label{e127.2}\begin{minipage}{14cm}\begin{itemize}
    \item[a)] Blow up the 15 points $Q_{ij}$ of (\ref{e124.2}); let
      $p_1:\~{\cN}_5\lra \cN_5$ denote this blow up.
    \item[b)] As each of the lines $L_{ijk}$ contains three points (see
      Lemma \ref{l124.1}), each $L_{ijk}$ can be blown down to an isolated
      singular point (the normal bundle is $\cO(-2)\oplus\cO(-2)$,
      cf.~(\ref{e110.1})). Let $p_2:\~{\cN}_5\lra
      \hat{\cN}_5$ denote this blow down.
\end{itemize}
\end{minipage}
\end{equation}
The following is easy to see (see Figures \ref{Figure1} and \ref{Figure2}).
\begin{lemma}\label{l127.3} The singular locus of $\hat{\cN}_5$ consists of
  the 20 isolated cusps from (\ref{e127.2}) b), and the ten cusps, the images
  of the singular points $P_{ij}$ of Proposition \ref{p124.1}, (ii).  The
  proper transforms of the $N_{ijkl}$ of Lemma \ref{l125.1}, (i) on $\cN_5$
  are $\fP^2$'s blown up in three points, a del Pezzo surface; the proper
  transforms of the $N_{ij}$ of Lemma \ref{l125.1}, (ii) are $\fP^2$'s
  blown up in six points, then the $L_{ijk}$ are blown down to four nodes,
  so this is the singular cubic surface with four nodes, the Cayley cubic.
\end{lemma}

\subsection{Moduli interpretation}
The Nieto quintic was discovered as the solution of a certain moduli
problem, and we briefly state the results of \cite{BN} describing this.

The point of departure is the action of the Heisenberg group $H_{2,2}$ on
$\fP^3$, and the study of quartics which are invariant under the action.
$H_{2,2}$ is a group of order 32 generated by the following linear
transformations of $\fP^3$ with coordinates $(z_0:z_1:z_2:z_3)$:
\begin{equation}\label{e128.1}\begin{array}{rclcl} \gs_1 & : &
    (z_0:z_1:z_2:z_3) & \mapsto & (z_2:z_3:z_0:z_1) \\ \gs_2 & : &
    (z_0:z_1:z_2:z_3) & \mapsto & (z_1:z_0:z_3:z_2) \\ \tau_1 & : &
    (z_0:z_1:z_2:z_3) & \mapsto & (z_0:z_1:-z_2:-z_3) \\ \tau_2 & : &
    (z_0:z_1:z_2:z_3) & \mapsto & (z_0:-z_1:z_2:-z_3) \end{array}
\end{equation}
The center of the group is $\pm1$ and $PH_{2,2}=H_{2,2}/\pm1$ has a nice
interpretation:
\begin{equation}\label{e128.2} PH_{2,2}\isom (\fZ/2\fZ)^4,
\end{equation}
which carries, as in (\ref{e121.1}), an induced symplectic form.  This
means that one can speak of isotropic elements of the {\em group}
$PH_{2,2}$. The normaliser of $H_{2,2}$ in $SL(4,\fC)$ maps surjectively to
$\gS_6\isom Sp(4,\fZ/2\fZ)$, which acts transitively on diverse geometric
loci of the sympectic form {\em inside} the group $PH_{2,2}$. These loci
are:
\begin{equation}\label{e128.3}\begin{minipage}{14cm}\begin{itemize}
    \item[a)] 15 pairs of skew lines
    \item[b)] 15 invariant tetrahedra
    \item[c)] ten fundamental quadrics.
\end{itemize}\end{minipage}\end{equation}

The moduli problem considered is a special set of quartics which are
invariant under (\ref{e128.2}). The set of {\em all} invariant quartics is
just a $\fP^4$, spanned for example by the five quartics:
$$\begin{array}{ccc} & g_0:=z_0^4+z_1^4+z_2^4+z_3^4 & \\
  g_1:=2(z_0^2z_1^2+z_2^2z_3^2) & g_2=2(z_0^2z_2^2+z_1^2z_3^2) &
  g_3:=2(z_0^2z_3^2+z_1^2z_2^2) \\ & g_4:=4z_0z_1z_2z_3. & \end{array}$$
Let $(A,B,C,D,E)$ denote the coordinates of a particular quartic
$Q_{(A,B,C,D,E)}=\{Ag_0+Bg_1+Cg_2+Dg_3+Eg_4=0\}$. The generic quartic
$Q_{(A,B,C,D,E)}$ is smooth, and the locus of singular quartics can be
determined as an equation in $(A,\ldots,E)$. Note that the $(A,\ldots,E)$
are functions of $(z_0:\cdots:z_3)$, so the answer as to whether
$Q_{(A,B,C,D,E)}$ is singular depends on the point $z\in \fP^3$. This is
discussed in detail in \cite{BN}. The result is given in Table
\ref{table19}.
\begin{table}
\begin{center}
\caption{\label{table19}Singular Heisenberg invariant quartics }
\begin{tabular}{|l|l|l|l|}\hline $z\in \fP^3$ & dim$Q^{sing}$ &
  $Q_{(A,B,C,D,E)}$ & $S_{(A,B,C,D,E)}$ \\ \hline\hline $\notin$ fix line &
  0 & Kummer surface & Segre cubic \\ \hline $\in$ one fix line & 2 &
  singular in four coordinate vertices & $A=0$ \\ \hline $\in$ the
  intersection of two fixed lines & 3 & singular along two fixed lines &
  $A=B=0$ \\ \hline
\end{tabular}
\end{center}
{\small Notations: $Q^{sing}$ denotes the space of quartics singular at
  $z$, $S_{(A,B,C,D,E)}$ denotes the equation of the locus $Q^{sing}$ in
  the coordinates $(A,B,C,D,E)$.}
\end{table}
As one sees, the first row of the table is equivalent to Theorem \ref{t119.1}
above!  The special class of quartics to be considered here is, however, a
quite different set, consisting of generically smooth quartics. This is the
set of Kummer surfaces of (1,3)-polarised abelian surfaces, which, as it
turns out, can be smoothly embedded in $\fP^3$. This was discovered
independently by Naruki and Nieto (see \cite{na} and \cite{Nie}). The 16
exceptional $\fP^1$'s resolving the 16 double points of the Kummer surface
are 16 disjoint {\em lines} on the quartics. Also, by a result of Nikulin
\cite{Ni}, the converse is true, i.e., any quartic containing 16 lines is a
Kummer surface. Furthermore, the quartic being invariant under $PH_{2,2}$,
if it contains one line, it contains all 16 transforms, so the moduli
involved is the condition:
\begin{equation}\label{e129.1} L\ \parbox[t]{12cm}{ is a line in $\fP^3$
    lying on a smooth Heisenberg invariant quartic surface}
\end{equation}
The equation describing this in the Grassmannian $\fG(2,2)=\{x_0^2+\cdots +
x_5^2=0\}$ is calculated in \cite{Nie}. It is
\begin{equation}\label{e130.1}
  \cM_{20}=\{\gs_5(x_0^2,\ldots,x_5^2)=0=\gs_1(x_0^2,\ldots,x_5^2)\}.
\end{equation}
Now one considers the natural 2-power map
\begin{eqnarray}\label{e130.2} m_2:\fP^5 & \lra & {\Bbb P}^5 \\
  (x_0,\ldots,x_5) & \mapsto & (x_0^2,\ldots,x_5^2)=(u_0,\ldots,u_5)
  \nonumber
\end{eqnarray}
and the image of $\cM_{20}$ in $\fP^5$. Comparing the equations
(\ref{e124.1}) and (\ref{e130.1}) we have
\begin{lemma}\label{l130.1} $m_2(\cM_{20})=\cN_5$.
\end{lemma}
The main results of \cite{BN} can be described as follows. First we define
a Zariski open subset $M^s\inn \cM_{20}$. The following 15 quadric surfaces
$q_{ij}\inn \fP^5$ actually lie on $\cM_{20}$, as is easily verified:
\begin{equation}\label{e130a.1} q_{ij}=\{x_i=x_j=0=\sum_{m\neq i,j}x_m^2\}.
\end{equation}
Under the squaring map $m_2$ (\ref{e130.2}) the quadric $q_{ij}$ maps to
the plane
\begin{equation} N_{ij}=\{u_i=u_j=0=\sum_{m\neq i,j}u_m\};
\end{equation}
so the image of ${\bf Q}:=\cup_{i,j}q_{ij}$ is ${\bf N}:=\cup_{i,j}N_{ij}$,
and the planes $N_{ij}$ are the 15 planes of Lemma \ref{l125.1} (ii).
Furthermore the $N_{ij}$ are contained in the branch locus of
${m_2}_{|\cM_{20}}:\cM_{20}\lra \cN_5$; this locus is {\em singular} on
$\cM_{20}$ because $\cN_5$ is tangent to $u_i=0$ and $u_j=0$ in all of
$N_{ij}$.

Next consider the inverse image under $m_2$ of the ten nodes; since the
nodes lie on {\em none} of the branch planes $u_i=0$, each node has
$\deg(m_2)=32$ inverse images, so $\cM_{20}$ has 320 singular points
(clearly also nodes), which are the $\gS_6$-orbit, call it ${\bf P}$, of
the points
\begin{equation}\label{e130a.2} (\pm1:\pm1:\pm1:\pm i:\pm i: \pm i).
\end{equation}
Finally consider the inverse images of the 15 Segre planes of Lemma
\ref{l125.1} (i). This locus is given by the 15 equations which are the
$\gS_6$-orbit of
\begin{equation}\label{e130a.3} x_0^2+x_1^2=x_2^2+x_3^2=x_4^2+x_5^2=0.
\end{equation}
Inspection shows that this degree 8 surface on $\cM_{20}$ splits into eight
planes, giving altogether 120=15.8 planes on $\cM_{20}$; let ${\bf R}$
denote their union. Now define:
\begin{equation}\label{e130a.5} M^s:=\cM_{20}-{\bf Q}-{\bf P}-{\bf N},\
  \cN_5^s:=m_2(M^s).
\end{equation}
Then the statement proved in \cite{BN} is
\begin{theorem}\label{t130.1} \begin{itemize}\item[a)] $M^s$ is isomorphic to
    a Zariski open subset of the moduli space $\cA_{(1,3)}(2)$ of abelian
    surfaces with a $(1,3)$ polarisation and a level $2$ structure;
  \item[b)] There is a double cover $p:\~{\cN}\lra \cN_5$ for which
    $p^{-1}(\cN_5^s)$ is isomorphic to a Zariski open set of the moduli
    space $\cA_{(2,6)}(2)$;
  \item[c)] $\cN_5^s$ is the moduli space of $PH_{2,2}$-invariant smooth
    quartic surfaces containing $16$ skew lines.
\end{itemize}
\end{theorem}
Since the varieties $\cM_{20}$, $\~{\cN}$ and $\cN_5$ are compactifications
of the Zariski open sets of (\ref{e130a.5}), we have the following:
\begin{corollary} \label{c130.1} There are birational equivalences:
  $$\cM_{20}- - \ra (\gG_{(1,3)}(2)\bs \fS_2)^*,\quad \~{\cN}- - \ra
  (\gG_{(2,6)}(2)\bs \fS_2)^*,\quad \cN_5- - \ra (\gG\bs \fS_2)^*,$$ where
  $ \gG_{(1,3)} \inn\gG_{(2,6)} \inn \gG,\
  [\gG:\gG_{(2,6)}(2)]=2.$
\end{corollary}
As is shown in \cite{BN}, the map $\~{\cN}\lra \cN_5$ is given in the
following way. It just happens to turn out the {\em any} of the
$PH_{2,2}$-invariant quartics with 16 skew lines actually contains 32
lines, the first skew set of 16 and a second set of 16 skew lines. The
second set of sixteen is found as the image of the first set under the
involution
\begin{equation}\label{e130.3} (x_0:\ldots:x_5)\mapsto \left({-1 \over x_0}:
    {1\over x_1}: \cdots :{1\over x_5}\right),
\end{equation}
which can be adjoined to the group $PH_{2,2}$ to form a group of order 32.
Altogether the 32 lines have the following properties.
\begin{equation}\label{e130.4} \begin{minipage}{14cm}\begin{itemize}
    \item[a)] The 32 lines intersect in 32 points;
    \item[b)] Each line contains ten of the 32 intersection points;
    \item[c)] Each intersection point is contained in ten of the 32 lines.
\end{itemize}\end{minipage}\end{equation}
A configuration with the properties (\ref{e130.4}) is called a
$(32_{10})$-configuration.

{}From Nikulin's results just mentioned, it follows that the second set of 16
lines are also the images of blown-up torsion points on another abelian
surface, so there are {\em two} abelian surfaces with $(2,6)$ polarisation
and level $2$ structure giving rise to the {\em same} resolved Kummer
surface, i.e., the map is given by
\begin{eqnarray}\label{e130.5} \~{\cN} & \lra & {\cal N}_5 \\
  (A_{\tau_1},A_{\tau_2}) & \mapsto & \-{(A_{\tau_1}/\{\pm 1\})}\isom
  \-{(A_{\tau_2}/\{\pm 1\})},\nonumber
\end{eqnarray}
where the isomorphism permutes the two sets of 16 skew lines.

The next step is to identify the modular subvarieties on the arithmetic
quotients of Corollary \ref{c130.1}. From the structure of the periods we
know that in terms of abelian surfaces, these modular subvarieties
parameterise the abelian surfaces which split. These loci are described to
some extent in \cite{BN}.
\begin{theorem}\label{t131.1}\begin{itemize}\item[a)] Points on $\cN_5$
    parameterise smooth quartic surfaces unless they lie on one of the 30
    planes of Lemma \ref{l125.1};
  \item[b)] points on $\cN_5$ parameterise quartic surfaces containing more
    that 32 lines if an only if the corresponding abelian surfaces are
    products. Furthermore, a line on a surface of this set of quartic
    surfaces has coordinates in $\fP^5$ which is in the $\gS_6$-orbit of

    $$x_0^4(x_1^2+x_2^2)+x_1^4(x_2^2+x_0^2)+x_2^4(x_0^2+x_1^2)-6x_0^2x_1^2x_2^2
    =0.$$
\end{itemize}
\end{theorem}
Unfortunately, these result do not allow us to explicitly describe the
relation between the compactification $\cM_{20}$ and compactifications of
$\cA_{(1,3)}(2)$, in particular the Baily-Borel embedding. This must be
considered an interesting open problem.

\subsection{A conjecture}
To end this section we make a conjecture on one of the birational models of
the variety $\cN_5$. Consider the birational map $\cN_5- - \ra \hat{\cN}_5$
of (\ref{e127.2}). Recalling now the Janus-like isomorphism between the
Picard modular variety $\-{X}_{\gG_{\sqrt{-3}}(\sqrt{-3})}$ and the Siegel
modular variety $\-{X}_{\gG(2)}$ (see \cite{J}), it is natural to ask
about an analogue
here, since the involved Siegel modular varieties of Corollary \ref{c130.1}
all are related to level
2, albeit with different polarisations. So consider abelian fourfolds with
complex multiplication by $\fQ(\sqrt{-3})$ of signature (3,1), with a level
$\sqrt{-3}$ structure, but with (1,1,1,3) polarisations.
\begin{problem}\label{p132.1} Is $\hat{\cN}_5$ the Satake compactification of
  $X_{(1,1,1,3)}(\sqrt{-3}):=\gG_{(1,1,1,3)}(\sqrt{-3})\bs \fB_3$, where
  $\gG_{(1,1,1,3)}(\sqrt{-3})$ denotes the arithmetic group giving
  equivalence of complex multiplication by $\fQ(\sqrt{-3})$, signature
  $(3,1)$, with a level
  $\sqrt{-3}$ structure and a $(1,1,1,3)$-polarisation?
\end{problem}
I conjecture that for {\em some} subgroup of $\gG_{(1,1,1,3)}(\sqrt{-3})$,
this does in fact hold.

Evidence for the conjecture:
\begin{itemize}\item[i)] The proper transforms of the 15 Segre planes are
  by Proposition \ref{p116.1} the moduli space of principally polarised
  abelian threefolds with complex multiplication by $\fQ(\sqrt{-3})$,
  signature (2,1), with a level
  $\sqrt{-3}$ structure, (although these moduli spaces are blown up in three
  points on $\hat{\cN}_5$). These could parameterize
  abelian fourfolds with said CM, signature (3,1) with a level
  $\sqrt{-3}$ structure and polarisation $(1,1,1,3)$ which split:
  $$A_4\isom A_3\times A_1,$$ where $A_3$ has CM, signature (2,1),
  polarisation (1,1,1), and $A_1$ has CM, but a polarisation 3.
\item[ii)] The proper transforms of the 15 planes $N_{ij}$ of Lemma
  \ref{l125.1}, (ii), are four nodal cubic surfaces (Lemma \ref{l127.3}).
  These surfaces occur also on the ball quotient $\cS_3$ above: pick any
  four of the nodes which are not coplanar; they determine a unique $\fP^3$
  in $\fP^4$, and its intersection with $\cS_3$, a cubic surface, has four
  nodes in the four nodes of $\cS_3$ in that $\fP^3$. (Note that there is a
  unique four-nodal cubic surface, as it is $\fP^2$ blown up in the six
  intersection points of four (general) lines, a complete quadrilateral in
  $\fP^2$, and any two such quadrilaterals are projectively equivalent.
  This cubic surface is usually called the Cayley cubic, mentioned above.)
\item[iii)] The singular locus consists of isolated singular points,
  resolved by quadric surfaces, so these singularities are rational.  Recall
  that at each $P_{ij}$, six of the 15 Segre planes meet. At each $Q_{ij}$
  (the 15 points (\ref{e124.2})), three of the Segre planes and six of the
  $N_{ij}$ of Lemma \ref{l125.1}, (ii) meet. In both cases, the quadric
  $\fP^1\times \fP^1$ can be covered equivariantly by a product
  $E_{\grr}\times E_{\grr}$
  of the elliptic curve $E_{\grr}$ with branching only at the intersection
  with the proper transforms of the 30 planes above, as follows:
\begin{itemize}\item $P_{ij}$: $E_{\grr}\times E_{\grr}\lra \fP^1\times
  \fP^1$ a Galois $\fZ/3\fZ$-quotient;
\item $Q_{ij}$: $E_{\grr}\times E_{\grr}\lra \fP^1\times\fP^1$ is the
  product of two double covers branched at $0,1,\grr, \grr^2$.
\end{itemize}
This supports by Lemma \ref{l115.1} the idea that this could be the
compactification locus of $X_{(1,1,1,3)}(\sqrt{-3})$.
\end{itemize}

\section{The Coble variety $\cal Y$}
\subsection{Arithmetic quotients of domains of type $\bf IV_n$}

Let $V$ be a $k$ vector space of dimension $n$, $k$ a totally real field,
and $b$ a bilinear symmetric form on $V$. Let $G(V,b)$ be the symmetry
group, and $G_{\fQ}=Res_{k|\fQ}G(V,b)$ the $\fQ$-group it defines. We
assume that $G_{\fQ}$ is of hermitian type, so that for every infinite prime
$\nu$ of $k$ the signature of
$b_{\nu}$ is $(n-2,2)$ or $b_{\nu}$ is definite.
$G_{\fQ}$ is (absolutely) simple (defines an irreducible
domain) only if $k=\fQ$ (or if $b_{\nu}$ is definite for all but a single
$\nu$, in which case the $\fQ$-group is anisotropic, but we will not
consider this situation), and the corresponding real group gives rise to a
bounded symmetric domain only if $b$ has Witt index 2. This is the case we
consider here. The classification of such forms is well-known; since we
require the Witt index to be 2, two such forms $b$ and $b'$ are equivalent
over $\fQ$ $\iff$ $det(b)=det(b')$, where $det(b)$ is to be viewed as an
element of $\fQ^{\times}/(\fQ^{\times})^2$.

Now let $\cL\inn V$ be a (maximal) lattice, and let $G_{\cL}$ be the
arithmetic group it defines, $\gG\inn G_{\cL}$ a subgroup of finite index.
We first remark on the moduli interpretation of the arithmetic quotient
$\xg$.
\begin{proposition}\label{p106.1} $\xg$ is a moduli space of (pure) Hodge
  structures of weight 2 on $V$ with $h^{2,0}=1$ (and $h^{1,1}=dim(\xg)$)
  with respect to the lattice $\cL\inn V$.
\end{proposition}
{\bf Proof:} The symmetry group of such a Hodge structure is of real type
$SO(n-2,2)$, and $G(V,b)$ is a $\fQ$-form in which $G_{\cL}$ is an
arithmetic subgroup. Since the corresponding ``period'' (i.e., position of
the varying complex subspace $H^{1,1}$ in $H^2_{\fC}$) is clearly the same
exactly when the two periods differ by an element of $G_{\cL}$, while $\gG$
defines a level structure of some kind, the result follows. \ende This
proposition is often used in the study of polarised K3-surfaces, which have
a pure Hodge structure of type (1,19,1). In fact, for each polarisation
degree (i.e., the number $C^2$ for the ample divisor $C$ on the K3-surface
which gives the projective embedding) $2e,\ (e\geq1)$ one has an arithmetic
group $\gG_e$ such that the arithmetic quotient $X_{\gG_e}$ is the moduli
space of K3 surfaces with the given polarisation. Recall the {\em Picard
  number} $\grr$ is the rank of the group of algebraic cycles, i.e.,
$\grr=rk_{\fZ}H^2(S,\fZ)\cap H^{1,1}$. Then one has the following.
\begin{proposition}\label{p106.2} Let $S$ be a K3 surface with $\grr$= the
  Picard number of $S$. Then the dimension of the moduli space of K3's
  which are in the family preserving the lattice of algebraic cycles
  $H^2(S,\fZ)\cap H^{1,1}$ is 20-$\grr$.
\end{proposition}
{\bf Proof:} Recall that for a K3 surface $H^1(S,\gT)\isom H^1(S,\gO^1)$,
so $H^{1,1}(S)$ may be viewed as the tangent space of the local deformation
space, which should be thought of as a varying complex subspace of
$H^2(S,\fC)_{prim}$, while $H^2(S,\fZ)$ is fixed. Let $\scA=H^2(S,\fZ)\cap
H^{1,1}$ be the lattice of algebraic cycles, $\scT=H^2(S,\fZ)\cap
(H^{2,0}(S)\oplus H^{0,2}(S))$ the lattice of transcendental cycles. We
have $rk_{\fZ}\scA=rk_{\fZ}(H^2(S,\fZ))-rk_{\fZ}\scT$ in general and
$rk_{\fZ}\scA=\grr$ by assumption, so $\grr=22-rk_{\fZ}\scT$, while the
moduli space is defined by the group $G(V',b')$, where
$V'=\scA^{\perp}\otimes \fQ$, since we are requiring $\scA$ to be
preserved. (Recall that for an algebraic cycle $\cC$ the integral
$\int_{\cC}\go$ over the holomorphic two-form $\go$ vanishes, hence the
algebraic cycles contribute nothing to the periods).  Thus
$G(\fR)=SO(20-\grr,2)$, giving rise to a domain of type $\bf IV_{20-\grr}$.
\ende Of course in this particular case, the lattice $\cL\inn V$ is very
special; the ``intersection form'' $b$ restricted to $\cL$ is even and
unimodular, and as is well-known, decomposes as
\begin{equation}\label{e106.1} \cL\isom <-2e>\oplus {\bf H}^2 \oplus {\bf
    E_8}^2,
\end{equation}
where ${\bf H}$ is the two-dimensional hyperbolic lattice, and ${\bf E_8}$
is the root lattice of type ${\bf E_8}$. Let us remark that the
compactification of these arithmetic quotients has been carried out in the
thesis \cite{scat}, but we will not need this.  We will very quickly
describe a particulary interesting family of K3 surfaces which has been
thoroughly studied by Yoshida and his collaborators, see \cite{MSY} for
details on all matters here.

\subsection{A four-dimensional family of K3's}\label{s106a.1}
The family of K3 surfaces to be described here is the set of surfaces which
are double covers of $\fP^2$ branched along the union of six disjoint
lines.  Recall that there is a 19-dimensional family of K3 surfaces which
are double covers of the plane branched along a sextic curve; they are
smooth as long as the sextic is smooth, and generically have Picard
number $\grr=1$. An arrangement of six lines in $\fP^2$ is a maximally singular
sextic; there are 15 intersection points of the six lines (if they are in
general position), and each such gives rise to an $A_1$-singularity on the
double cover. Resolving the 15 double points introduces 15 exceptional
curves with self intersection number $-2$, so together with the pullback of
the generic line, this gives 16 independent cycles on the surface:
$\grr=16$. Hence the transcendental lattice $\scT$ has rank 4, so by
Proposition \ref{p106.2}, the moduli space is four-dimensional.  Let
\begin{equation}\label{e106a.1} \gG=\{g\in G(\scT_{\fR},Q)(\isom
  SO(4,2))\Big| g(\scT)\inn \scT\},
\end{equation}
where $Q$ is the intersection form on $H^2(S,\fZ)$, extended to $\fR$, then
restricted to $\scT_{\fR}$. This is clearly an arithmetic subgroup, and by
Proposition \ref{p106.1}, the arithmetic quotient $\xgeq$ is the
four-dimensional moduli space. We list some of the interesting loci for
this family. Let $L$ be the given arrangement, $L=l_1\cup \ldots \cup l_6$,
and let $t_p:=$ the number of $p$-fold points of the arrangement, i.e., the
number of points at which $p$ of the lines meet (see (\ref{e109.1})), and
let $\pi:S\lra \fP^2$ denote the (singular) double cover.

\subsubsection{Three-dimensional loci}
\begin{itemize}
\item[1)] Suppose there is a conic which is {\em tangent} to all six lines.
  Then the inverse image of this quadric is a $\fP^1$, which, as is easily
  checked, has self-intersection number $4-6=-2$, so the double cover has
  16 exceptional cycles, hence $\grr=17$. It is in fact easy to see that
  the surface $S$ is in this case a classical Kummer surface, i.e., a
  quartic surface in $\fP^3$ with 16 nodes which is the Kummer variety of a
  principally polarised abelian surface $A_S$. The projection from a node
  gives the double cover $\pi:S\lra \fP^2$, and the tangent conic is the
  image of the (blown up) node used to project. The abelian surface is the
  Jacobian of a genus 2 curve, and this curve is the double cover of the
  conic, {\em branched at the six points of tangency}. This is well-known.

\item[2)] If $t_3=1,\ t_2=12$, then the threefold point induces an
  $A_2$-singularity on the double cover which is resolved by two $\fP^1$'s,
  so there are now 2+12 exceptional $\fP^1$'s and the hyperplane section.
  We have the following picture. \\ \fbox{
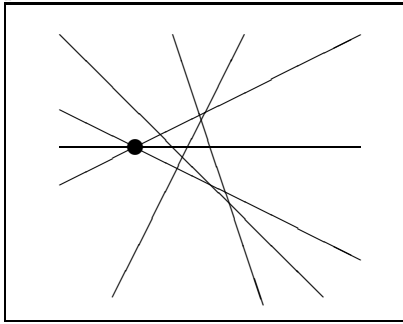
\begin{minipage}{5cm}
  \unitlength1cm
\begin{picture}(5,4)(0.5,0.3)
  \put(1,3){\line(2,-1){4}} \put(2,2.5){\circle*{.2}}
  \put(1,2.5){\line(1,0){4}} \put(1,2){\line(2,1){4}}
  \put(1,4){\line(1,-1){3.5}} \put(2.5,4){\line(1,-3){1.2}}
  \put(1.7,.5){\line(1,2){1.75}}
\end{picture}
\end{minipage}}
\begin{minipage}{10cm} There are in fact three more exceptional $\fP^1$'s,
  which are the inverse image on the double cover of the three lines which
  pass through the triple point and one of the three double points not
  lying on a line through the triple point. It is easy to see that these
  three double points are independent parameters of such arrangements, so
  this defines a three-dimensional family, so by Proposition \ref{p106.2},
  we have $\grr=17$ for the generic member of this family.
\end{minipage}

\end{itemize}

\subsubsection{Two-dimensional loci}
\begin{itemize}
\item[3)]

  \fbox{
\begin{minipage}{5cm}
  \unitlength1cm
\begin{picture}(5,4)(0.5,0.3)
  \put(1,2.5){\line(2,-1){4}} \put(2,2){\circle*{.2}}
  \put(1,2){\line(1,0){4}} \put(1,1.5){\line(2,1){4}}
  \put(2.5,.5){\line(1,2){1.75}} \put(3.35,2.3){\circle*{.2}}
  \put(1,3.5){\line(2,-1){4}} \put(2.55,4){\line(1,-2){1.75}}
\end{picture}
\end{minipage}}
\begin{minipage}{10cm}If $t_3=2, t_2=9$, there are two possibilities.
  Suppose first that the two threefold points do {\em not} lie on one of
  the six lines.  Then we have the picture to the left. This gives rise to
  two isolated $A_2-$singularities.  The inverse image of the line joining
  the two threefold points is also an exceptional $\fP^1$. In this case the
  generic double cover has $\grr=18$, and as parameters one can take two
  double ratios: consider two of the lines $l_1, l_2$, both passing through one
  of the threefold points $p$; the three intersection points with the other
  lines, together with $p$, give four points on each line -- hence two
  double ratios.
\end{minipage}

\item[4)] It may also occur that both threefold points lie on a line, but
  in this case we also have $\grr=18$, i.e., a two-dimensional family.
\end{itemize}

\subsubsection{One-dimensional loci}
\begin{itemize}
\item[5)] If $t_4=1$, then the double cover has an {\em elliptic}
  singularity over the point, so is not K3. Hence this is a genuine {\em
    degeneration} of the K3, i.e., belongs to the boundary of the
  compactification. It turns out that then a line must also be double, so
  that the double cover has two components.
\item[6)] As a further specialisation of 4) it may happen that there are
  three triple points. Since four of the lines may be choosen fixed (for
  example $x_0=x_1=x_2=0,\ x_0-x_1=0$), there is only one modulus, given
  for example by the intersection point of the two variable lines. Here we
  have $\grr=19$.

\end{itemize}

\subsubsection{Zero-dimensional loci}
\begin{itemize}
\item[7)] If three lines are taken, each {\em double}, then the double
  cover splits into two copies of $\fP^2$. This is in the closure of the set
  of degenerations of type 5).
\item[8)] The arrangement is the complete quadrilateral. The picture is: \\
  \fbox{
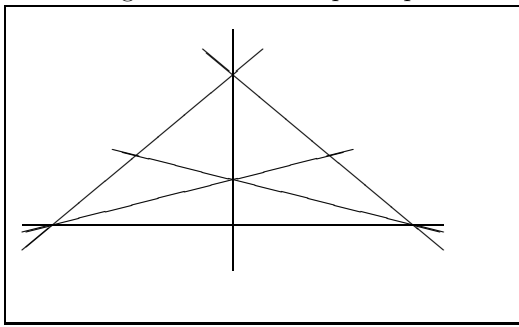
\begin{minipage}{6.5cm}
  \unitlength.8cm
\begin{picture}(8,5)(1.5,0)
  \thinlines \put(2,1.5){\line(1,0){6.5}} \put(2,1.5){\line(4,1){5}}
  \put(2,1.5){\line(6,5){3.5}}

  \put(2,1.5){\line(-1,0){0.5}} \put(2,1.5){\line(-4,-1){0.5}}
  \put(2,1.5){\line(-6,-5){0.5}}

  \put(5,.75){\line(0,1){4}} \put(5,4){\line(6,-5){3.5}}
  \put(5,4){\line(-6,5){0.5}} \put(8,1.5){\line(-4,1){5}}
  \put(8,1.5){\line(4,-1){0.5}}

\end{picture}
\end{minipage}}
\begin{minipage}{9cm}It is known that the {\em Fermat} cover (not the double
  cover) of this arrangement is Shioda's elliptic modular surface of level
  4, $S(4)$, so it follows that the double cover is isogenous to $S(4)$,
  i.e., a quotient of $S(4)$ by a group isomorphic to $(\fZ/2\fZ)^4$. This
  is the most special K3 surface in the family and has $\grr=20$.
\end{minipage}
\end{itemize}

\subsubsection{Level 2 structure}
Now consider, in addition to the above data, a level 2 structure.
Geometrically this amounts to fixing an {\em order} of the six lines. In
terms of the lattice $\scT$ it is not so easy to see what it means. In
\cite{MSY} it is shown by explicit computation that the subgroup $\gG(2)$
is the group generated by reflections on the ``roots'' of $\scT$, that is
the integral elements of norm $-2$. Furthermore it is shown there that
$\gG$ is generated by the reflections on the elements of norm $-2$ or $-4$,
and that $\gG/\gG(2)\isom \gS_6\times \fZ/(2)$. Hence by the results 2.7.1,
2.7.7, 2.8.2 of \cite{MSY} we have
\begin{proposition}\label{p106b.1} The arithmetic quotient $\gG(2)\bs \cD$ is
  the moduli space for K3 surfaces which are double covers of $\fP^2$,
  branched over an {\em ordered} set of six lines.
\end{proposition}
\begin{table}
\caption{\label{table18} Loci of a four-dimensional family of K3 surfaces }

\vspace*{.5cm}
\begin{minipage}{16.5cm}
  \hspace*{2.5cm}\fbox{\begin{minipage}{5cm} Locus 1)

      Igusa quartic \unitlength1cm
\begin{picture}(2,2)(0,.2)
  \put(0,1){\circle{.2}} \put(1,0.5){\circle{.2}} \put(2,1){\circle{.2}}
  \put(0,2){\circle{.2}} \put(1,2.5){\circle{.2}} \put(2,2){\circle{.2}}
  \bezier{100}(0,1)(-.2,1.5)(0,2) \bezier{100}(0,2)(0.5,2.5)(1,2.5)
  \bezier{100}(1,2.5)(1.5,2.5)(2,2) \bezier{100}(2,2)(2.3,1.5)(2,1)
  \bezier{100}(2,1)(1.5,0.45)(1,.5) \bezier{100}(1,.5)(0.5,0.45)(0,1)
\end{picture}

{\hspace*{\fill} \fbox{1} \hfill}
\end{minipage}}
\hspace*{\fill} \fbox{\begin{minipage}{5cm} Locus 2)

    $X^{\{ijk\}}$ \unitlength1cm
\begin{picture}(2,2)(0,.2)
  \put(1,1){\circle{.2}} \put(2,1.5){\circle{.2}} \put(1,2){\circle{.2}}
  \put(3,1.8){\circle{.2}} \put(1,1){\circle{.2}}
  \put(1.1,1){\line(1,0){.9}} \put(2,1){\circle{.2}}
  \put(2.1,1){\line(1,0){.9}} \put(3,1){\circle{.2}}
\end{picture}

{\hspace*{\fill} \fbox{20} \hfill}
\end{minipage}} \hspace{\fill}

\begin{minipage}{15cm}\unitlength1cm \begin{picture}(15,2)
    \put(5,.7){\vector(-1,1){.8}} \put(8,.7){\vector(1,1){.8}}
    \put(13,.7){\vector(-1,1){.8}}
\end{picture}
\end{minipage}

\hspace*{3cm} \fbox{\begin{minipage}{5cm} Locus 3)

    $X^{\{ijk;lmn\}}$ \unitlength1cm
\begin{picture}(2,2)(.2,.3)
  \put(0,1){\put(1,1){\circle{.2}} \put(1.1,1){\line(1,0){.9}}
    \put(2,1){\circle{.2}} \put(2.1,1){\line(1,0){.9}}
    \put(3,1){\circle{.2}} } \put(1,1){\circle{.2}}
  \put(1.1,1){\line(1,0){.9}} \put(2,1){\circle{.2}}
  \put(2.1,1){\line(1,0){.9}} \put(3,1){\circle{.2}}
\end{picture}

{\hspace*{\fill} \fbox{10} \hfill}
\end{minipage}}
\hspace{\fill} \fbox{\begin{minipage}{5cm} Locus 4)

    $X^{\{ijk;imn\}}$ \unitlength1cm
\begin{picture}(2,2)
  \put(.5,1){\circle{.2}} \put(.6,1.1){\line(2,1){.9}}
  \put(.6,1.1){\line(2,-1){.9}} \put(1.5,.5){\circle{.2}}
  \put(1.5,1.5){\circle{.2}} \put(1.6,.6){\line(2,-1){.9}}
  \put(1.6,1.6){\line(2,1){.9}} \put(2.5,1.4){\circle{.2}}
  \put(2.5,0){\circle{.2}} \put(2.5,2){\circle{.2}}
\end{picture}

{\hspace*{\fill} \fbox{90} \hfill}
\end{minipage}}

\begin{minipage}{15cm} \unitlength1cm \begin{picture}(15,2)
    \put(5,.7){\vector(0,1){.8}} \put(9,.7){\vector(1,1){.8}}
    \put(13,.7){\vector(0,1){.8}}
\end{picture}
\end{minipage}

\hspace*{3cm} \fbox{\begin{minipage}{5cm} Locus 5)

    $X^{\{ij\}}$ \unitlength1cm
\begin{picture}(2,2)
  \put(.5,1){\circle{.2}}\put(.6,1){\line(1,0){.9}}
  \put(1.5,1){\circle{.2}}\put(1.6,1){\line(1,0){.9}}
  \put(2.5,1){\circle{.2}}\put(2.6,1){\line(1,0){.9}}
  \put(2,2){\circle{.3}}\put(2,2){\circle*{.2}}
\end{picture}

{\hspace*{\fill} \fbox{15} \hfill}
\end{minipage}}
\hspace{\fill} \fbox{\begin{minipage}{5cm} Locus 6)

    $X^{\{ijk;klm;mni\}}$ \unitlength1cm
\begin{picture}(2,2)(0,.2)
  \put(.5,.5){\circle{.2}}\put(.6,.5){\line(1,0){.9}}
  \put(.5,.5){\line(1,2){.4}}
  \put(1.5,.5){\circle{.2}}\put(1.6,.5){\line(1,0){.9}}
  \put(1,1.5){\circle{.2}}\put(1,1.5){\line(1,2){.4}}
  \put(1.5,2.5){\circle{.2}}\put(1.5,2.5){\line(1,-2){.4}}
  \put(2,1.5){\circle{.2}} \put(2,1.5){\line(1,-2){.4}}
  \put(2.5,.5){\circle{.2}}
\end{picture}

{\hspace*{\fill} \fbox{120} \hfill}
\end{minipage}}

\begin{minipage}{15cm} \unitlength1cm \begin{picture}(15,2)
    \put(5,.7){\vector(0,1){.8}} \put(9,.7){\vector(1,1){.8}}
    \put(13,.7){\vector(0,1){.8}}
\end{picture}
\end{minipage}

\hspace*{3cm} \fbox{\begin{minipage}{5cm} Locus 7)

    $X^{\{ijk;kl;mn\}}$ \unitlength1cm
\begin{picture}(2,2)
  \put(1,1){\circle{.3}}\put(1,1){\circle*{.2}}
  \put(2,1){\circle{.3}}\put(2,1){\circle*{.2}}
  \put(1.5,2){\circle{.3}}\put(1.5,2){\circle*{.2}}
\end{picture}

{\hspace*{\fill} \fbox{15} \hfill}
\end{minipage}}
\hspace{\fill} \fbox{\begin{minipage}{5cm} Locus 8)

    $X^{\{ijk;klm;mni;jln\}}$ \unitlength1cm
\begin{picture}(2,2)(.9,.3)

  \put(.5,1.5){\circle{.2}}\put(1.6,.9){\circle{.2}}\put(1.6,2.1){\circle{.2}}
  \put(2,1.5){\circle{.2}}\put(2.7,.3){\circle{.2}}\put(2.7,2.7){\circle{.2}}
  \put(0,1.25){\line(2,1){3}}\put(0,1.75){\line(2,-1){3}}
  \put(1.33,.5){\line(2,3){1.5}}\put(1.33,2.5){\line(2,-3){1.5}}
\end{picture}

{\hspace*{\fill} \fbox{30} \hfill}
\end{minipage}}

\end{minipage}
\vspace*{.5cm} {\small The notations $X^{\{ijk;klm;mni;jln\}}$, etc, are
  taken from \cite{MSY}; the arrows indicate inclusions among the various
  loci. \hfill\newline

\vspace*{-1.2cm}The symbol \unitlength1cm
\begin{picture}(.15,.1)\put(0.2,0.1){\circle{.3}}\put(0.2,0.1){\circle*{.2}}
\end{picture} \quad means a double line. The number of each kind of loci is
indicated by \fbox{x}; the dimensions are three in the top row down to zero
in the last row. Locus 1) is where the six points lie on a conic, while the
20 $X^{\{ijk\}}$ are the loci where there are three of the six points on a
line. The 15 $X^{\{ijk;kl;mn\}}$ lie on the boundary of the moduli space,
while the 30 $X^{\{ijk;klm;mni;jln\}}$ lie in ``the farthest interior'' of
the domain. A more complete description is given in Corollary
\ref{c133.1}.}

\end{table}

We refer the reader to \cite{MSY} for a detailed description of the loci
described above, of the periods and of the corresponding Picard-Fuchs
equations (and much more). We give in Table \ref{table18} a description of
the loci, giving the dual graph of the six lines (i.e., a vertex for each
line, two vertices lying on a line $\iff$ the corresponding lines meet), as
well as the number of loci, and the names given to them in \cite{MSY}.

We now give an explicit projective description of the
Baily-Borel compactification of the arithemetic quotient $\gG(2)\bs \cD$ of
Proposition \ref{p106b.1}. All the facts presented here were proved
originally by Coble \cite{C} or by Yoshida and his collaborators in
\cite{MSY}. We have the four-dimensional family of K3 surfaces just
discussed, defined in terms of a set of (ordered) six lines in
the plane. Dual to the six lines are six points, and so the relation with
the moduli space of cubic surfaces is evident.  Let two ordered sets of six
lines, $(l_1,\ldots,l_6),\ (l_1',\ldots,l_6')$ be given.
\begin{definition}\label{d133.0} The two sets of lines $(l_1,\ldots,l_6),\
  (l_1',\ldots,l_6')$ are said to be {\em associated}, if the following
  relation holds. Since the set $(l_1,\ldots,l_6)$ is ordered, we can form
  two triangles, $$\gD(l_1,l_2,l_3),\quad \gD(l_4,l_5,l_6);$$ these two
  triangles have together six vertices, which come equipped with a
  numbering, say $(p_1,\ldots,p_6)$, and these correspond dually to another
  ordered set of six lines, $(l_{p_1},\ldots,l_{p_6})$. Then
  $(l_1,\ldots,l_6)$ and $(l_1',\ldots,l_6')$ are associated, if:
  $(l_{p_1},\ldots,l_{p_6})=(l_1',\ldots,l_6')$, as a set of six ordered
  lines.
\end{definition}
Of course, starting with two sets of ordered six points, one can define in
the same way the notion of association. Since, as abstract moduli spaces,
the space of ordered sets of six lines is the ``same'' (by duality) as the
set of ordered sets of six points, we see that we are dealing here with the
space of sets of six ordered points in $\fP^2$. This problem was dealt with
in the papers of Coble \cite{C}, and has been given a modern treatment in
\cite{DO}.  It can be described as follows.  The relevant moduli space is
easy to describe: let $(p_1,...,p_n)$ be a set of $n$ points in $\fP^k$;
this is represented by $M$, the $n\times (k+1)$ matrix whose $i^{th}$
column gives the coordinates of the point $p_i$. The moduli space is then
the GIT quotient
\begin{equation}\label{e133.2}
  \hbox{\pkn} = GL(k+1) \bs M(n,k+1)/(\fC^*)^n.
\end{equation}
By taking the set of semistable points in $M(n,k+1)$ the above quotient is
compact, although singular. It is classical that \pos\ is a threefold whose
compactification can be identified with a cubic threefold in $\fP^4$ with
ten ordinary double points, which is just the Segre cubic $\cS_3$. Note
that the similar moduli problem, namely six points on a conic in $\fP^2$,
is realised by the Igusa quartic $\cI_4$, so these are very closely
related, but not identical moduli problems.

Our interest here is in \pts, a fourfold.  In this case we may represent
elements by matrices
\begin{equation}\label{e133.3}
  \hbox{\pkn} \ni M = \left[ \begin{array}{c c c c c c} 1 & 0 & 0 & 1 & x &
      w \\ 0 & 1 & 0 & 1 & y & z \\ 0 & 0 & 1 & 1 & u & u \\
\end{array} \right] ,
\end{equation}
and as Coble shows, the map \pts $\lra \fP^4$, $M \mapsto [x:y:w:z:u]$ is a
birational map (it is clear that $\pts$ is rational, this map simply gives
an explicit birationalisation). The GIT theory here consists of finding
$G$-invariant functions on \pts, and these turn out to be generated by
$3\times 3$ minors of $M$.

In terms of the matrix $M$ the process of association can be described as
follows. Each such matrix $M$ determines a second one: since the six points
are ordered, one can define six lines by $l_{12}=\ove{p_1p_2}$,
$l_{13}=\ove{p_1p_3}$, $l_{23}=\ove{p_2p_3}$, $l_{45}=\ove{p_4p_5}$,
$l_{46}=\ove{p_4p_6}$, $l_{56}=\ove{p_5p_6}$; these six lines determine
dually six points, whose coordinates are then brought into the normal form
given above. It turns out that the entries of the second matrix are
determined by the fact that the maximal minors are proportional to the
maximal minors of the first. More precisely, if we let $(ijk)$ denote the
$3\times 3$ minor of $M$ which is given by the columns $i,j,k$, and if we
let $M'$ be the associated matrix, $(ijk)'$ the corresponding minor, then
the minors of $M$ and $M'$ are related by:
\begin{equation}\label{e133.1} (123)(145)(246)(356)=(124)'(135)'(236)'(456)'.
\end{equation}
Now association is an involution on $\pts$, and one can take the {\em
  quotient} by this involution.
\begin{definition}\label{d133.1} Let $\cY$ be the double cover of $\fP^4$
  branched along the Igusa quartic $\cI_4$, $\pi:\cY\lra \fP^4$.
\end{definition}
Clearly $\cY$ will be singular precisely along the singular locus of
$\cI_4$, i.e.,
\begin{lemma}\label{l133.1} The singular locus of $\cY$ consists of 15 lines,
  the inverse images of the 15 singular lines of $\cI_4$.
\end{lemma}
\begin{theorem}[\cite{DO},Example 4, p.~37]\label{t133.0}
  The moduli space of six ordered points in $\fP^2$ is equal to the double
  cover $\cY$, and the double cover involution on $\cY$ coincides with the
  association involution on $\pts$.
\end{theorem}
In other words, a set $(p_1,\ldots,p_6)$ is {\em associated to itself}, if
and only if the six points lie on a conic in $\fP^2$.

Consider one of the hyperplanes $H$ in $\fP^4$, $H=\{\gt_m^4=0\}$ of
Proposition \ref{p120.1}. Since $H$ is {\em tangent} to $\cI_4$, the
inverse image $\pi^{-1}(H)$ in $\cY$ will {\em split into two copies of
  $\fP^3$}. In this way, we get a union of 20 $\fP^3$'s on $\cY$,
\begin{lemma}\label{l133.2} The inverse images $\pi^{-1}(H)$ of the tangent
  hyperplanes $H=\{\gt_m^4=0\}$ consist of two copies each of $\fP^3$, and
  these two $\fP^3$'s on $\cY$ meet in the quadric surface which is the
  inverse image under $\pi$ of the quadric on $\cI_4$ to which $H$ is
  tangent. This gives a total of 20 such $\fP^3$'s on $\cY$.
\end{lemma}

A resolution of singularities of $\pts$ is affected by resolving the Igusa
quartic by blowing up the ideal of the 15 lines; this is the map $\grr_2$
of Theorem \ref{t123.1}.  Let $\hat{\hbox{\pts}}$ denote this
desingularisation $\hat{\hbox{\pts}} \lra \hbox{\pts}$. On
$\hat{\hbox{\pts}}$ we have a set of 36 divisors, the {\em discriminant
  locus}, the proper transforms of the Igusa quartic, the 20 $\fP^3$'s and
the 15 exceptional divisors of the blow up.

It is clear how this variety is the moduli space of cubic surfaces: blow up
$\fP^2$ in the six points, and embed by the linear system of cubic curves
through the six points. The ordering of the six points of course determines
a marking of the 27 lines in the well-known manner. The symmetry group of
$\hat{\hbox{\pts}}$ is $\gS_6\times \fZ/2\fZ$; although the Weyl group
$W(E_6)$ acts birationally on it, the action is not regular. For that it is
neccessary to modify $\hat{\hbox{\pts}}$ even more. Dolgachev mentions in
\cite{DO} that he suspects it is sufficient to blow up $\hat{\hbox{\pts}}$
in the intersection of the 36 divisors.

One of the many things proved in \cite{MSY} is the following.
\begin{theorem}\label{t133.1} The variety $\cY$ is the Baily-Borel
  compactification of the arithmetic quotient $\gG(2)\bs \cD$ of
  Proposition \ref{p106.1}.
\end{theorem}
The proof given in \cite{MSY} of this fact simply (!)  calculates the image
of the period map, and in determining when the periods lie on the boundary
of the period domain $\cD$, the authors find that this locus
coincides with the set of
K3 surfaces whose set of six lines correspond to those singularities of
Lemma \ref{l133.1} of $\cY$.
\begin{corollary}\label{c133.1} The Loci 5) and 7) of Table \ref{table18}
  are the inverse images on $\cY$ of the 15 singular lines and 15 singular
  points, respectively, of the branch locus $\cI_4$. The Loci 3) of Table
  18 are the inverse images of the ten special hyperplane sections of Lemma
  \ref{l133.2}, i.e., the quadrics.  The loci 2) of Table 18 are the 20
  $\fP^3$'s of Lemma \ref{l133.2}, and Locus 1) is just the branch locus of
  the double cover.
\end{corollary}

\part{A Gem of the modular universe}\label{chapter13}
\renewcommand{\arraystretch}{1}

\section{The Weyl group $W(E_6)$}

\subsection{Notations}
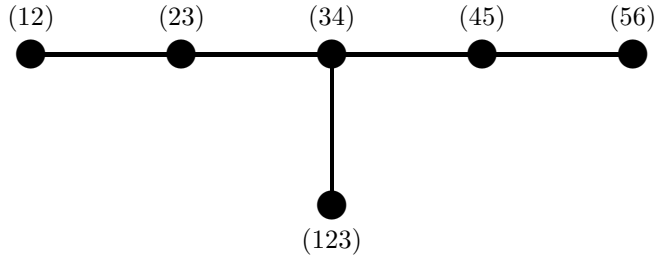
\begin{figure}[hbt]
  $$ \unitlength0.8cm
\begin{picture}(10.5,3.5)

  \put(0.25,3.0){\circle*{0.5}} \put(2.75,3.0){\circle*{0.5}}
  \put(5.25,3.0){\circle*{0.5}} \put(7.75,3.0){\circle*{0.5}}
  \put(10.25,3.0){\circle*{0.5}} \put(5.25,0.5){\circle*{0.5}}

  \linethickness{0.4mm}

  \put(0.25,3.0){\line(1,0){2.5}} \put(2.75,3.0){\line(1,0){2.5}}
  \put(5.25,3.0){\line(1,0){2.5}} \put(7.75,3.0){\line(1,0){2.5}}
  \put(5.25,0.5){\line(0,1){2.5}}

  \put(0.25,3.6){\makebox(0,0){(12)}} \put(2.75,3.6){\makebox(0,0){(23)}}
  \put(5.25,3.6){\makebox(0,0){(34)}} \put(7.75,3.6){\makebox(0,0){(45)}}
  \put(10.25,3.6){\makebox(0,0){(56)}}
  \put(5.25,-0.1){\makebox(0,0){(123)}}

\end{picture}
$$
\caption[The graph of the group of the 27 lines]{\label{Figure4}\small
  The graph of the group of the permutations of the 27 lines}
\end{figure}

We use the same notation as above for the 27 lines on a cubic surface in
$\fP^3$: $a_1,...a_6,\ b_1,...,b_6,\ c_{12},...,c_{56}$.  The 36 double
sixes are:
$$N= \left[ \matrix{ a_1 & a_2 & a_3 & a_4 & a_5 & a_6 \cr b_1 & b_2 & b_3
    & b_4 & b_5 & b_6 \cr } \right], \hspace{1cm} (1) $$
$$N_{ij}=\left[ \matrix{ a_i & b_i & c_{jk} & c_{jl} & c_{jm} & c_{jn} \cr
    a_j & b_j & c_{ik} & c_{il} & c_{im} & c_{in} \cr } \right],
\hspace{1cm} (15)$$
$$N_{lmn}\footnote{here we switch notations from $N_{ijk}$ in equation
  (\ref{eB2.2}) to $N_{lmn}$ for convenience} =\left[ \matrix{ a_i & a_j &
    a_k & c_{mn} & c_{ln} & c_{lm} \cr c_{jk} & c_{ik} & c_{ij} & b_l & b_m
    & b_n \cr} \right] \hspace{1cm} (20).  $$ The 45 tritangents
are:\begin{equation}\label{eQ1.1}\begin{minipage}{6cm}
\begin{center}
  $(ij)=<a_i\ b_j\ c_{ij}>,\ \ i\neq j \ \ (30)$

  $(ij.kl.mn)=<c_{ij}\ c_{kl}\ c_{mn}>\ \ (15).$
\end{center}
\end{minipage}
\end{equation}
Two double sixes are {\em syzygetic} it they contain four lines in common,
for example: $$N\ \ \hbox{and} \ \ N_{12}\ \ \hbox{have } a_1,\ a_2,\ b_1,\
b_2\ \hbox{in common},$$ and {\em azygetic} if they have six lines in
common, for example: $$N\ \ \hbox{and} \ \ N_{456}\ \ \hbox{have } a_1,\
a_2,\ a_3,\ \ b_4,\ b_5,\ b_6\ \hbox{in common}.$$ Two azygetic double
sixes have six lines in common and contain 12 other lines; these 12 lines
form another double six, azygetic with respect to both, for example $N,\
N_{123},\ N_{456}$.  Such triples are refered to as triples of azygetic
double sixes or, because of the interpretation in terms of tritangents, a
trihedral pair. Each double six is syzygetic to 15 others, forming 270 such
pairs, and azygetic to 20 others, forming 120 triples.  Our notation for
the 120 triples are:
\begin{equation}\label{eQ1.2}
\begin{array}{rclr}
  \{ijk\} & = & <N,N_{ijk},N_{lmn}>, & (10) \\ \{ij.jk\} & = &
  <N_{ij},N_{ik},N_{jk}>, & (20) \\ \{ij.kl\} & = &
  <N_{ij},N_{ikl},N_{jkl}> & (90).\\
\end{array}
\end{equation}
We recognize these as the trihedral pairs of (\ref{eB2.1}) under the
correspondence
$$\left[\begin{array}{ccc} a_i & b_j & c_{ij} \\ b_k & c_{jk} & a_j \\
    c_{ik} & a_k & b_i \end{array}\right] \llra
<N_{ij},N_{ik},N_{jk}>,\quad \left[\begin{array}{ccc} c_{il} & c_{jm} &
    c_{kn} \\ c_{mn} & c_{ik} & c_{jl} \\ c_{jk} & c_{ln} & c_{im}
  \end{array}\right] \llra <N,N_{ijk},N_{lmn}>, $$
$$\left[\begin{array}{ccc} a_i & b_j & c_{ij} \\ b_l & a_k & c_{kl} \\
    c_{il} & c_{jk} & c_{mn} \end{array}\right] \llra
<N_{ij},N_{ikl},N_{jkl}>.$$ Hence the triads of trihedral pairs discussed
there are expressed in condensed form as follows:
\begin{equation}\label{eQ1.3}\begin{minipage}{14cm}
    $$ [ijk.lmn]=\left[ \matrix{ N_{ij} & N_{jk} & N_{ik} \cr N_{lm} &
        N_{mn} & N_{ln} \cr N & N_{ijk} & N_{lmn} \cr }
    \right],\hspace{1cm} (10)$$

    $$ [ij.kl.mn]=\left[ \matrix{ N_{ij} & N_{ikl} & N_{jkl} \cr N_{kl} &
        N_{kmn} & N_{lmn} \cr N_{mn} & N_{nij} & N_{mij} \cr } \right],
    \hspace{1cm} (30).$$
\end{minipage}
\end{equation}
The group of incidence preserving permutations of the 27 lines, a group of
order 51840, can be generated by the following six operations:
$$(i,i+1),\ i=1,...,5:\hbox{transposition of the indices},$$ and
$$ (123) :\hbox{map } N \mapsto N_{123},$$ and the graph of this
presentation is shown in Figure \ref{Figure4}.
This is the graph whose vertices correspond to generators, two vertices A,
B being connected if ABA=BAB and not connected if AB=BA.

\subsection{Roots}

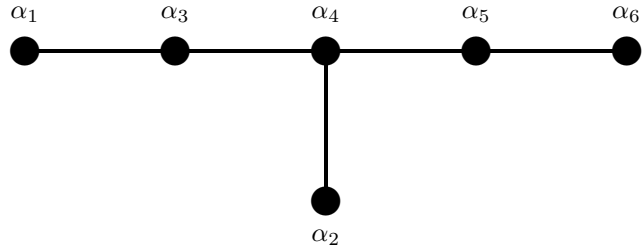
\begin{figure}[hbt]
  $$ \unitlength0.8cm
\begin{picture}(10.5,3.5)

  \put(0.25,3.0){\circle*{0.5}} \put(2.75,3.0){\circle*{0.5}}
  \put(5.25,3.0){\circle*{0.5}} \put(7.75,3.0){\circle*{0.5}}
  \put(10.25,3.0){\circle*{0.5}} \put(5.25,0.5){\circle*{0.5}}

  \linethickness{0.4mm}

  \put(0.25,3.0){\line(1,0){2.5}} \put(2.75,3.0){\line(1,0){2.5}}
  \put(5.25,3.0){\line(1,0){2.5}} \put(7.75,3.0){\line(1,0){2.5}}
  \put(5.25,0.5){\line(0,1){2.5}}

  \put(0.25,3.6){\makebox(0,0){$\ga_1$}}
  \put(2.75,3.6){\makebox(0,0){$\ga_3$}}
  \put(5.25,3.6){\makebox(0,0){$\ga_4$}}
  \put(7.75,3.6){\makebox(0,0){$\ga_5$}}
  \put(10.25,3.6){\makebox(0,0){$\ga_6$}}
  \put(5.25,-0.1){\makebox(0,0){$\ga_2$}}

\end{picture}
$$
\caption[The Dynkin diagram of $E_6$]{\label{Figure5}\small
  The Dynkin diagram of the Weyl group of $E_6$}
\end{figure}

Let $\tt$ be a maximal abelian subalgebra of the compact Lie algebra
$\ee_{6,u}$ over $\fR$, i.e., $\tt \isom \fR^6$.  Let $x_1,...,x_6$ be
coordinates such that the root forms of $E_6$ are:
\begin{eqnarray*}
  (40) & & \pm(x_i\pm x_j), \hspace{1cm} 1\leq i <j\leq 5 \\ (32) & & \pm{1
    \over 2}(\pm x_1\pm x_2\pm x_3\pm x_4 \pm x_5 + x_6), \ \ \hbox{even
    number of ``$-$'' signs inside the parenthesis.}
\end{eqnarray*}
(Note that in Bourbaki notation, our variables $x_i=\ge_i,\ i=1,..,5$,
while our coordinate $x_6$ is denoted $\ge_8-\ge_7-\ge_6$ there).  The 36
positive root forms are given by $\pm x_i+x_j$ and ${1\over 2}(\pm x_1\pm
x_2\pm x_3\pm x_4 \pm x_5 + x_6)$, and they correspond to the 36 double
sixes of the 27 lines on a cubic surface.  We use the following notations
for these forms
\begin{equation}\label{eQ2.1}
\begin{array}{lcl} h & = & {1 \over 2}(x_1+...+x_6), \\
  h_{1j} & = & x_{j-1}-{1\over 2}(x_1+...+x_5-x_6),\ \ j=2,...,6 \\ h_{jk}
  & = & -x_{j-1}+x_{k-1},\ \ 1\neq j<k \\ h_{1jk} & = & x_{j-1}+x_{k-1}, \
  \ j,k=2,...,6 \\ h_{jkl} & = &+x_{j-1}+x_{k-1}+x_{l-1}-{1 \over
    2}(x_1+...+x_5-x_6),\ \ j,k,l\neq 1. \\
\end{array}
\end{equation}
The Weyl group of $E_6$ is generated by the reflections on these 36
hyperplanes; we denote these reflections by $s,\ s_{ij},$ and $s_{ijk}$.
As a system of simple roots we take :

\begin{equation}\label{eQ2.q}
\begin{array}{cclcl}
  \ga_1 & = & -{1 \over 2}(-x_1+...+x_5-x_6) & = & h_{12} \\ \ga_2 & = &
  x_1+x_2 & = & h_{123} \\ \ga_3 & = & -x_1+x_2 & = & h_{23} \\ \ga_4 & = &
  -x_2+x_3 & = & h_{34} \\ \ga_5 & = & -x_3+x_4 & = & h_{45} \\ \ga_6 & = &
  -x_4+x_5 & = & h_{56}. \\
\end{array}
\end{equation}
Then the Dynkin diagram is as shown in Figure \ref{Figure5}; we recover
Figure \ref{Figure4} by replacing $\ga_i$ by the corresponding {\it
  reflection} $s, s_{ij}, s_{ijk}$ on the hyperplanes where $h, h_{ij},
h_{ijk}$, respectively, vanish. This shows clearly the isomorphism of
$W(E_6)$ and the group of the permutations of the 27 lines, $$Aut(\cL)\isom
W(E_6).$$

The action of the reflections on the root forms can be described as
follows:\par \vspace{.5cm} \renewcommand{\arraystretch}{1.4}
$\begin{array}{ccc|ccc|ccc|ccc} s(h_{ij}) & = & h_{ij} & s(h_{ijk}) & = &
  h_{lmn} & s_{ijk}(h) & = & h_{lmn} & s_{ij}(h) & = & h \\ \hline
  s_{ijk}(h_{lmn}) & = & h & s_{ijk}(h_{kmn}) & = & h_{kmn} &
  s_{ijk}(h_{jkn}) & = & h_{in} & s_{ijk}(h_{ln}) & = & h_{ln} \\ \hline
  s_{ijk}(h_{kn}) & = & h_{ijn} & s_{ijk}(h_{jk}) & = & h_{jk} &
  s_{ij}(h_{klm}) & = & h_{klm} & s_{ij}(h_{jlm}) & = & h_{ilm} \\ \hline
  s_{ij}(h_{ijm}) & = & h_{ijm} & s_{ij}(h_{ij}) & = & h_{ij} &
  s_{ij}(h_{jk}) & = & h_{ik} & s_{ij}(h_{lk}) & = & h_{lk} \\ \hline
\end{array}$\renewcommand{\arraystretch}{1}

\subsection{Vectors}
The Killing form of $E_6$, a quadratic invariant, can be calculated as the
sum of the squares of all roots, and evaluates to (a constant times):
$$I_2=x_1^2+x_2^2+x_3^2+x_4^2+x_5^2+{1 \over 3}x_6^2.$$ With respect to the
Killing form we have the vectors dual to the root forms:
\begin{equation}\label{eQ3.1}
\begin{array}{lcll}
  H & = & {1 \over 2}(1,1,1,1,1,3); & \\ H_{1j} & = & -{1 \over
    2}(1,...,-1,...,-3), & 1 \hbox{ in the $(j-1)^{st}$ spot},\ j=2,\ldots,6;
  \\ H_{jk} & = & -{1 \over 2}(0,..1,..,-1,..,0), & \pm 1 \hbox{ in the
    $(j-1)^{st}, (k-1)^{st}$ spot}, 1<j<k\leq 6; \\ H_{1jk} & = & {1 \over
    2}(0,..1,..,1,..,0), & 1 \hbox{ in the $(j-1)^{st}, (k-1)^{st}$ spot},\
  1<j<k\leq 6; \\ H_{jkl} & = & -{1 \over 2}(1,-1..,-1,..,-1,..,-3), & 1's
  \hbox{ in the $(j-1),\ (k-1),\ (l-1)$ spots}, \\ & & & \ 1<j<k<l\leq 6; \\
\end{array}
\end{equation}
which may be thought of as the root vectors (of the positive roots; the
negative roots have a ``$-$'' sign in front).

\begin{table}
\caption{\label{table23}
  The arrangement in ${\Bbb P}^5$} \renewcommand{\tabcolsep}{4pt}
\begin{tabular}{|l|r|rr|rrr|rrrrr|rrrr|} \hline $N(\cO)$ & $A_1$ & $A_1^2$
  & $A_2$ & $A_1^3$ &$A_{1,2}$ & $A_3$ & $A_{1^2,2}$ & $A_2^2$ & $A_{1,3}$
  & $A_4$ & $D_4$ & $A_{1,2^2}$ & $A_{1,4}$ & $A_5$ & $D_5$ \\ \hline

  $\cO$ & k & $t_2(3)$ & $t_3(3)$ & $t_3(2)$ &$t_4(2)$ & $t_6(2)$ &
  $t_5(1)$ & $t_6(1)$ & $t_7(1)$ & $t_{10}(1)$ & $t_{12}(1)$ & $t_7$ &
  $t_{11}$ & $t_{15}$ & $t_{20}$ \\ \hline $\#$ & 36 & 270 & 120 & 540 &
  720 & 270 & 1080 & 120 & 540 & 216 & 45 & 360 & 216 & 36 & 27 \\ \hline
  $t(4)$ & 1 & 15 & 10 & 45 & 80 & 45 & 150 & 20 & 105 & 60 & 15 & 70 & 66
  & 15 &15 \\ \hline $t_2(3)$ & & 1 & 0 & 6 & 8 & 3 & 28 & 4 & 18 & 12 & 3
  & 20 & 20 & 6 & 7 \\ $t_3(3)$ & & & 1 & 0 & 6 & 9 & 9 & 2 & 18 & 18 & 6 &
  6 & 18 & 6 & 9 \\ \hline $t_3(2)$ & & & & 1 & 0 & 0 & 6 & 0 & 3 & 0 & 1 &
  6 & 6 & 1 & 3 \\ $t_4(2)$ & & & & & 1 & 0 & 3 & 1 & 3 & 3 & 0 & 4 & 6 & 3
  & 3 \\ $t_6(2)$ & & & & & & 1 & 0 & 0 & 2 & 4 & 2 & 0 & 4 & 2 & 5 \\
  \hline $t_5(1)$ & & & & & & & 1 & 0 & 0 & 0 & 0 & 2 & 2 & 0 & 1 \\
  $t_6(1)$ & & & & & & & & 1 & 0 & 0 & 0 & 3 & 0 & 3 & 0 \\ $t_7(1)$ & & &
  & & & & & & 1 & 0 & 0 & 0 & 2 & 1 & 1 \\ $t_{10}(1)$ & & & & & & & & & &
  1 & 0 & 0 & 1 & 1 & 2 \\ $t_{12}(1)$ & & & & & & & & & & & 1 & 0 & 0 & 0
  & 3 \\ \hline
\end{tabular}

\end{table}

As is well-known, there is also a set of 27 fundamental weights which form
an orbit of $W(E_6)$, namely:
\begin{equation}\label{eQ3.3}
\begin{array}{cclccl}
  a_1 & = & -{2 \over 3}x_6; & a_j & = & x_{j-1}-{1 \over 2}(x_1+...+x_5+{1
    \over 3}x_6); \\ b_1 & = & {1 \over 2}(x_1+...+x_5-{1 \over 3}x_6); &
  b_j & = & x_{j-1}+{1 \over 3}x_6; \\ c_{1j} & = & -x_{j-1}+{1 \over 3}x_6;
  & c_{ij} & = & -x_{j-1}-x_{i-1}+{1 \over 2}(x_1+...+x_5-{1 \over 3}x_6)
  .\\
\end{array}
\end{equation}
These form the $W(E_6)$ orbit of the fundamental weights denoted
$\ove{\go}_1$ and $\ove{\go}_6$ in Bourbaki, which are just our $-a_1$ and
$b_6$, respectively.  Note that the following relation holds:
\begin{equation}\label{eQ3.4} \sum_{i=1}^6 a_i = -3h =-3(\sum_{i=1}^6 x_i)
  =-\sum_{i=1}^6 b_i.
\end{equation}
Also note that the $a_i$ and $b_i$ are related by
\begin{equation}\label{eQ3.5} b_i=a_i-{1\over 3}(a_1+\cdots +a_6).
\end{equation}
The corresponding vectors which are dual with respect to the Killing form
are:
\begin{equation}\label{eQ3.2}
\begin{array}{cclcclr}
  A_1 & = & (0,...,0,-2); & A_j & = & {1 \over 2}(-1,...,+1,..,-1) & +1\
  \hbox{in the $j-1$ spot;} \\ B_1 & = & {1 \over 2}(1,...,1,-1); & B_j & =
  & (0,...,1,..,1) & +1\ \hbox{in the $j-1$ spot;} \\ C_{1j} & = &
  (0,..,-1,...,1); & C_{ij} & = & {1 \over 2}(1,..,-1,..,-1,..,-1) & -1\
  \hbox{in the $j-1$, $i-1$ spots.} \\
\end{array}
\end{equation}

\subsection{The arrangement defined by $W(E_6)$}
The 36 hyperplanes in $\fP^5$ defined by the vanishing of the 36 root forms
form the arrangement $\cA({\bf E_6})$ of (\ref{e108.2}).
For later reference we give the combinatorial data of
the arrangement here. We denote as in (\ref{e109.1}) a $\fP^m$ through
which $k$ of the hyperplanes pass by $t_k(m)$. For the normalisers we use
the notation $A_{i^k,j^l}$ for $A_i^k\times A_j^l$. The data of the
arrangement is given in Table \ref{table23}.

\subsection{Special Loci}
In Table \ref{table24} we give a list of special loci which will be
particularly important in what follows, so we give a brief description of
each.
\begin{table}[htb]
\caption{\label{table24} Special loci in ${\Bbb P}^5$ \hfill}
$$\begin{array}{|r|l| l| c| l|}\hline \# & \hbox{space} & \hbox{Symmetry} &
  N(\cO) & \hbox{notation in Table \ref{table23}} \\
  \hline \hline 36 & \fP^4 & A_5 &
  A_1 & - \\ \hline 120 & \fP^3 & D_4 & A_2 & t_3(3) \\ \hline 120 &
  \fP^1 & A_2 & A_2\times A_2 & t_6(1) \\ \hline 216 & \fP^1 & A_2 & A_4 &
  t_{10}(1) \\ \hline 45 & \fP^1 & A_1 & D_4 & t_{12}(1) \\ \hline 36 &
  \hbox{point} & - & A_5 & t_{15} \\ \hline 27 & \hbox{point} & - & D_5 &
  t_{20} \\ \hline \end{array}$$
\end{table}

\subsubsection{36 $\fP^4$'s}\label{i1}
In each of the 36 hyperplanes given by the vanishing of one of the 36 forms
(\ref{eQ2.1}), $h$ say, the induced group is $\gS_6$, and as a reflection
group on $\fP^4$ it defines a projective arrangement of 15 planes; since
each double six is syzygetic to 15 and azygetic to 20 others, there are 15
hyperplanes through which one of the other 35 intersect $h$, and ten planes
through which two others of the 35 meet $h$. We immediately recognize this
geometry as that in $\fP^4$ discussed in the first part of the paper. The
15 hyperplanes are the 15 $\cH_{ij}$ of (\ref{e112b.3}), each of which cuts
out three planes on $\cS_3$, and the ten are the hyperplanes mentioned in
(\ref{e117a.1a}) and Proposition \ref{p120.1}. These in turn are the dual
hyperplanes to the ten nodes on $\cS_3$.

\subsubsection{120 $\fP^3$'s}\label{i2}
These $\fP^3$'s correspond to the 120 triples
 of azygetic double sixes,
i.e., each is cut out by three of the 36 hyperplanes of \ref{i1}. In each
such hyperplane, these $\fP^3$'s correspond to the ten hyperplanes in $h$
just mentioned, given by the $K_{ijk}$ of (\ref{e117a.1a}). Each of these
contains 15 planes, and one can check that these are just the faces and
symmetry planes of a cube. The six lines in $K_{ijk}$ which are the
singular locus $\cI_4\cap K_{ijk}$, are easily identified with the six
12-fold lines $t_{12}(1)$ which are contained in $K_{ijk}$\footnote{The
  arrangement is $\cA({\bf D_4})$, minus the plane at infinity. Of the 16
  $t_3(1)$ of $\cA({\bf D_4})$ (see (\ref{e109.2})), four lie in the plane
  at infinity.}, and the nine points $t_{20}$ contained in $K_{ijk}$ are
the intersection points of those six lines\footnote{Likewise, nine of the
  12 $t_6$ of $\cA({\bf D_4})$ lie in the plane at infinity}.  Equations of
the 120 $\fP^3$'s are given by a triple of azygetic double sixes, e.g., by
$<h,h_{ijk},h_{lmn}>$.

\subsubsection{120 $\fP^1$'s} \label{i3} The
120 lines correspond exactly to $A_2$ subroot systems, each containing
three (positive) roots, so that each line contains three of the 36 points.
The 120 lines are determined as follows.  Consider a triad of triples of
azygetic double sixes and the corresponding matrix of linear forms (see
(\ref{eQ1.3})), say
$$[ijk.lmn]=\left[ \matrix{ h_{ij} & h_{jk} & h_{ik} \cr h_{lm} & h_{mn} &
    h_{ln} \cr h & h_{ijk} & h_{lmn} \cr } \right].$$ Taking the ideal
defined by the vanishing of two rows defines the corresponding line, i.e.,
\begin{equation}\label{eQ5.1}
\begin{array}{lcl}
  L_{\{ij.jk\}} & = & <h_{lm}, h_{mn}, h_{ln}, h, h_{ijk}, h_{lmn}> \\
  L_{\{lm.mn\}} & = & <h_{ij}, h_{jk}, h_{ik}, h, h_{ijk}, h_{lmn}> \\
  L_{\{ijk\}} & = & <h_{ij}, h_{jk}, h_{ik},h_{lm}, h_{mn}, h_{ln}>.\\
\end{array}
\end{equation}
Each of the 120 lines contains three of the nodes, so for example,
\begin{equation}\label{eQ5.2}
  H_{ij},\ H_{jk},\ H_{ik} \in L_{\{ij.jk\}}.
\end{equation}
There are 40 such triples of the 120 lines, which have the characterising
property that they span $\fP^5$. These correspond to subroot systems of the
type $A_2 \times A_2 \times A_2$, where all three copies are orthogonal to
one another.  Note that given an $A_2$ subroot system, there is a unique
$A_2 \times A_2$ subroot system orthogonal to it. Thus the $A_2$ subroot
system is defined by the vanishing of the six root forms of the
complementary $A_2 \times A_2$. There are 120 of each of both types of
subroot systems. Summing up, there are six of the 36 hyperplanes passing
through each of these 120 lines while each such line contains three of the
36 nodes.

The {\it induced arrangement} is as follows.
Blowing up along the line introduces
an exceptional $\fP^3$ over each point of the line; the intersection of it
with the proper transforms of the six planes passing through it is the
induced arrangement. It is of type $A_2\times A_2$, i.e., is given by two
skew lines in $\fP^3$ and two sets of three hyperplanes through each line.

\subsubsection{216 $\fP^1$'s}\label{i4}
Consider a pair of {\em skew} lines, say $a_1, a_2.$ There is a unique
double six containing the given pair as a column, e.g.,
$$N_{12}=\left[ \matrix{ a_1 & b_1 & c_{23} & c_{24} & c_{25} & c_{26} \cr
    a_2 & b_2 & c_{13} & c_{14} & c_{15} & c_{16} \cr } \right]. $$ There
are 216 lines in $\fP^5$ which join points such as $A_1, A_2, H_{12}$ (see
(\ref{eQ3.1}) and (\ref{eQ3.2})).  The
ideal of these 216 lines is generated by 24 sextics, forming the
irreducible representation denoted 24$_p$ in \cite{BL}. We can exhibit
these sextics explicitly, as follows.  The 216 lines are given by the
equations:
\begin{equation}\label{eQ7.1}
\begin{array}{lclc}
  <A_i,A_j,H_{ij}> & = & <h_{kl}|k,l\neq i,j; h_{ijk}> & (15) \\
  <B_i,B_j,H_{ij}> & = & <h_{kl}|k,l\neq i,j; h_{klm}|k,l,m\neq i,j> & (15)
  \\ <A_i,B_i,H> & = & <h_{kl}|k,l\neq i> & (6) \\ <A_i,C_{jk},H_{lmn}> & =
  & <h_{jk}, h_{\gl \gm}|\gl,\gm\neq i,j,k; h_{ij\gl}, h_{ik\gl}|\gl\neq
  i,j,k> & (60) \\ <B_i,C_{jk},H_{ijk}> & = & <h_{jk}, h_{\gl
    \gm}|\gl,\gm\neq i,j,k; h_{\gl\gm\gn}|\gl\neq i,j,k,\gm\neq i,k,
  \gn\neq i> & (60) \\ <C_{ik},C_{jk},H_{ij}> & = &
  <h_{ijm},h_{mn},h,h_{kmn}|m,n\neq i,j,k> & (60), \\
\end{array}
\end{equation}
i.e., each is defined by the vanishing of ten of the $h$'s; these lines are
the $t_{10}(1)$ listed in the table of the arrangement. We claim the
sextics are the products of the six root forms of an $A_2\times A_2$
subroot system. To see this, pick one, say $\Phi=h_{12}\cdot h_{13}\cdot
h_{23}\cdot h_{45}\cdot h_{46}\cdot h_{56}$. It will suffice to check that
for any of the 216 lines listed in (\ref{eQ7.1}), at least one of the
hyperplanes on the right hand side is among the set $h_{12}, h_{13},
h_{23},h_{45}, h_{46}, h_{56}$. This is at most a tedious, but
straightforward task.

The dual $\fP^3$'s, which are defined by the vanishing of the forms which
are dual to the points of
the left-hand sides, each {\em contain} the ten points which are
 dual to the forms on
the right, for example
$$P_{<A_i,A_j>}=\{a_i=a_j=h_{ij}=0\}\ni (H_{kl},H_{ijk}).$$

The induced arrangement over each line is the arrangement $\cA(W({\bf
  A_4}))$ of (\ref{e109.2}). Among the ten hyperplanes defining the line,
say $<A_1,B_1,H>$, there are ten triples of azygetic double sixes,
$\{ij.jk\}$ in (\ref{eQ1.2}), with $i\neq 1$: $\{23.34\},\ \{23.35\},\
\{23.36\},\ \{24.45\},\ \{24.46\},\ \{25.56\},\ \{34.45\},\ \{34.46\},\
\{35.56\}$, and $\{45.56\}$, and these determine the ten $t_3(1)$ of
(\ref{e109.2}).

\subsubsection{45 $\fP^1$'s}\label{i5}
The 45 lines are the lines joining the 27 points of (\ref{eQ3.2}) in
threes, for example, $$L_{(ij)}=<A_i,B_j,C_{ij}>.$$ These lines are defined
by the vanishing of 12 of the $h$'s, so for example
\begin{equation}\label{eQ9.1}
  L_{(12)}=<h_{34},h_{35},h_{36},h_{45},h_{46},h_{56},h_{234},h_{235},h_{236},
  h_{245},h_{246},h_{256}>;
\end{equation}
these are the hyperplanes corresponding to the 12 double sixes {\em not}
containing any of $a_1, b_2, c_{12}$.

The induced arrangement is $\cA(W({\bf D_4}))$, with 12 planes
corresponding to the 12 hyperplanes through the line. Again there will be
$t_2(1)$'s and $t_3(1)$'s, corresponding to triples of azygetic double
sixes (respectively to syzygetic double sixes).

The ideal of the 45 lines is generated by 15 quartics which form the
irreducible representation denoted 15$_q$ in \cite{BL}. It is easy to see
that this space of quartics is generated by a product of four pairwise
azygetic $h$'s, for example by $h_{24}\cdot h_{124}\cdot h_{35}\cdot
h_{135}$. In fact, each hyperplane of type $h_{ij}$ contains the 15 lines
numbered (like the tritangents) by: \par

$$\begin{array}{lccl} (kl) & \hbox{for} & k,l\neq i,j & \hbox{(12 of
    these)}\\ (ij.kl.mn) & \hbox{for} & k,l,m,n \neq i,j\ & \hbox{(3 of
    these)}
\end{array}$$ \\
while the hyperplanes of type $h_{ijk}$ contain the 15 lines numbered by:
$$\begin{array}{lccl} (mn) & \hbox{for} & n=i\ or \ j, m\neq i,j & \hbox{(9
    of these)}\\ (il.jm.kn) & & & \hbox{(6 of these).}
\end{array}$$\\
It is now easy to check that every line is contained in at least one of the
four hyperplanes. Alternatively we can argue as follows: each $h$ contains
15 of the lines; there are six $\fP^3$'s which are the intersection of two
of the four, three of which are contained in each $h$. These three meet in
a common line in the $h$, so the number of lines contained in the union is:
$4\cdot(15-7)+2\cdot 6+1=45$, where the 7=number of lines in the union of
the three $\fP^2$'s in each $h$, 2=3-1 is the number of lines in each such
$\fP^2$, not in the others, and one is the common line. Note that this is
the Macdonald representation corresponding to the four roots of an
$A_1\times A_1\times A_1\times A_1\inn D_4$ subroot system. Five of these
lines meet at each of the 27 points, corresponding to the five tritangents
through each of the 27 lines.

\subsubsection{36 points}\label{i6}
These are the 36 points (\ref{eQ3.1}) dual to the 36 hyperplanes of
\ref{i1}. The induced arrangement is of course just the arrangement in
$\fP^4$ above. There are 15 hyperplanes passing through each of the 36
points, and these are just the hyperplanes which are coded by the double
sixes which are syzygetic to the one with the notation of the point as in
(\ref{eQ3.1}). So, for example, the 15 $\fP^4$'s through the point $H$ are
the 15 $h_{ij}$.

These points correspond to ($\pm$) the roots of $E_6$. The orthogonal
complement in $\fR^6$ of the root $\ga$ is projectively equivalent to the
{\em dual} hyperplane to the point. For example, $H$ is dual to $h$, and
one of the hyperplanes $P$ will contain $H$ $\iff$ the dual point $p$ is
contained in $h$. The ideal of the 36 points is generated by 20 cubics,
forming the irreducible representation of $W(E_6)$ denoted 20$_p$ in
\cite{BL}. We can find these cubics explicitly as follows. Consider the
hyperplanes $a_1, b_2, c_{12}$ corresponding to a tritangent. From Table
\ref{table23} above we see that each of these hyperplanes contains 20 of
the 36 points (actually, the table contains the dual information: there are
20 of the $h_{ij}$, etc, passing through each of the 27 points), and the
$\fP^3$ which is the common intersection of these three contains 12 of the
36 (the dual information is contained in the table: the 45 lines are
12-fold lines).  Hence the product $a_1\cdot b_2\cdot c_{12}$ contains
3.(20-12)+12 = 36 of the 36 points.

Through each of the 36 points, also 15 of the 27 hyperplanes of
(\ref{eQ3.3}) pass, corresponding to the 15 lines {\em not} contained in
the double six whose notation the point has. For example, the point
$H={1\over 2}(1,\ldots,1,3)$ is contained in all the $c_{ij}$. In the
exceptional $\fP^4$ at the point, both sets of 15 hyperplanes (coming from
the 36, respectively 27 hyperplanes) {\it coincide}.

\subsubsection{27 points}\label{i7}
These are the points $A_i,B_i,C_{ij}$ of (\ref{eQ3.2}). There are 20 of the
36 hyperplanes meeting at each, so the induced arrangement is one of 20
$\fP^3$'s in $\fP^4$, and one sees easily that it is $\cA(W({\bf D_5}))$.
This arrangement is also induced in any of the hyperplanes $a_i,
b_i,c_{ij}$ of (\ref{eQ3.3}); we note that there are two kinds of $\fP^2$,
namely $t_2(2)$'s, corresponding to pairs of skew lines, and $t_3(2)$'s,
corresponding to tritangents. Since each line is contained in five
tritangents, there are five of the latter and 15 of the former (in each
$a_i$, etc.). These 15 form an arrangement of type $\cA(W({\bf A_5}))$ as
discussed above.  The ideal of these 27 points is generated by 30 cubics,
forming the irreducible representation denoted 30$_p$ in \cite{BL}. It is
easy to see that this space of cubics is generated by a product of three
members of a syzgetic triple as $h_{12}\cdot h_{13}\cdot h_{23}$, for
example: Each of the hyperplanes contains 15 of the 27, the $\fP^3$ which
is thier common intersection contains nine, so the union contains
3.(15-9)+9=27, or all of the points. Note that this is just the Macdonald
representation corresponding to the (3) roots of an $A_2$ subroot system.

We need, in addition to the above, certain information on the dual spaces.
\subsubsection{45 $\fP^3$'s} \label{i8}
Consider one of the 45 $\fP^3$'s which is dual to one of the 45 lines of
\ref{i5}; it is cut out by three of the forms (\ref{eQ3.3}), and can be
denoted as one of the 45 tritangents, for example, if $L_{(ij)}$ denotes
the line $<A_i,B_j,C_{ij}>$ as in (\ref{eQ9.1}), the dual $\fP^3$ may be
denoted by $l_{(ij)}$, and \begin{eqnarray}\label{eq3.1} l_{(ij)} & = &
  a_i\cap b_j \cap c_{ij} \\ l_{(ij:kl:mn)} & = & c_{ij}\cap c_{kl}\cap
  c_{mn}.\nonumber
\end{eqnarray}
Consider the $\fP^3$ $l_{(12)}$, given by $a_1=b_2=c_{12}=0$, or
$x_1=x_6=0$. Then one checks easily that the hyperplanes (\ref{eQ2.1})
reduce in $l_{(12)}$ to the arrangement $\cA({\bf F_4})$ of (\ref{e108.2}).
Considering how the 27 hyperplanes (\ref{eQ3.3}) intersect $l_{(12)}$, we
find that these restrict to the set of short roots, that is, give a
subarrangement of type $\cA({\bf D_4})$. See also Proposition \ref{pQ7.1}
below.

\subsection{Invariants}
Since the 27 forms (\ref{eQ3.3}) are (as a set) invariant under the Weyl
group the expression
\begin{equation}\label{invariantsQ3}
  I_k:=\sum_{i,j} \left\{ a_i^k + b_i^k + c_{ij}^k \right\},
\end{equation}
if non-vanishing, is an invariant of degree $k$. The ring of invariants of
$W(E_6)$ is generated by elements in degrees 2, 5, 6, 8, 9 and 12, which
can be taken to be $I_2,\ldots, I_{12}$.  We note that while $I_2$ and
$I_5$ are {\em unique}, the other invariants are only defined up to
addition of terms coming from lower degrees.

\section{The invariant quintic}
\subsection{Equation}
There is a unique (up to scalars) $W(E_6)$-invariant polynomial of degree
5.  Written with integer coefficients in the variables $x_i$ it is
\begin{equation}\label{eiq1.1} f(x_1,\ldots,x_6)=x_6^5-6x_6^3\gs_1(x)
  -27x_6\left(\gs_1^2(x)-4\gs_2(x)\right) -648\sqrt{\gs_5(x)},
\end{equation}
where $\gs_i(x)$ is the $i$th elementary symmetric polynomial of the
$x_1^2,\ldots,x_5^2$, so in particular $\sqrt{\gs_5(x)}=x_1x_2x_3x_4x_5$.
The polynomial $f(x)$ displays manifestly the $W(D_5)$-invariance of the
quintic. Under the change of variables from the $x_i$ to the $a_i$ of
(\ref{eQ3.3}), the equation $g(a)$ can be derived as follows. By
(\ref{eQ3.5}), we have $b_i=a_i-{1\over3}(a_1+\cdots+a_6)$, which by
equation (\ref{eQ3.4}) can be written $b_i=a_i-h$. The following trick was
shown to me by I. Naruki. Consider $\prod_{i=1}^6 a_i -\prod_{i=1}^6 b_i$.
This sextic divides the root $h$, and the quotient is $W(E_6)$-invariant.
To see this, calculate
\begin{eqnarray} a_1\cdot \cdots \cdot a_6 -(b_1\cdot \cdots \cdot b_6) & = &
  \prod a_i -\prod (a_i-h) \label{eiq1.2} \\ & = & \gs_6(a)-\left[
    \gs_6(a)-h\gs_5(a) +h^2\gs_4(a) -h^3 \gs_3(a) +h^4\gs_2(a)
    -h^5\gs_1(a)\right], \nonumber \end{eqnarray} where here $\gs_i(a)$ are
the elementary symmetric functions of the $a_i$.  Consequently,
\begin{eqnarray*} a_1\cdot \cdots \cdot a_6 -(b_1\cdot \cdots \cdot b_6) & = &
  h\left(\gs_5(a) -h\gs_4(a) +h^2\gs_3(a)-h^3\gs_2(a) +h^4\gs_1(a)\right),
\end{eqnarray*}
and since by (\ref{eQ3.4}) $h=-{1\over 3}\gs_1(a)$, this yields
\begin{eqnarray}
    g(a) & = & 81\gs_5(a)+27\gs_4(a)\gs_1(a) +9\gs_3(a)\gs_1^2(a)
    +3\gs_2(a)\gs_1^3(a) +\gs_1^5(a),\label{eiq1.3}
\end{eqnarray}
giving the
expression of the invariant quintic expressing manifestly the
$W(A_5)=\gS_6$-invariance.
\begin{definition} \label{diq1.1} The {\em invariant quintic} $\cI_5$ is the
  hypersurface of degree 5
  $$\cI_5:=\{x\in \fP^5\Big| f(x)=0\}\isom \{a\in \fP^5\Big| g(a)=0\},$$
  where the isomorphism is given by the change of coordinates from the
  $x_i$ to the $a_i$.
\end{definition}

\subsection{Singular locus}
Because of the equivalence of the coordinates $x_i,\ i=1,\ldots,5$, there
are essentially two different partial derivatives of $f$, namely
\begin{eqnarray}\label{eiq2.1} j_1: & = & {\del f \over \del x_1}\isom
  \cdots \isom {\del f \over \del x_5} \\ j_2: & = & {\del f \over \del
    x_6}. \nonumber \end{eqnarray} Calculating these forms gives
\begin{eqnarray}\label{eiq2.2} -{\del f \over \del x_i} & = & 12x_6^3x_i
  +54x_6x_i(x_j^2+x_k^2+x_l^2+x_m^2)+648 x_jx_kx_lx_m, \\ {\del f\over \del
    x_6} & = & 5x_6^4-18x_6^2\gs_1(x) +27(\gs_1^2(x) -4\gs_4(x)). \nonumber
\end{eqnarray} These are quartics with manifest $W(D_4)$ and $W(D_5)$
symmetry, respectively.
\begin{theorem}\label{tiq4.1} The singular locus of $\cI_5$ consists of the 120
  lines of \ref{i3}, which meet ten at a time in the 36 points of \ref{i6}.
\end{theorem}
{\bf Proof:} Consider first the hyperplane section $x_6=0$. Then the
equations to be solved are
\begin{eqnarray}\label{eiq4.1} x_ix_jx_kx_l & = & 0,\quad (i,j,k,l<6); \\
\label{eiq4.2} \gs_1^2(x)-4\gs_2(x) & = & 0. \end{eqnarray}
{}From (\ref{eiq4.1}) we get: two of the $x_i$ must vanish, say $x_4,\ x_5$,
and then (\ref{eiq4.2}) takes the form \begin{eqnarray}
  \left(x_1^2+x_2^2+x_3^2\right)^2
  -4\left(x_1^2x_2^2+x_1^2x_3^2+x_2^2x_3^2\right) & = & 0, \nonumber \\
  (x_1+x_2+x_3)(x_1-x_2-x_3)(x_1+x_2-x_3)(x_1-x_2+x_3) & = & 0
  \label{eiq4.3}
\end{eqnarray}
which splits into a product of four lines. Since the product $x_1\cdot
\cdots \cdot x_5=0$ is a coordinate simplex in $\fP^4=\{x_6=0\}$, it
follows that the 2-simplices of this simplex correspond to planes where two
of the coordinates vanish, hence there are ${ 5 \choose 2}=10 $ such
2-simplices; in each we have the four lines given by (\ref{eiq4.3}).  This
gives the 40 of the 120 lines contained in $x_6=0$. This implies that, in
the union of the 27 hyperplanes (\ref{eQ3.3}), the singular locus of
$\cI_5$ consists of 120 lines.

Suppose that $x_6\neq 0$. Then the simultaneous vanishing of the partials
${\del f \over \del x_i},\ i=1,\ldots,5$ imply that four of the $x_i$ must
vanish, say $x_2=x_3=x_4=x_5=0$. But the intersection of $\cI_5$ with the
line $\{x_2=x_3=x_4=x_5=0\}$ is given by
\begin{equation}\label{eiq4.4}
  x_6^5-6x_6^3x_1^2-27x_6x_1^4=x_6(x_6+i\sqrt{3}x_1)(x_6-i\sqrt{3}x_1)
  (x_6+3x_1) (x_6-3x_1),
\end{equation}
and the last two terms are the equations of $b_2$ and $c_{12}$, two other
of the 27 hyperplanes of (\ref{eQ3.3}). From this we conclude that any
singular point is contained in one of the 27 hyperplanes, and by the above,
that the singular locus of $\cI_5$ consists of the 120 lines, as stated.
\ende

The types of singularities are given by the following.
\begin{proposition}\label{piq5.1} The singularities of $\cI_5$ can be
  described as follows.
\begin{itemize}\item[i)] At a generic point $x\in$ one of the 120 lines, a
  transversal hyperplane section has an ordinary $A_1$-singularity, so the
  singularity is of type disc$\times A_1$.
\item[ii)] At one of the 36 intersection points $p$, the singularity has
  multiplicity 3, and the tangent cone is of the form
  $$s_5+s_4t+s_3t^2,$$ where $s_3$ is the cone over the Segre cubic,
  $s_4=s_3\cdot h_p$, where $h_p$ is the hyperplane of \ref{i1} dual to
  $p$, and $s_5$ is the cone over the intersection $\cI_5\cap h_p$.
\end{itemize}
\end{proposition}
{\bf Proof:} i) follows from consideration of generic hyperplane sections
of $\cI_5$, which are quintic threefolds with 120 isolated singularities,
so singularities worse than $disc\times A_1$ are impossible. ii) is just a
computation, done as follows. Suppose the point is
$p=H_{23}=(1,-1,0,0,0,0)$. Then inhomogenizing by setting
$t_i=x_i/x_1-tp_i$ (where $p_i$ denotes the $i$th coordinate of $p$),
inserting into the equation of $\cI_5$ gives the stated result. The fact
that $s_3$ is the cone over $\cS_3$ can be seen as follows. We can write
\begin{equation}\label{eiq6.1} f=s_5+s_3(th_p +t^2)
\end{equation}
and it follows that blowing up $\cI_5$ at $p$ is given by setting
$t=\infty$ and that the proper transform of $\cI_5$ in the exceptional
$\fP^4$ (of the blow up of $\fP^5$ at the point $p$) is a cubic $S_3=0$,
where $s_3=0$ is the cone over $S_3=0$.  Since there are ten of the 120
lines meeting at $p$, the resolving divisor of the blow up, which is a
cubic threefold, will have ten isolated singularities. As mentioned already
above, this implies the cubic threefold is isomorphic to $\cS_3$. One can
also see the 15 special hyperplane sections of $\cS_3$: these are the
proper transforms, under the blowing up of $p$, of the 15 of the 36
hyperplanes \ref{i1} passing through the point. The rest is calculation.
\ende

We have (using Macaulay) calculated the ideal $\scI(120)$ of the 120 lines,
and it turns out to be just the Jacobian ideal of $\cI_5$. I know of no
simple proof of this fact.

\subsection{Resolution of singularities} It turns out to be very easy to
desingularize $\cI_5$. By the proof of Proposition \ref{piq5.1}, we know
the 36 triple points can be resolved by blowing up each such point $p$. Let
$\grr^{(1)}:\cI_5^{(1)}\lra \cI_5$ denote this blow up of $\cI_5$. This has
the effect of seperating all 120 lines of $\cI_5$, and the singularities
along the lines are just $A_1$, by Proposition \ref{piq5.1}, i).  Hence
a desingularisation is achieved by resolving each of the 120 lines. There
are two possible ways to do this. First, one can blow up the lines in
${\fP^5}^{(1)}$=$\fP^5$ blown up in the 36 points, and take the proper
transform of $\cI_5^{(1)}$; this has the effect of replacing each singular
line by a quadric surface bundle, a $\fP^1\times \fP^1$-bundle, over the
line. Hence there are 120 exceptional divisors, each isomorphic to
$\fP^1\times \fP^1\times \fP^1$. We call this resolution of $\cI_5$ the
{\em big resolution} and denote it by $\~{\cI}_5$.  Secondly, we can take a
small resolution by blowing down one of the fiberings in the exceptional
$\fP^1\times \fP^1$ over a point of the line. In this way, each of the
singular lines is replaced by a $\fP^1$-bundle over the line, in other
words by a $\fP^1\times \fP^1$. We call this the {\em small resolution} and
denote it by $\cI_5^{(s)}$. Here no further (beyond the 36 on
$\cI_5^{(1)}$) exceptional divisors are introduced.

\begin{lemma}\label{liq7.1} The quintic $\cI_5$ has two resolutions, which we
  denote by $\~{\cI}_5$ and $\cI_5^{(s)}$. On $\~{\cI}_5$ there are 36+120
  exceptional divisors, 36 copies of the resolution of the Segre cubic, and
  120 copies of $\fP^1\times \fP^1\times \fP^1$. On $\cI_5^{(s)}$ there are
  only 36 exceptional divisors, each a small resolution of the Segre cubic.
\end{lemma}

\subsection{${\cal I}_5$ is rational}

Quite generally, in $\fP^5$, taking four $\fP^3$'s which meet only in
lines, there is a unique line which meets each of them and passes through a
given point $P\in \fP^5$, namely the line $<\ga,P>\cap<\gb,P>\cap
<\gg,P>\cap <\gd,P>$, if $\ga,\gb,\gg,\gd$ denote the $\fP^3$'s and
$<\ga,P>$ is the hyperplane spanned by $\ga$ and $P$.  Now let $P\in
\cI_5$, and choose four of the 15 of the 45 $\fP^3$'s through one of the
triple points $p$, such that the four $\fP^2$'s on $(\cS_3)_p$ meet each
other only in points; then the four $\fP^3$'s meet only in lines, and we
may apply this reasoning to conclude:
\begin{center}
  \parbox{14cm}{for each $P\in \cI_5-\{4\ \fP^3$'s \}, there is a unique
    line $L_p$, which joins $P$ and $\ga,\gb,\gg,\gd$.}
\end{center}
Then, fixing a generic hyperplane $F\inn\fP^5$. the line $L_p$ intersects
$F$ in a single point; we get a rational map:
$$\cI_5 -\ -\ -\ \rightarrow F$$
$$P\mapsto L_p\cap F.$$

We now carry out this argument to derive an explicit rationalisation. I am
indebted to B. v. Geemen for help in performing this. We {\em choose} four
convienient $\fP^3$'s which only meet in lines (although these do not all
pass through a point). The four $\fP^3$'s will be defined as follows:
\begin{equation}\label{eQR.1} \begin{array}{clccc}
    P_1 & = & \{a_1=b_4=0\} & = & \{l_1=m_1=0\} \\ P_2 & = & \{a_4=b_5=0\}
    & = & \{l_2=m_2=0\} \\ P_3 & = & \{a_5=b_6=0\} & = & \{l_3=m_3=0\} \\
    P_4 & = & \{c_{35}=c_{12}=0\} & = & \{l_4=m_4=0\}
\end{array}
\end{equation}
Letting $F$ be an auxilliary $\fP^4$ with homogenous coordinates
$(y_0:\ldots:y_4)$, the intersection of the line $<P_1,\ga>\cap \cdots \cap
<P_4,\ga>$ with $F$ is given by
$$y_0l_i-y_im_i=0,$$ which leads to \begin{equation}\label{eQR.2}
  \begin{array}{ccl} y_0 & = & m_1\cdots m_4 \\ y_i & = & l_i\cdot
    m_1\cdots \hat{m}_i \cdots m_4,
\end{array}\end{equation}
a system of quartics in $\fP^5$, which, when restricted to $\cI_5$, give
the rational map $\gff:\cI_5 - - \ra \fP^4(=F)$. Inverting the equations
for $x_1,\ldots,x_6$ we get
\begin{equation}\label{eQR.3} \begin{array}{ccl}
    x_1 & = & y_0^6y_1y_3+y_0^5y_1^2y_3-2y_0^6y_2y_3-y_0^5y_1y_2y_3
    +y_0^4y_1^2y_2y_3+2y_0^4y_1y_2y_3^2+2y_0^3y_1^2y_2y_3^2 \\ & & \
    -y_0^6y_1y_4-y_0^5y_1^2y_4-y_0^5y_1y_2y_4-2y_0^4y_1^2y_2y_4
    -y_0^3y_1^2y_2^2y_4-y_0^5y_1y_3y_4-y_0^4y_1^2y_3y_4 \\ & & \ \
    -2y_0^4y_1y_2y_3y_4-4y_0^3y_1^2y_2y_3y_4-2y_0^3y_1y_2^2y_3y_4
    -3y_0^2y_1^2y_2^2y_3y_4-2y_0^3y_1y_2y_3^2y_4 \\ & & \ \ \
    -3y_0^2y_1^2y_2y_3^2y_4-2y_0^2y_1y_2^2y_3^2y_4
    -2y_0y_1^2y_2^2y_3^2y_4+y_0^2y_1^2y_2y_3y_4^2 \\ & & \ \ \ \
    +y_0y_1^2y_2^2y_3y_4^2+y_0y_1^2y_2y_3^2y_4^2+y_1^2y_2^2y_3^2y_4^2 \\

    x_2 & = & y_0^6y_1y_3+y_0^5y_1^2y_3+y_0^5y_1y_2y_3+y_0^4y_1^2y_2y_3
    +2y_0^5y_1y_3^2+2y_0^4y_1^2y_3^2+2y_0^4y_1y_2y_3^2 \\ & & \
    +2y_0^3y_1^2y_2y_3^2-y_0^6y_1y_4-y_0^5y_1^2y_4-y_0^5y_1y_2y_4
    -2y_0^4y_1^2y_2y_4-y_0^3y_1^2y_2^2y_4-5y_0^5y_1y_3y_4 \\ & & \ \
    -5y_0^4y_1^2y_3y_4-6y_0^4y_1y_2y_3y_4-8y_0^3y_1^2y_2y_3y_4
    -2y_0^3y_1y_2^2y_3y_4-3y_0^2y_1^2y_2^2y_3y_4-4y_0^4y_1y_3^2y_4 \\ & & \
    \ \ -4y_0^3y_1^2y_3^2y_4-6y_0^3y_1y_2y_3^2y_4
    -7y_0^2y_1^2y_2y_3^2y_4-2y_0^2y_1y_2^2y_3^2y_4
    -2y_0y_1^2y_2^2y_3^2y_4+2y_0^5y_1y_4^2 \\ & & \ \ \ \

+2y_0^4y_1^2y_4^2+2y_0^4y_1y_2y_4^2+4y_0^3y_1^2y_2y_4^2+2y_0^2y_1^2y_2^2y_4^2
    +4y_0^4y_1y_3y_4^2+4y_0^3y_1^2y_3y_4^2 \\ & & \ \ \ \ \
    +6y_0^3y_1y_2y_3y_4^2+9y_0^2y_1^2y_2y_3y_4^2+2y_0^2y_1y_2^2y_3y_4^2
    +5y_0y_1^2y_2^2y_3y_4^2+2y_0^3y_1y_3^2y_4^2 +2y_0^2y_1^2y_3^2y_4^2 \\ &
    & \ \ \ \ \ \
    +4y_0^2y_1y_2y_3^2y_4^2+5y_0y_1^2y_2y_3^2y_4^2+2y_0y_1y_2^2y_3^2y_4^2
    +3y_1^2y_2^2y_3^2y_4^2 \\

    x_3 & = & 2y_0^7y_3+3y_0^6y_1y_3+y_0^5y_1^2y_3+2y_0^6y_2y_3
    +3y_0^5y_1y_2y_3+y_0^4y_1^2y_2y_3-2y_0^7y_4-3y_0^6y_1y_4 \\ & & \
    -y_0^5y_1^2y_4-2y_0^6y_2y_4-5y_0^5y_1y_2y_4-2y_0^4y_1^2y_2y_4
    -2y_0^4y_1y_2^2y_4-y_0^3y_1^2y_2^2y_4-2y_0^6y_3y_4 \\ & & \ \
    -3y_0^5y_1y_3y_4-y_0^4y_1^2y_3y_4-4y_0^5y_2y_3y_4
    -6y_0^4y_1y_2y_3y_4-2y_0^3y_1^2y_2y_3y_4-2y_0^3y_1y_2^2y_3y_4 \\ & & \
    \ \ -y_0^2y_1^2y_2^2y_3y_4+2y_0^3y_1y_2y_3^2y_4
    +y_0^2y_1^2y_2y_3^2y_4-2y_0^3y_1y_2y_3y_4^2 \\ & & \ \ \ \
    -y_0^2y_1^2y_2y_3y_4^2-2y_0^2y_1y_2^2y_3y_4^2
    -y_0y_1^2y_2^2y_3y_4^2-2y_0^2y_1y_2y_3^2y_4^2-y_0y_1^2y_2y_3^2y_4^2 \\
    & & \ \ \ \ \ -2y_0y_1y_2^2y_3^2y_4^2-y_1^2y_2^2y_3^2y_4^2 \\

    x_4 & = & 2y_0^7y_3+3y_0^6y_1y_3+y_0^5y_1^2y_3+y_0^5y_1y_2y_3
    +y_0^4y_1^2y_2y_3-2y_0^5y_1y_3^2-2y_0^4y_1^2y_3^2-2y_0^7y_4 \\ & & \
    -3y_0^6y_1y_4-y_0^5y_1^2y_4-3y_0^5y_1y_2y_4-2y_0^4y_1^2y_2y_4
    -y_0^3y_1^2y_2^2y_4-2y_0^6y_3y_4-y_0^5y_1y_3y_4 \\ & & \ \

+y_0^4y_1^2y_3y_4-2y_0^4y_1y_2y_3y_4-y_0^2y_1^2y_2^2y_3y_4+4y_0^4y_1y_3^2y_4
    +4y_0^3y_1^2y_3^2y_4+2y_0^3y_1y_2y_3^2y_4 \\ & & \ \ \
    +3y_0^2y_1^2y_2y_3^2y_4-2y_0^4y_1y_3y_4^2-2y_0^3y_1^2y_3y_4^2
    -2y_0^3y_1y_2y_3y_4^2-3y_0^2y_1^2y_2y_3y_4^2 \\ & & \ \ \ \
    -y_0y_1^2y_2^2y_3y_4^2-2y_0^3y_1y_3^2y_4^2-2y_0^2y_1^2y_3^2y_4^2
    -2y_0^2y_1y_2y_3^2y_4^2-3y_0y_1^2y_2y_3^2y_4^2-y_1^2y_2^2y_3^2y_4^2 \\

    x_5 & = & -y_0^6y_1y_3-y_0^5y_1^2y_3-y_0^5y_1y_2y_3-y_0^4y_1^2y_2y_3
    +y_0^6y_1y_4+y_0^5y_1^2y_4+2y_0^6y_2y_4+3y_0^5y_1y_2y_4 \\ & & \
    +2y_0^4y_1^2y_2y_4+2y_0^4y_1y_2^2y_4+y_0^3y_1^2y_2^2y_4
    +3y_0^5y_1y_3y_4+3y_0^4y_1^2y_3y_4+2y_0^4y_1y_2y_3y_4 \\ & & \ \
    +4y_0^3y_1^2y_2y_3y_4+2y_0^3y_1y_2^2y_3y_4+y_0^2y_1^2y_2^2y_3y_4
    +y_0^2y_1^2y_2y_3^2y_4-2y_0^5y_1y_4^2-2y_0^4y_1^2y_4^2 \\ & & \ \ \
    -2y_0^4y_1y_2y_4^2-4y_0^3y_1^2y_2y_4^2-2y_0^2y_1^2y_2^2y_4^2
    -2y_0^4y_1y_3y_4^2-2y_0^3y_1^2y_3y_4^2-2y_0^3y_1y_2y_3y_4^2 \\ & & \ \
    \ \
    -5y_0^2y_1^2y_2y_3y_4^2-3y_0y_1^2y_2^2y_3y_4^2-y_0y_1^2y_2y_3^2y_4^2
    -y_1^2y_2^2y_3^2y_4^2 \\ x_6 & = &
    -3y_0^6y_1y_3-3y_0^5y_1^2y_3-3y_0^5y_1y_2y_3-3y_0^4y_1^2y_2y_3
    +3y_0^6y_1y_4+3y_0^5y_1^2y_4+3y_0^5y_1y_2y_4 \\ & & \
    +6y_0^4y_1^2y_2y_4+3y_0^3y_1^2y_2^2y_4+3y_0^5y_1y_3y_4+3y_0^4y_1^2y_3y_4
    +6y_0^4y_1y_2y_3y_4+6y_0^3y_1^2y_2y_3y_4 \\ & & \ \
    +3y_0^2y_1^2y_2^2y_3y_4-3y_0^2y_1^2y_2y_3^2y_4+3y_0^2y_1^2y_2y_3y_4^2
    +3y_0y_1^2y_2^2y_3y_4^2+3y_0y_1^2y_2y_3^2y_4^2 \\ & & \ \ \
    +3y_1^2y_2^2y_3^2y_4^2
\end{array}
\end{equation}
a system of octics in $\fP^4$, yielding the rational map $$ \psi:\fP^4\lra
\cI_5.$$ These rational maps are morphisms outside of the base locus.
\begin{lemma}\label{lQR.1} The base locus of $\gff$ consists of the four
  $\fP^3$'s $P_1,\ P_2,\ P_3,\ P_4$. The base locus of $\psi$ is a surface
  of degree 32.\ende
\end{lemma}
The first statement is clear from construction, while the second is a
computation. We performed this with the help of Macaulay to calculate a
standard basis of the ideal; the base locus is the intersection of the six
octics.

\section{Hyperplane sections}
\subsection{Reducible hyperplane sections}
Consider the hyperplane section $H_5:=\cI_5\cap \{x_6=0\}$; it is the union
of five $\fP^3$'s which form a coordinate simplex $x_1\cdot x_2\cdot
x_3\cdot x_4\cdot x_5$ in the $\fP^4$ given by $\{x_6=0\}$. Now
$a_1=-2/3x_6$ and invariance implies that the 27 hyperplane sections
$a_i=0,\ b_i=0,\ c_{ij}=0$ all have the same property. Each such hyperplane
contains 40 of the 120 $\fP^1$'s, which meet six at a time in 20 of the 36
points.  Consider three lines in a tritangent, say $(a_1,\ b_2,\ c_{12})$.
These three hyperplanes pass through a common $\fP^3$, namely
$$l_{(12)}:=\{x_6=0,\ x_1=0\}.$$ Such $\fP^3$'s therefore correspond to the
tritangents and there are 45 such on $\cI_5$; these are just the 45
$\fP^3$'s of \ref{i8}.  Hence we have
\begin{proposition}\label{pQ7.1} The quintic $\cI_5$ contains 45 $\fP^3$'s,
  which are cut out by the 27 hyperplane sections (\ref{eQ3.3}), and each
  such hyperplane section meets $\cI_5$ in the union of five of the 45
  $\fP^3$'s. These can be numbered in terms of the tritangents of a cubic
  surface, i.e., for any 3 lines in a tritagent plane of a cubic surface,
  the corresponding hyperplanes of (\ref{eQ3.3}) intersect in a common
  $\fP^3$, and this $\fP^3$ lies on $\cI_5$.
\end{proposition}
Also the intersections of the 45 $\fP^3$'s can be described.  Each such
$\fP^3$ contains 16 of the 120 lines which meet in 12 of the 36 points;
these 12 points are the vertices of a triad of desmic tetrahedra. Consider
the $\fP^3$ $l_{(12)}$; the corresponding tritangent meets 12 others,
namely (13), (14), (15), (16), (32), (42), (52), (62), (12.34.56),
(12.35.46), (12.36.45) and (21), and the 12 $\fP^3$'s corresponding to them
meet $l_{(12)}$ in a $\fP^2$ (the generic intersection has dimension 1).
These 12 planes in $l_{(12)}$ form the arrangement $\cA({\bf D_4})$ in
$\fP^3$.

\subsection{Special hyperplane sections}
We now consider the intersections of $\cI_5$ with the 36 reflection
hyperplanes \ref{i1}. Take for example the reflection hyperplane $\{h=0\}$;
since $h$ is just a multiple of $\gs_1(a)$ (see (\ref{eQ3.4})), it follows
from the equation (\ref{eiq1.3}) that the intersection $\{h=0\}\cap \cI_5$
is a quintic hypersurface in $\fP^4$ with the equation:
\begin{equation}\label{eiq8.1} Q_1:=\{h=0\}\cap \cI_5 = \left\{
\begin{array}{c} \gs_1(a) = 0 \\ \gs_5(a) = 0 \end{array} \right.
\end{equation}
Comparing with the equation (\ref{e124.1}), we see that this is a copy of
the Nieto quintic! By symmetry, each of the 36 hyperplane sections is
isomorphic to this one, and we denote them by
\begin{equation}\label{eiq8.2} T=\{h=0\}\cap \cI_5,\
  T_{ij}=\{h_{ij}=0\}\cap \cI_5,\ T_{ijk}=\{ h_{ijk}=0\}\cap \cI_5.
\end{equation}
So we have:
\begin{proposition}\label{piq8.1} There are 36 copies of the Nieto quintic
  (\ref{e124.1}) on $\cI_5$.
\end{proposition}
We can determine the singular locus of these hyperplane sections,
independently of the discussion given in section \ref{section4.1}.  The
reflection hyperplane contains 20 of the 120 lines, which meet in 15 of the
36 points (corresponding to the 15 roots of an ${\bf A_5}$ subsystem), so
the quintic threefold has 20 singular lines, with 15 singular points of
multiplicity 3. In fact, the resolving divisor of each of these 15 points
is a four-nodal cubic surface, which is a hyperplane section of the Segre
cubic $\cS_3$ (see the discussion following Problem \ref{p132.1}, ii)).
Furthermore, recalling that there are ten of the 120 lines which pass
through the triple point which is {\em dual} to the given reflection
hyperplane, each such intersects the reflection hyperplane transversally,
giving the ten isolated ordinary double points on that quintic (see
Proposition \ref{p124.1}), and in some sense ``explains'' these isolated
singularities.

\subsection{Generic hyperplane sections}
A generic hyperplane section is a quintic threefold in $\fP^4$ with 120
nodes. This is a fascinating family of Calabi-Yau threefolds, which has a
beautiful geometric configuration associated with it, in some sense
``dual'' to the configuration of the 27 lines on a cubic surface.
\begin{proposition} Let $H\in  \fP^5$ be a generic hyperplane and let
  $Q_H=\cI_5\cap H$ be the hyperplane section. Then we have
\begin{itemize}\item[1)] There are 45 $\fP^2$'s on $Q_H$, which are cut out
  by 27 hyperplanes; these could appropriately be called {\em quintangent
    planes}.
\item[2)] The group of incidence preserving permutations of the 45
  $\fP^2$'s is $W(E_6)$; this is also the group of incidence preserving
  permutations of the 27 hyperplanes.
\item[3)] There are 36 hyperplane sections of $Q_H$, each of which is a
  20-nodal quintic surface.
\item[4)] The 120 nodes of $Q_H$ form an orbit under $W(E_6)$.
\end{itemize}\end{proposition}
{\bf Proof:} For any of the 45 $\fP^3$'s in $\cI_5$ and hyperplane section
$H$, it holds that $\fP^3\cap H=\fP^2\inn H\cap \cI_5=Q_H$, showing 1). The
second point is evident, and in a sense ``dual'' to the situation with
cubic surfaces.  We have seen that a special hyperplane section as in
(\ref{eiq8.2}) is isomorphic to the Nieto quintic and has 20 singular lines
in its singular locus; therefore any generic hyperplane section has exactly
20 nodes. 4) follows since the $W(E_6)$ orbit consisting of the 120 lines,
restricted to the hyperplane section is still an orbit. \ende We now
consider some of the invariants of the nodal quintic threefolds. Let $V$
denote a nodal quintic, $\hat{V}\lra V$ a small resolution and
$\tilde{V}\lra V$ a big resolution. Letting $s$ denote the number of nodes,
the betti numbers are
\begin{equation}\label{e8.3.1}\begin{minipage}{15cm}$$\begin{array}{ll}
      b_1(V)=1=b_1(\hat{V}), & b_2(\tilde{V})=1+d+s; \\
      b_2(V)=1+d=b_2(\hat{V}), & b_4(\tilde{V})=1+d+s; \\
      b_3(V)=b_3(V_t)-s+d, & b_3(\hat{V})=b_3(V)-s+d=b_3(\tilde{V}),
    \end{array}$$
\end{minipage}\end{equation}
where $V_t$ is a smooth hypersurface of same degree as $V$ and $d$ is the
{\it defect}. The defect may be calculated by the following result.
\begin{theorem}[\cite{W}, p.~27]\label{twerner}
  Let $V\in \fP^4$ be a nodal hypersurface of degree $n\geq 3$. Then
  $$\dim(\scP_{2n-5}(V))=\dim\left\{ \parbox{6cm}{homogenous polynomials of
      degree $2n-5$ in $\fP^4$, containing all nodes of $V$}\right\} =
  {2n-1 \choose 4} -s +d.$$
\end{theorem}
Applied to the case at hand, we need the dimension of the space of {\it
  quintics} vanishing at all the nodes. Clearly this is the degree five
component in the ideal of the 120 points. As we mentioned above, we {\it
  know} the ideal of the 120 lines (it is the Jacobian ideal of $\cI_5$
$\scJ ac(\cI_5)$), so we know also the ideal of the 120 points; it is the
restriction of $\scJ ac(\cI_5)$ to the hyperplane, generated by six
quartics.
\begin{proposition} The dimension of the space $\scP_5(Q_H)$ is 30.
\end{proposition}
{\bf Proof:} Each of the quartics (which are clearly independent for a
generic hyperplane $H$) of $\scJ ac(\cI_5)$ can be multiplied by any
hyperplane, giving a quintic which contains the 120 nodes. The set of
hyperplanes is $(\fP^5)^{\vee}$, so the dimension of $\scP_5(Q_H)$ is
$6\cdot 5=30$. \ende We can now apply Theorem \ref{twerner} to calculate
the defect $d$ for $Q_H$. The formula is $126-120+d=30$, from which is
follows that $d=24$.
\begin{corollary} The small resolutions $\hat{Q}_H$ of the quintic
  threefolds $Q_H$ have the following betti and Hodge numbers:
  $$\begin{array}{lll} b_2(\hat{Q}_H) = 25, & h^{1,1}=25, & \\
    b_3(\hat{Q}_H)=12=2+2h^{2,1}, & h^{2,1} =5, & e=2h^{1,1}-2h^{2,1}=40.
\end{array}$$
\end{corollary}
{\bf Proof:} Insertion of $d=24$ in (\ref{e8.3.1}).\ende In the well known
manner for Calabi-Yau threefolds the isomorphism $H^2(V,\gO^1)\isom
H^1(V,\Theta)$ identifies the Hodge space $H^{2,1}$ with the space of
infinitesimal deformations of $V$, $H^1(V,\Theta)$. This is by the above
five-dimensional, hence the moduli space of these 120 nodal quintics (a
Zariski open subset of $(\fP^5)^{\vee}$) is also a global space of
complex deformations of the small resolution. We can describe the space
$H^{2,1}$ more concretely as follows. Consider the space $\scP_5(Q_H)$; let
$\scJ\inn \scP_5(Q_H)$ be the subspace generated by the Jacobi ideal of
$Q_H$; since $Q_H$ has five partial derivatives, $\scJ$ is $5\cdot 5=25$
dimensional, and $\scJ$ cannot contribute to infinitesimal deformations, so
we have
$$H^{2,1}(\hat{Q}_H)\isom \scP_5(Q_H)/\scJ.$$

As a final remark consider the Picard group $\Pic(\hat{Q}_H)$ and the
orthocomplement of the hyperplane section $\Pic^0(\hat{Q}_H)$. Then
$rk_{\fZ}\Pic^0(\hat{Q}_H)=24$, and the 45 $\fP^2$'s give us privledged
representatives in $\Pic^0$; the 27 hyperplanes represent relations, so we
have an exact sequence
$$\fZ^{27}\lra \fZ^{45}\lra \Pic^0(\hat{Q}_H) \lra 1,$$ and the kernel is
six-dimensional. The sum sequence is then
\begin{equation}\label{e48a.1} 1 \lra \fZ^6 \lra \fZ^{27}\lra \fZ^{45} \lra
  \fZ^{24}\lra 1,
\end{equation}
and this is really dual to the sequence (\ref{eB3.2}) for cubic surfaces.

\begin{remark}
  The period map for this five-dimensional family of Calabi-Yau threefolds
  maps to the domain $\cD = Sp(6,\fR)/U(1)\times U(5)$. Note that any
  hyperplane passing through one of the 45\ $\fP^3$'s will intersect
  $\cI_5$ in the union of that $\fP^3$ and a residual quartic; clearly
  these constitute the set of cusps for the period map, i.e., on the 45
  lines in $(\fP^5)^{\vee}$ (the dual $\fP^5$) which parameterise the set of
  hyperplanes passing through one of the 45 $\fP^3$'s, the period map maps
  to the boundary of the domain $\cD$ above.  These 45 one-dimensional
  cusps meet in 27 points, i.e., zero-dimensional cusps, which correspond
  to the 27 hyperplane sections which split into the union of five
  $\fP^3$'s. But we can say more. Noting that, excepting the hyperplanes
  above, all hyperplane sections of $\cI_5$ are irreducible quintics, the
  worst that can happen is that the hyperplane passes through one of the 36
  triple points of $\cI_5$. We will see below that these are still
  Calabi-Yau (Proposition \ref{p158.1}), hence {\it not} contained in the
  boundary.
\end{remark}

\subsection{Tangent hyperplane sections}
We now consider the case of a hyperplane tangent to $\cI_5$ at a point
$p\in \cI_5$. In this case the section $Q_p$ aquires an additional node.
Note that the 121 nodes fall into two ``orbits'', one set of 120 on which
$W(E_6)$ acts as a permutation group, and the additional point $p$. For a
121-nodal quintic the same calculation as above gives $e(Q_p)=42$,
$h^{2,1}=4,\ h^{1,1}=25$.  It follows that the $H_4(Q_p,\fQ)$ is the same
as for $Q_x,\ x\in \fP^5$ generic.  The difference to the generic case is
in $H_3$, more precisely in $H^{2,1}$.  Indeed, we now require
$\scP_5(Q_p)$, that is, quintics through all 121 nodes, so as opposed to
the general case, we now only have, for each of the five quartics in the
Jacobi ideal of $Q_p$, since each contains $p$, a five-dimensional family
of quintics, as above. But for the quartics through the 120 nodes which are
{\it not} in the Jacobi ideal, we must take a hyperplane {\it through the
  point} $p$, so
\begin{proposition} For a 121-nodal quintic $Q_p$, $p\in \cI_5$, we have
  $\dim\scP_5(Q_p)=5\cdot 5 + 1\cdot 4 =29$. \ende
\end{proposition}
We can now apply Theorem \ref{twerner} to calculate the defect:
$$d=29-126+121 = 24.$$
\begin{corollary} The betti numbers for the small resolution $\hat{Q}_p$
  are
  $$\begin{array}{ll}b_2(\hat{Q}_p) =25, & h^{1,1} = 25, \\ b_3(\hat{Q}_p) =
    10, & h^{2,1}=4,\ \ \ \ e=42. \\
\end{array}$$ \end{corollary}
We remark that since $h^{1,1}$ is still 25, the sequence (\ref{e48a.1})
still holds for $Q_p$.

\section{Birational maps and the projection from a triple point}
\subsection{The cuspidal model}
First we recall the notations $\cI_5^{(1)}$ for the blow up of $\cI_5$ at
the 36 triple points, $\~{\cI}_5$ for the big resolution of $\cI_5$, and
$\cI_5^{(s)}$ for the small resolution. Note that on $\cI_5^{(1)}$, each of
the 120 lines has normal bundle $\cO(-2)^{\oplus 3}$, hence each line can
be blown down to an isolated singular point.
\begin{definition}\label{dq4.1} Consider the following birational
  transformation of $\cI_5$:
\begin{itemize}\item[i)] Blow up the 36 triple points,
  $\grr^{(1)}:\cI_5^{(1)} \lra \cI_5$;
\item[ii)] Blow down the proper transforms of the 120 lines to 120 isolated
  singularities, $\grr^{(2)}:\cI_5^{(1)}\lra \hat{\cI}_5$.
\end{itemize}
Step ii) defines the {\em cuspidal model} $\hat{\cI}_5$.
\end{definition}
This is a four-dimensional analogue of $\hat{\cN}_5$ of (\ref{e127.2}).
Indeed, for each of the 36 hyperplane sections of Proposition \ref{piq8.1},
the proper transform on $\hat{\cI}_5$ is isomorphic to $\hat{\cN}_5$:
\begin{lemma}\label{lq4.1} Let $T\isom \cN_5$ be one of the 36 special
  hyperplane sections of (\ref{eiq8.2}), and let $\hat{T}$ denote its
  proper transform on $\hat{\cI}_5$. Then $\hat{T}\isom \hat{\cN}_5$.
\end{lemma}
{\bf Proof:} Just check that the steps i) and ii) of Definition
\ref{dq4.1}, when restricted to $T$, coincide with those of
(\ref{e127.2}).\ende Let us mention that $\hat{\cI}_5$ ``looks like'' a
ball quotient too, at least assuming a positive answer to Problem
\ref{p132.1}.  We explain what ``looks like'' means in the following items.
\begin{itemize}
\item[I1] Each isolated singularity is resolved by a $\fP^1\times
  \fP^1\times \fP^1$; the arrangement induced in each by the proper
  transforms of the 36 hyperplanes and 36 exceptional divisors is a {\em
    product}, consisting of three fibres in each fibering (i.e., $\{\hbox{3
    points}\}\times \fP^1\times \fP^1,\ \fP^1\times \{\hbox{3
    points}\}\times \fP^1,\ \fP^1\times \fP^1\times \{\hbox{3 points}\}$,
  see \ref{i3}). Hence this can be covered in an equivariant way by
  $E_{\grr}\times E_{\grr}\times E_{\grr}$ (see Lemma \ref{l115.2}).
\item[I2] The proper transforms of the 36 hyperplane sections of
  Proposition \ref{piq8.1} are by Lemma \ref{lq4.1} isomorphic to
  $\hat{\cN}_5$, so, if the Problem \ref{p132.1} has an affirmative
  solution, these are ball quotients, with cusps being those isolated
  singularities of $\hat{\cI}_5$ which are contained in the given
  $\hat{T}$.
\item[I3] Consider the 45 $\fP^3$'s of Proposition \ref{pQ7.1}. These are
  (the proper transforms of) the 45 $\fP^3$'s of \ref{i8}.  These $\fP^3$'s
  are also ball quotients, in fact in two different ways.
\begin{itemize}\item[1)] There is a cover $Y\lra \fP^3$, branched over the
  arrangement $\cA(W({\bf D_4}))$ in $\fP^3$, which is a ball quotient.
  This example can be derived from the solution 4) of (\ref{e111d.1}) by
  means of the natural squaring map $m_2:\fP^3\lra \fP^3,\ (x_0:\ldots
  :x_3)\mapsto (x_0^2:\ldots : x_3^2)$. Then the arrangement $\cA(W({\bf
    D_4}))$ is the pullback under $m_2$ of the six symmetry planes of the
  tetrahedron in the arrangement $\cA(W({\bf A_4}))$, and pulling back the
  solution 4), we get the cover $Y\lra \fP^3$, branched along $\cA(W({\bf
    D_4}))$ (with branching degree 3 at each hyperplane), which is a ball
  quotient by a fix point free group.
\item[2)] There is a cover $Z\lra \fP^3$, branched along the arrangement
  $\cA(W({\bf F_4}))$ in $\fP^3$ (but not a Fermat cover), which is a ball
  quotient; this example is explained in \cite{hunt}, Thm.~7.6.5, and is
  the {\it only} known ball quotient related to a plane arrangement in
  $\fP^4$ which does {\it not} derive from those given by solutions of the
  hypergeometric differential equation.
\end{itemize}
Both of the arrangements mentioned, $\cA(W({\bf D_4}))$ and $\cA(W({\bf
  F_4}))$, arise naturally on the 45 $\fP^3$'s: the first is the
intersection with the 27 hyperplanes, the second is the intersection with
the 36 hyperplanes.
\end{itemize}

\subsection{Projection from a triple point}
Let $p\in \cI_5$ be one of the 36 triple points, and let $h_p$ be the dual
hyperplane (one of the 36 $\fP^4$'s of \ref{i1}).  The projection of
$\fP^5$ from $p$ is defined as follows. Consider the $\fP^4$ of all lines
through $p$; this is just the dual $h_p$, and each line $l_p$ through $p$
corresponds to a unique point of $h_p$ (its intersection with $h_p$). Since
any point $x$ of $\fP^5$ is on a unique line $(l_x)_p$ through $p$, the map
\begin{eqnarray}\label{e154.1} \pi_p:\fP^5 & \lra & h_p \\
  x & \mapsto & (l_x)_p\cap h_p \nonumber
\end{eqnarray} gives the {\em projection from $p$}. Restricting to $\cI_5$
this gives a generically finite (rational) map, which we also denote by
$\pi_p$, $\pi_p:\cI_5- - \ra h_p$.
\begin{lemma}\label{l154.1} $\pi_p:\cI_5- - \ra h_p$ is generically
  a double cover.
\end{lemma}
{\bf Proof:} Since the triple point has multiplicity 3, a generic line will
meet $\cI_5$ in $(5-3)=2$ further points. \ende
\begin{lemma}\label{l154.2} $\pi_p:\cI_5- - \ra h_p$ is a quotient map by the
  group $G_p\isom \fZ/2\fZ$ generated by the reflection $\gs_p$ on the root
  $p$.
\end{lemma}
{\bf Proof:} The reflection $\gs_p$ fixes $h_p$; it is the inversion
($(z_0:z_1) \mapsto (z_1:z_0)$) on any line $l_p$ through $p$, where the
homogenous coordinates are choosen such that $l_p\cap h_p=(1:1)$.  Since
$\cI_5$ is mapped by $\gs_p$ onto itself, it follows that two points of
$\cI_5\cap l_p$ are related by inversion on $l_p$. So the group action is
manifest. \ende We now describe how to make $\pi_p$ into a {\em morphism}.
First of all, one must blow up $p$; let $\grr_p:\cI_{5,p}\lra \cI_5$ denote
this blow up. Let $(\cS_3)_p$ be the copy of $\cS_3$ which is the
exceptional divisor at $p$.  For any $x\in (\cS_3)_p$, the line $(l_x)_p$
through $p$ and intersecting $h_p$ in the Segre cubic there, is tangent to
$\cI_5$ {\em at the triple point} $p$. Secondly, certain subvarieties get
{\em blown down}. Indeed, suppose $(l_x)_p$ is {\em contained in} $\cI_5$
for some $x\in \cI_5$. Then, clearly, $(l_x)_p\mapsto (l_x)_p\cap h_p$, the
whole line maps to a point, or in other words, gets blown down.
\begin{lemma}\label{l155.1} The projection $\pi_p:\cI_5- - \ra h_p$, which is
  well-defined on $\cI_{5,p}$, blows down all linear subspaces on $\cI_5$
  which pass through $p$, and is a double cover outside the union $\scL_p$
  of all such linear subspaces on $\cI_5$ passing through $p$.\ende
\end{lemma}
We now describe $\scL_p$. Recall that the linear subspaces on $\cI_5$ are
the 45 $\fP^3$'s and their intersections. Hence $\scL_p$ consists of all
the $\fP^3$'s and their intersections, which pass through $p$. Recall from
\ref{i6} that this is the set of 15 of the 45 $\fP^3$'s of Proposition
\ref{pQ7.1}. Therefore, we get
\begin{lemma}\label{l155.2} The projection $\pi_p,p:\cI_{5,p}\lra h_p$ blows
  down the union of 15 $\fP^3$'s to the 15 planes in $h_p$ which are the
  intersection of $\cS_3$ and $\cN_5$.
\end{lemma}
Now let $X=\cI_{5,p}^{\%}$, the double cover of $h_p$ branched along the
union of $\cS_3$ and $\cN_5$ (which is of degree 8, so a double cover
exists). $X$ is clearly {\em singular along} the 15 planes. Indeed:
\begin{lemma}\label{l155.3} $\pi_p,p:\cI_{5,p}\lra h_p$ factors over
  $\cI_{5,p}^{\%}$, and $\Pi:X\lra h_p$ is the double cover of $\fP^4$
  branched along the union $\cS_3\cup \cN_5$.
\end{lemma}
{\bf Proof:} This follows from the discussion above; the branch locus $\cR$
is the set:
$$\cR=\{x\in \cI_{5,p}\Big| (l_x)_p \hbox{ is tangent to $\cI_5$ at
  $x$}\}.$$ This happens if either
\begin{itemize}\item[i)] $x\in h_p$, since then $x$ is fixed by $\gs_p$;
\item[ii)] $x\in (\cS_3)_p$, the exceptional divisor over $p$.
\end{itemize}
Therefore $\cR=\cS_3\cup \cN_5$. By Lemmas \ref{l155.1} and \ref{l155.2},
15 $\fP^3$'s are blown down to $\fP^2$'s, and outside of this locus, $\Pi$
is 2:1. \ende

\subsection{Double octics and quintic hypersurfaces}
With the result of Lemma \ref{l155.3} at hand, we can get a new slant on
the quintic threefolds which are hyperplane sections of $\cI_5$. For this,
consider a hyperplane section of the cover $\Pi:X\lra h_p$, that is, let
$H\inn h_p$ be a hyperplane, and let $X_H$ be its inverse image in $X$:
$$\Pi_H:X_H\lra H,$$ a double cover of $\fP^3$. The branch locus is $H\cap
(\cS_3\cup \cN_5)$, which is the union of a cubic and a quintic surface in
$\fP^3$. Note the $H\cap \{\hbox{ one of the 15 $\fP^2$'s $\inn \cS_3\cap
  \cN_5$}\}$ is a {\em line}, contained in both $H\cap \cS_3$ and in $H\cap
\cN_5$. In other words, $H\cap (\cS_3\cup \cN_5)=S_H\cup Q_H$, where $S_H$
is the cubic surface, $Q_H$ is the quintic surface, and $$S_H\cap
Q_H=\{\hbox{15 lines}\}.$$
\begin{proposition}\label{p158.1} Let $X_H=\Pi^{-1}(H)$, the double cover of
  $\fP^3$ branched along $S_H\cup Q_H$. Then there is a canonical model
  $\-X_H$ of $X_H$ which is Calabi-Yau.
\end{proposition}
{\bf Proof:} We know the resolution of $X$; it is given by ``inverting''
the projection from the node, by blowing up along the 15 planes $\cS_3\cap
\cN_5$, yielding $\cI_{5,p}$. Let $\-X_H$ be the proper transform of $X_H$
in $\cI_{5,p}$. Assuming $H$ to be sufficiently general, $\-X_H$ clearly
has canonical singularities (as $\cI_{5,p}$ does), so we must only show
that it is Calabi-Yau. We note, however, that $\-X_H$ is (the proper
transform on $\cI_{5,p}$ of) a hyperplane section of $\cI_5$! This is
because the degree is invariant under projection, hence under $\Pi$. But
this is a hyperplane section of $\cI_5$ through the triple point $p$. Hence
the proper transform on $\-X_H$ of the exceptional divisor $(\cS_3)_p$ is a
hyperplane section of $\cS_3$, i.e., a (generically smooth) cubic surface.
This singularity is known to be canonical, and $\-X_H$ is, just as a nodal
quintic, canonically Calabi-Yau. \ende

\begin{table}
\caption{\label{table25}
  Degenerations of double octics and quintic hypersurfaces}
\begin{minipage}{16.5cm}
\begin{center}
  \fbox{\begin{minipage}{12cm}\begin{center} Space of all quintic
        hypersurfaces

        101-dimensional
\end{center}
\end{minipage}}

$$\cup$$

\fbox{\begin{minipage}{6cm}\begin{center} 120-nodal quintics

      5-dimensional
\end{center}
\end{minipage}}

$$\cup$$

\fbox{\begin{minipage}{5cm}\begin{center} quintic hypersurfaces with 111
      nodes and one multiplicity 3 singular point

      4-dimensional
\end{center}\end{minipage}}

$$\|$$

\fbox{\begin{minipage}{5cm}\begin{center} double cover $Y\lra \fP^3$,
      branched over $S\cup Q$

      $S\cap Q=\{\hbox{ 15 lines}\}$
\end{center}
\end{minipage}}

$$\cap$$

\fbox{\parbox{6cm}{double cover branched over cubic and quintic, such that
    $S\cup Q$ is stable}}

$$\cap$$ \fbox{\begin{minipage}{12cm}\begin{center} Space of all double
      octics

      149-dimensional
\end{center}\end{minipage}}
\end{center}
\end{minipage}
\end{table}

\begin{corollary}\label{c158.1} The family of hyperplane sections of $\cI_5$
  passing through one of the 36 triple points $p$ is, via projection, a
  family of Calabi-Yau threefolds which are degenerations of double octics.
\end{corollary}
It is natural to ask the meaning of this in terms of variations of Hodge
structures. Recalling that a Type III degeneration of a K3 surface,
corresponding to a zero-dimensional boundary component of the period
domain, is one like a quartic degenerating into four planes, it is natural
to ask
\begin{question} Is a double cover of $\fP^3$ branched over the union of a
  cubic and a quintic a semistable degeneration of a double octic?
\end{question}
\begin{remark} There is a notion of ``connecting''
  moduli spaces of CY threefolds by degenerations, and the Corollary
  \ref{c158.1} shows that the moduli space of quintic hypersurfaces in
  $\fP^4$ and the moduli space of double octics are connected; the
  birational transformations which are required for such ``connections''
  are given here by projection in projective space, very geometric.
\end{remark}
In Table \ref{table25} we give a rough description of these relations.

\subsection{The dual picture}
Now consider $\ga(\cS_3\cup \cN_5)$, with $\ga$ the map (\ref{e126.1})
given by the quadrics on the ten nodes of $\cS_3$. By Theorem \ref{t122a.1}
and by definition of $\cW_{10}$ (\ref{e127.1}), we have
\begin{equation}\label{e156.1} \ga(\cS_3\cup \cN_5) = \cI_4\cup \cW_{10},
\end{equation}
and by Theorem \ref{t127.1}, the intersection $\cI_4\cap \cW_{10}$ consist
of 10 quadric surfaces. Define $\cW$ to be the double cover of $\fP^4$
branched along $\cW_{10}$:
\begin{equation}\label{e156.2} \tau:\cW\lra \fP^4=(h_p)^{\vee}.
\end{equation}
We may consider the fibre square:
\begin{equation}\label{e156.3} \begin{array}{rcl} {\cal Z} & \lra & {\cal Y} \\
    \downarrow & & \downarrow \pi \\ \tau:{\cal W} & \lra & {\Bbb P}^4
  \end{array}
\end{equation}
where $\pi:\cY\lra \fP^4$ is defined in Definition \ref{d133.1}. Then
$\pi_{\cZ}:\cZ\lra \fP^4$ is a Galois cover with Galois group $G_{\cZ}\isom
\fZ/2\fZ\times \fZ/2\fZ$. Let $H\isom \fZ/2\fZ\inn G_{\cZ}$ be the diagonal
subgroup; it is a normal subgroup, and we may form the quotient
$$\eta:\cZ\lra \cZ',\quad \cZ'=\cZ/H.$$
\begin{lemma}\label{l156.1} $\pi_{\cZ}:\cZ\lra \fP^4$ factors over $\eta$,
  and $\eta':\cZ'\lra \fP^4$ is a double cover, hence Galois.
\end{lemma}
{\bf Proof:} This is a general fact about fibre squares of double covers
like (\ref{e156.3}). \ende
\begin{theorem}\label{t156.1} The rational map $\ga$ induces a rational map
  of the double covers $\Pi:X\lra \fP^4$ of Lemma \ref{l155.3} and
  $\eta':\cZ'\lra \fP^4$ of Lemma \ref{l156.1}. Furthermore, the rational
  map
  $$\Xi:X- - \ra \cZ'$$ is $\gS_6$-equivariant.
\end{theorem}
{\bf Proof:} Recall from Lemma \ref{l25aux} that $\ga$ blows up the ten
nodes and blows down the tangent cones of the nodes to the quadric surfaces
(on $\cI_4$). $\Pi:X\lra h_p$ is a double cover branched along $\cS_3\cup
\cN_5$, and we can calculate the image of the branch locus under $\ga$.
The ten nodes get blown up, the ten quadric cones (in $\fP^4$) get blown
down to quadric surfaces (in the exceptional $\fP^3$'s). Let $\~C\inn X$ be
the inverse image in $X$ of the union of the ten quadric cones; then on
$X\bs \~C$, $\ga$ is {\it biregular}. On the other hand, $\ga(\~C)$ is just
the union of the ten quadric surfaces of the intersection $\cI_4\cap
\cW_{10}$. Consequently
$$\ga_{|X\bs \~C}:X\bs\~C \lra \cZ'\bs(\eta')^{-1}(\cI_4\cap \cW_{10})$$ is
a reguar morphism of double covers, and letting $C\inn h_p$ denote the ten
quadric cones, $X\bs \~C\lra h_p\bs C$ is a double cover, as is also
$$\cZ'\bs(\eta')^{-1}(\cI_4\cap \cW_{10}) \lra (h_p)^{\vee}\bs \cI_4\cap
\cW_{10}, $$ while $\ga(C)=\cI_4\cap \cW_{10}$. Hence in the diagram
$$\begin{array}{rcl} X & \stackrel{\Xi}{\lra} & \cZ' \\ \downarrow & &
  \downarrow \\ h_p & \stackrel{\ga}{\lra} & (h_p)^{\vee} \end{array}$$
$\Xi$ is regular outside of $\~C$ and maps $\~C$ to $(\eta')^{-1}(\cI_4\cap
\cW_{10})$. Furthermore, everything is defined $\gS_6$-equivariantly. This
proves the Theorem. \ende

\begin{corollary}\label{c157.1} $\cI_5$ sits $\gS_6$-equivariantly
  birationally in the center of the diagram \unitlength1cm
  $$\begin{picture}(2,2)\put(0,0){${\cal W}$}
    \put(.5,0.15){\vector(1,0){1.3}} \put(1.9,0){${\Bbb P}^4$.}
    \put(0,1.6){${\cal Z}$} \put(.5,1.75){\vector(1,0){1.3}}
    \put(1.9,1.6){${\cal Y}$}
    \put(.25,1.5){\vector(0,-1){1}}\put(1.95,1.5){\vector(0,-1){1}}
    \put(.4,1.6){\vector(1,-1){.5}} \put(.9,.9){${\cal Z}'$}
    \put(1.2,.8){\vector(1,-1){.5}}
\end{picture}$$
This shows the relation between the quintic $\cI_5$ and the Coble variety
$\cY$. \end{corollary} {\bf Proof:} We have the series of modifications
$$\begin{array}{ccccccc} \cI_5 & \lra & \cI_{5,p} & \stackrel{\gb}{\lra} &
  X & \stackrel{\Xi}{\lra} & \cZ' \\ & & & & \downarrow & & \downarrow \\ &
  & & & \fP^4 & \stackrel{\ga}{\lra} & \fP^4, \end{array}$$ where
$\cI_5\lra \cI_{5,p}$ blows up the node $p$, $\gb$ blows down the 15
$\fP^3$'s through the node to the 15 $\fP^2$'s of the intersection
$\cS_3\cap \cN_5$, $\ga$ and $\Xi$ are as described above. Since $\cZ'$
clearly sits in the center of the diagram and all modifications are
$\gS_6$-equivariant, the Corollary follows. \ende

\section{${\cal I}_5$ and cubic surfaces}

\subsection{The Picard group}
Let $A_1(\cI_5)$ be the Chow group of Weil divisors modulo algebraic
equivalence. Clearly a generic hyperplane section yields an element in
$A_1(\cI_5)$, which we denote by $n$. Recall the reducible hyperplane
sections of Proposition \ref{pQ7.1} which split each into the union of five
copies of $\fP^3$. These subvarieties are divisors on $\cI_5$, hence also
yield classes in the Chow group. These 45 divisors are related by 27
relations, the sum of the five classes in the Chow group being equivalent
to $n$. Since $\cI_5$ is normal, we have an injection $\Pic(\cI_5)\hra
A_1(\cI_5)$. Let $\Pic^0(\cI_5)$ denote the orthogonal complement of the
class $n$ in $A_1(\cI_5)$ with respect to this injection. Then we have
\begin{lemma}\label{liq6.1} We have an exact sequence of $\fZ$-modules,
  $$0\lra \fZ^6\lra \fZ^{27}\lra \fZ^{45} \lra {\em\Pic}^0(\cI_5) \lra 0.$$
\end{lemma}
{\bf Proof:} The 45 $\fP^3$'s are classes in $A_1(\cI_5)$ which generate
$\Pic^0(\cI_5)$ (as they contain all singularities), and the 27 relations
are those just mentioned, given by the 27 hyperplane sections. So the
sequence is clear as soon as we have shown that the rank of $\Pic^0(\cI_5)$
is 24 (see the sequence (\ref{e48a.1})). This now follows from the
Lefschetz hyperplane theorem, as the dimension of $\cI_5$ is four, so there
is an isomorphism between the $H^2$'s of $\cI_5$ and a hyperplane section.
We may apply the Lefschetz theorem because the singularities of $\cI_5$ and
of a hyperplane section are local complete intersections (see the book by
Goresky \& MacPherson for details). \ende Note that this sequence displays
$\Pic^0(\cI_5)$ as an {\em irreducible} $W(E_6)$-module. Furthermore, we
see that just as in (\ref{e48a.1}), this sequence is dual to the
corresponding sequence for cubic surfaces.

\subsection{${\cal I}_5$ and cubic surfaces: combinatorics}
We collect the facts relating the combinatorics of the 27 lines with those of
$\cI_5$ in Table \ref{table26}.
\begin{table}[htb]
\caption{\label{table26} Combinatorics of ${\cal I}_5$ and the 27 lines}
\vspace*{.5cm}
\renewcommand{\arraystretch}{1.5}
\begin{tabular}{|l|l|}\hline
Locus on a cubic surface (see Table \ref{table20}) &
Locus on ${\cal I}_5$ \\ \hline \hline
27 lines $a_i, b_i, c_{ij}$ & 27 hyperplane sections $\{a_i=0\}\cap {\cal
I}_5$, etc. \\ \hline
2 lines are skew & \parbox{7cm}{the hyperplanes intersect in one of 216
$\fP^3$'s dual to the lines of \ref{i4}; this $\fP^3$ intersects ${\cal I}_5$
in the union of three planes and a quadric (see Lemma \ref{l126.1})} \\ \hline
\parbox{7cm}{two lines are in a tritangent} & \parbox{7cm}{the
hyperplanes intersect in one of the 45 $\fP^3$'s of \ref{i8} } \\ \hline
45 tritangents & the 45 $\fP^3$'s of \ref{i8} \\ \hline
\parbox{7cm}{Two tritangents meet in a line of the cubic surface} &
\parbox{7cm}{two of the 45 $\fP^3$'s meet in a $\fP^2$; this is one of the
planes in the $\fP^3$ defining the arrangement ${\cal A}(W({\bf D_4}))$ as
discussed in \ref{i8}} \\ \hline
\parbox{7cm}{Two tritangents meet in a line outside of the cubic surface} &
\parbox{7cm}{two of the 45 $\fP^3$'s are {\em skew}, i.e., meet only in a
line; this line is part of the singular locus of the arrangement ${\cal
A}(W({\bf D_4}))$ just mentioned} \\ \hline
36 double sixes & \parbox{7cm}{36 triple points of ${\cal I}_5$ AND 36 copies
of the Nieto quintic ${\cal N}_5$} \\ \hline
\parbox{7cm}{Two double sixes are azygetic} & \parbox{7cm}{two of the triple
points lie on one of the 120 lines of the singular locus of ${\cal I}_5$} \\
\hline
\parbox{7cm}{Two double sixes are syzygetic} & \parbox{7cm}{two of the triple
points do not lie on one of the 120 lines} \\ \hline
\parbox{7cm}{A line is {\em not} contained in a double six} &
\parbox{7cm}{the hyperplane dual to the line contains the triple point which
corresponds to the double six} \\ \hline \end{tabular}

\end{table}

\section{The dual variety}
Let $\cI_5^{\vee}$ be the projective dual variety to $\cI_5$; since
$\cI_5$ is invariant under $W(E_6)$, so is $\cI_5^{\vee}$. Although we
do not have explicit equations for $\cI_5^{\vee}$, we can say quite a
bit about its geometry, just from the fact that it is dual to $\cI_5$.

\subsection{Degree}

First we show that {\em degree of\ $\cI_5^{\vee}$=10m+4k}. Quite
generally, one can say the following. Suppose we are given a variety
$X\inn \fP^n$ which has singular locus consisting of a set of {\em lines},
meeting each other in a set of {\em points}, and let us further assume
that the situation is symmetric, i.e., each line contains the same
number of points, each point being hit by the same number of lines;
let us denote these numbers by $N=\#$ lines, $M=\#$ points, $\gn=\#$
points in each line and $\gm=\#$ lines meeting at each point. Consider
the dual variety $X^{\vee}$.  We claim:
\begin{itemize}
\item[-] There are $N$ $\fP^{n-2}$'s $\inn X^{\vee}$.
\item[-] Each $\fP^{n-2}$ is cut out by $\gn$ hyperplanes.
\item[-] There are $M$ such special hyperplane sections of $X^{\vee}$.
\item[-] Each of the $M$ hyperplanes meets $X^{\vee}$ in $\gm$ of the
$N\ \fP^{n-2}$'s.
\item[-] Hence, deg($X^{\vee}$)=$m\gm+rest$,
\end{itemize}
where the $rest$ is given in terms of the local geometry around the given
point.
The proofs of these are immediate: each of the points corresponds to a
hyperplane (=set of all hyperplanes through the point), each line defines
dually a $\fP^{n-2}$, and since $X$ is singular along the line, each hyperplane
through the line is {\em tangent} to $X$ there $\Rightarrow$ the dual
$\fP^{n-2}\inn X^{\vee}$.  The other statements are then clear.
To determine $rest$, consider the following. The set theoretic
image of the given point in the dual variety is the {\em total} transform
({\em not} the proper transform) of the given point. This is set
theoretically easy to compute, but there may be a multiplicity coming in.

We apply these considerations to $\cI_5$ and $\cI_5^{\vee}$: on $\cI_5$
we have singular lines, $N$=120, $M$=36, $\gn$=3, $\gm$=10, and hence
$deg(\cI_5^{\vee})=m10+rest.$ In our case $rest$ is easy to figure out:
recall that we resolved the singularities of \cI$_5$ by blowing up the 36
points, then the 120 lines; the resolving divisors over the points were copies
of the Segre cubic. The variety dual to the Segre cubic is the Igusa quartic,
and the image of the ten nodes on the Segre cubic are ten quadric
surfaces (\ref{e117a.1}) which are {\em tangent hyperplane sections},
i.e., the hyperplanes which meet the Igusa quartic in one such quadric
and are tangent to it there. These ten hyperplanes are of course just the
10 $\fP^3$'s on the dual variety being cut out by the chosen hyperplane
section (see Proposition \ref{p160.1} below).
This hyperplane section of $\cI_5^{\vee}$ may be {\em tangent}
to $\cI_5^{\vee}$ along the Igusa quartic, hence
$$deg(\cI_5^{\vee})=10m+k\cdot4.$$

\subsection{Singular locus}

Consider the 45 $\fP^3$'s on \cIf; since there is a pencil of hyperplanes
through each, the dual variety $\cI_5^{\vee}$\ will have 45 singular
lines, which meet in 27 points (which are dual to the 27 hyperplanes
cutting out the 45 $\fP^3$'s).
These 27 points are of course $A_i, B_i, C_{ij}$.
Applying our reasoning from above to this we see that
deg(${\cI_5^{\vee}}^{\vee}$)=5+ rest. We conclude rest=0, or in other words,
{\em a resolution of singularities of ${\cI_5^{\vee}}$ is affected by blowing
up the 45 lines simultaneously; there is no exceptional divisor over the 27
points.}

However, since we are dealing with fourfolds, $\cI_5^{\vee}$ could even be
normal and still have a singular locus of dimension two. For example, it is
reasonable to believe that the ten quadrics on each copy of the Igusa quartic
$\cI_4$ on the reducible hyperplane sections discussed below might be {\em
singular} on $\cI_5^{\vee}$, but that is of course just a guess.
Furthermore, there is no reason whatsoever why the dual variety should be
normal. In fact, it is a case of great exception when the dual variety is
normal, the general case being that there is a singular parabolic divisor
(coming from the intersection $\Hess(X)\cap X$), as well as a double point
locus, also (in general) a divisor, coming from the set of bitangents.
In our case, however, since $\Hess(\cI_5)\cap \cI_5$ consists of the union of
the 45 $\fP^3$'s, all of which get blown down, there is no parabolic {\em
divisor}. But there is no easy way to exclude a double point divisor.

\subsection{Reducible hyperplane sections}

As already mentioned, $\cI_5^{\vee}$\ contains 120 $\fP^3$'s, each being cut
out by three of the $h$'s, (in fact by a triple of azygetic double sixes), and
each
such intersection $h\cap \cI_5^{\vee}$ consists of ten such $\fP^3$'s, plus a
copy of the Igusa quartic. There are 36 such hyperplane sections which
decompose
into ten $\fP^3$'s and a copy of the Igusa quartic:
\begin{proposition}\label{p160.1} The 36 hyperplane sections ${\bf h}\cap
\cI_5$, ${\bf h}=h,\ h_{ij},\ h_{ijk}$, are reducible, consisting of ten
$\fP^3$'s and a copy of the Igusa quartic $\cI_4$. The ten $\fP^3$'s are just
the $K_{ijk}$ of (\ref{e117a.1a}), each a bitangent plane to $\cI_4$.
\end{proposition}
{\bf Proof:} These are the 36 hyperplanes dual to the 36 triple points of
$\cI_5$; at each such $p$ ten of the 120 lines meet, and the triple point
itself yields the copy of $\cI_4$ (it is blown up with exceptional divisor
$(\cS_3)_p$, which is dual to $(\cI_4)_p$, a copy of $\cI_4$). \ende
So restricted to the triple point, the duality $\cI_5- - \ra \cI_5^{\vee}$
yields precisely the dual map $\ga$ of (\ref{e126.1})!

The 120 $\fP^3$'s meet two at a
time in 270 $\fP^2$'s, each of which is cut out by six of the $h$'s (2 triples
of
azygetic double sixes, two rows in a triple).  Note that these 270 $\fP^2$'s
are
the $t_6(2)$ of Table \ref{table23}.
Each $\fP^2$ contains two nodes
 and five of the 27 points, as well as two of the 45 lines. Through each
such line two of these $\fP^2$'s pass (as each line is cut out by 12 of the
$h$'s).  Therefore in each $h$ we have ten $\fP^3$'s meeting in ${10 \choose
2}=45\ \fP^2$'s which meet in 15 of the 45 $\fP^1$'s, and contain 15 of the 27
points. The 15 lines and 15 points are just the singular locus of the Igusa
quartic, and the ten $\fP^3$'s are tangent to the Igusa quartic along quadrics,
as mentioned earlier.

\subsection{Special hyperplane sections}

Inspection of the 27 forms and 27 points in $\fP^5$ shows that each of the 27
hyperplanes contains {\em none} of the 27 points and {\em none} of the 45
lines;
it follows that hyperplane sections such as $\cK:=\cI_5^{\vee}\cap \{a_1=0\}$
are irreducible hypersurfaces in $\fP^4$ with 45 isolated singularities,
coming from the intersections with the singular lines of $\cI_5^{\vee}$.
As mentioned above, there may also be a singular locus coming from other
singularities on $\cI_5^{\vee}$.
Furthermore, there are 40 $\fP^2$'s lying on this threefold, and 16 hyperplanes
in
$a_i$ which cut out ten of these on $\cK$.  The 16 hyperplanes are those 16 of
the 216 $\fP^3$'s which lie in $a_i$, corresponding to the 16 lines which $a_i$
is skew to.  The symmetry group of this threefold is $W(D_5)$. This is a {\em
degeneration} of a {\em generic} hyperplane section, which will contain 120
$\fP^2$'s.

\end{document}